\documentclass[12pt,leqno]{article}
\usepackage{amssymb}
\usepackage[mathscr]{eucal}
\usepackage{amsmath,amssymb,latexsym,theorem,bbm}
\usepackage{amsmath,amssymb,latexsym,theorem}
\usepackage{color,url}
\usepackage{appendix}

\newcommand{\CC}{\mathbb{C}}

\newcommand{\NN}{\mathbb{N}}

\newcommand{\RR}{\mathbb{R}}

\newcommand{\ZZ}{\mathbb{Z}}

\newcommand{\ba}{{\boldsymbol{a}}}
\newcommand{\bA}{{\boldsymbol{A}}}

\newcommand{\bB}{{\boldsymbol{B}}}

\newcommand{\tb}{\widetilde{b}}

\newcommand{\tbB}{{\widetilde{\bB}}}

\newcommand{\bc}{{\boldsymbol{c}}}

\newcommand{\bD}{{\boldsymbol{D}}}

\newcommand{\be}{{\boldsymbol{e}}}

\newcommand{\bI}{{\boldsymbol{I}}}

\newcommand{\bP}{{\boldsymbol{P}}}

\newcommand{\br}{{\boldsymbol{r}}}

\newcommand{\bS}{{\boldsymbol{S}}}

\newcommand{\bT}{{\boldsymbol{T}}}

\newcommand{\bv}{{\boldsymbol{v}}}

\newcommand{\tbv}{\widetilde{\bv}}
\newcommand{\bx}{{\boldsymbol{x}}}
\newcommand{\bX}{{\boldsymbol{X}}}
\newcommand{\by}{{\boldsymbol{y}}}
\newcommand{\bY}{{\boldsymbol{Y}}}

\newcommand{\bz}{{\boldsymbol{z}}}
\newcommand{\bZ}{{\boldsymbol{Z}}}

\newcommand{\Bbeta}{{\boldsymbol{\beta}}}

\newcommand{\tBbeta}{\widetilde{\Bbeta}}

\newcommand{\blambda}{{\boldsymbol{\lambda}}}
\newcommand{\tlambda}{\widetilde{\lambda}}
\newcommand{\tblambda}{\widetilde{\blambda}}

\newcommand{\bxi}{{\boldsymbol{\xi}}}

\newcommand{\bmu}{{\boldsymbol{\mu}}}

\newcommand{\bbeta}{{\boldsymbol{\beta}}}

\newcommand{\tbD}{ \boldsymbol{\widetilde{\bD}} }

\newcommand{\bphi}{{\boldsymbol{\phi}}}
\newcommand{\bvarphi}{{\boldsymbol{\varphi}}}
\newcommand{\bzero}{{\boldsymbol{0}}}
\newcommand{\bone}{{\boldsymbol{1}}}

\newcommand{\cB}{{\mathcal B}}

\newcommand{\cF}{{\mathcal F}}

\newcommand{\cR}{{\mathcal R}}

\newcommand{\cU}{{\mathcal U}}

\newcommand{\dd}{\mathrm{d}}
\newcommand{\ee}{\mathrm{e}}

\newcommand{\EE}{\operatorname{\mathbb{E}}}

\newcommand{\PP}{\operatorname{\mathbb{P}}}

\newcommand{\card}{\operatorname{\mathrm{card}}}

\newcommand{\tN}{\widetilde{N}}

\newcommand{\tv}{\widetilde{v}}

\newcommand{\tW}{\widetilde{W}}

\newcommand{\vare}{\varepsilon}

\renewcommand{\mid}{\,|\,}

\renewcommand{\leq}{\leqslant}
\renewcommand{\geq}{\geqslant}

\newcommand{\bbone}{\mathbbm{1}}

\newcommand{\proofend}{\hfill\mbox{$\Box$}}

\numberwithin{equation}{section}

\theoremstyle{change} \theorembodyfont{\em}
\newtheorem{Lem}{Lemma.}[section]
\newtheorem{Thm}[Lem]{Theorem.}
\newtheorem{Pro}[Lem]{Proposition.}
\newtheorem{Cor}[Lem]{Corollary.}
\newtheorem{Def}[Lem]{Definition.}

\theorembodyfont{\rm}
\newtheorem{Rem}[Lem]{Remark.}

\begin{document}

\begin{center}
 {\bfseries\Large Distributional properties of jumps of multi-type CBI processes}

\vspace*{3mm}

 {\sc\large
 M\'aty\'as $\text{Barczy}^{*,\diamond}$,
  \ Sandra $\text{Palau}^{**}$}

\end{center}

\vskip0.2cm

\noindent
 * HUN-REN--SZTE Analysis and Applications Research Group,
   Bolyai Institute, University of Szeged,
   Aradi v\'ertan\'uk tere 1, H--6720 Szeged, Hungary.

\noindent
 ** Instituto de Investigaciones en Matem\'aticas Aplicadas y en Sistemas,
    Universidad Nacional Aut\'onoma de M\'exico, CDMX, 04510, Ciudad de M\'exico, M\'exico.

\noindent E-mails: barczy@math.u-szeged.hu (M. Barczy),
                  sandra@sigma.iimas.unam.mx (S. Palau).

\noindent $\diamond$ Corresponding author.

\vskip0.2cm


{\renewcommand{\thefootnote}{}
\footnote{\textit{2020 Mathematics Subject Classifications\/}:
 60J80, 60G55.}
\footnote{\textit{Key words and phrases\/}:
 multi-type continuous state and continuous time branching processes with immigration, first jump time,
 supremum for jumps, total L\'evy measure.}
\footnote{M\'aty\'as Barczy was supported by the project TKP2021-NVA-09.
Project no.\ TKP2021-NVA-09 has been implemented with the support
 provided by the Ministry of Culture and Innovation of Hungary from the National Research, Development and Innovation Fund,
 financed under the TKP2021-NVA funding scheme.}
}

\vspace*{-10mm}

\begin{abstract}
We study the distributional properties of jumps of multi-type continuous state and continuous
 time branching processes with immigration (multi-type CBI processes).
We derive an expression for the distribution function of the first jump time of a multi-type CBI process
 with jump size in a given Borel set having finite total L\'evy measure, which is defined
 as the sum of the measures appearing in the branching and immigration mechanisms of the multi-type CBI process in question.
Using this we derive an expression for the distribution function of the local supremum of the norm of the jumps
 of a multi-type CBI process.
Further, we show that if $A$ is a nondegenerate rectangle anchored at zero and with total L\'evy measure zero,
 then the probability that the local coordinate-wise supremum of jumps of the multi-type CBI process belongs to $A$ is zero.
We also prove that a converse statement holds.
\end{abstract}

\tableofcontents

\section{Introduction}
\label{section_intro}

Multi-type continuous state and continuous time branching processes with immigration (multi-type CBI processes)
 are special Markov processes with values in $[0,\infty)^d$  that can be well-applied for describing the random evolution
 of a population consisting of individuals having different (finitely many) types.
In this paper, we study the distributional properties of jumps of such processes.
Given a multi-type CBI process $(\bX_u)_{u\geq 0}$, a Borel set $A$ of $[0,\infty)^d\setminus\{\bzero\}$ and $t>0$,
 let us introduce $\tau_A:= \inf\{u > 0: \Delta \bX_u \in A\}$ and $J_t(A):= \card(\{u \in (0, t] : \Delta \bX_u  \in A\})$,
 with the convention $\inf(\emptyset) := \infty$, where $\Delta \bX_u := \bX_u - \bX_{u-}$, $u >0$,
 and $\card(H)$ denotes the cardinality of a set $H$.
One can call $\Delta \bX_u\in[0,\infty)^d$ and its Euclidean norm, the size vector and the size of the jump of $\bX$ at time $u>0$,
 respectively.

In case of a single-type CBI process $(X_u)_{u\geq 0}$ (i.e., $d=1$),
 He and Li \cite[Theorem 3.1]{HeLi} studied the distribution of $J_t(A)$ and $\tau_A$, namely,
 they showed that $\PP(J_t(A)<\infty)=1$, $t>0$, and derived an expression for
 $\PP(\tau_A>t \mid X_0=x)$, $t\geq 0$, $x\geq 0$, in terms of a solution of
 a deterministic differential equation provided that the set $A$ has a finite total L\'evy measure, which is defined
 as the sum of the measures appearing in the branching and immigration mechanisms of $(X_u)_{u\geq 0}$.
He and Li \cite[Theorem 4.1]{HeLi} also proved that the total L\'evy measure in question is equivalent to
 the probability measure  $\PP\big( \sup_{s\in(0,t]} \Delta X_s\in \cdot \mid X_0=x\big)$
 for any $t>0$ and $x+\gamma>0$, where $\gamma$ is the drift in the immigration
 mechanism of $(X_t)_{t\geq 0}$.
The above recalled results of He and Li \cite{HeLi} on the distributional properties of jumps of single-type CBI processes
 are also contained in Sections 10.3 and 10.4 in the new second edition of Li's book \cite{Li}.

In this paper, we generalize some results of He and Li \cite{HeLi} on the distributional properties of jumps of single-type CBI processes.
Before presenting our new results, we recall and summarize some related research on the distributional properties of jumps of CBI processes.
Jiao et al.\ \cite{JiaMaSco} considered a so-called $\alpha$-CIR process, which is a natural extension of
 the standard Cox-Ingersoll-Ross (CIR) process
 by adding a jump part driven by an $\alpha$-stable L\'evy process with index $\alpha \in (1,2]$.
This is a new short interest rate model.
They derived expressions for the Laplace transform of the number of large jumps in a finite time interval
 whose jump sizes are larger than a given value (see Jiao et al.\ \cite[Proposition 5.1]{JiaMaSco}),
 and for the distribution function and expectation of the first large jump time
 (see Jiao et al.\ \cite[Corollary 5.2 and Proposition 5.4]{JiaMaSco}).
Ji and Li \cite[Proposition 4.5]{JiLi} considered the expectation $\EE(f(X_t)\bone_{\{ \zeta_k > t\}})$, $t\geq 0$,
 where $k$ is a positive integer, $(X_t)_{t\geq 0}$ is a single-type CBI process starting at $0$,
 the random variable $\zeta_k$ is the $k$-th jump time of $(X_t)_{t\geq 0}$ with jump size in $(1,\infty)$
 and $f:[0,\infty)\to[0,\infty)$ is a convex and nondecreasing function satisfying some further properties
 (including integrability ones).
They derived an upper estimation for $\EE(f(X_t)\bone_{\{ \zeta_k > t\}})$ in terms of
 $\EE(f(X_s)\bone_{\{ \zeta_{k-1} > t-s \}})$, $s\in[0,t]$, and the distribution function of $\zeta_1$.
Chen and Li \cite[(1.7) and (1.8)]{CheLi} introduced the notion of a continuous time mixed state branching process,
 which is a two-dimensional branching Markov process $(Y^{(1)}_t,Y^{(2)}_t)_{t\geq 0}$ taking values
 in $[0,\infty)\times \{0,1,2,\ldots\}$ given as a pathwise unique strong solution of
 a two-dimensional SDE.
Among others, the authors derived an expression for the distribution function of
 $\inf\{s\geq 0 : \Delta Y^{(1)}_s > r_1 \;\;\text{or}\;\;  \Delta Y^{(2)}_s > r_2\}$
 in terms of the solution a two-dimensional system of deterministic differential equations,
 where $r_1,r_2\geq 0$, see Chen and Li \cite[Theorem 5.3]{CheLi}.
Horst and Xu \cite[Subsection 5.1]{HorXu} used Theorem 3.1 in He and Li \cite{HeLi} and the arguments in its proof
 to describe the joint distribution of the arrival time of the last jump whose magnitude belongs to a given set
 and the number of jumps with magnitudes in the same set up to a given time point
 for a two-dimensional continuous time stochastic process of which the second coordinate process is a single-type
 continuous state and continuous time branching process (single-type CB process).
The two-dimensional continuous time stochastic process in question plays an important role in financial mathematics,
 and its second coordinate process is called a volatility process.

According to our knowledge, multi-type CB processes were first introduced in 1969 by Watanabe \cite{Wa}
 as special measure-valued branching processes,
 for more details on the history, see Li \cite{Li}.
There are different characterizations of multi-type CBI processes, for example, using their Laplace transform (see Duffie et al.\ \cite{DufFilSch}),
 or viewing the process as the pathwise unique strong solution of an SDE (see Ma \cite{Ma} and Barczy et al.\ \cite{BarLiPap2}),
 or as the unique solution to a time-change equation (see Caballero et al.\ \cite{CabPerUri}).
Through these approaches, many authors have addressed several questions for multi-type CBI processes, among others,
 their asymptotic behaviour \cite{BarPap,KypPalRen,BarPalPap,BarPalPap2},
 extinction times \cite{KypPal,ChaMar2}, existence of densities \cite{FriJinRud2}, and behavior at the boundary \cite{FriJinRud}.

Now we turn to summarize the content of our paper, it is structured as follows.
At the end of the Introduction, we collect all the notations that will be used throughout the paper.
In Section \ref{section_CBI}, we recall the definition of multi-type CBI processes, a result on their representation
 as pathwise unique strong solutions of some SDEs (see \eqref{SDE_atirasa_dimd}),
 and the notion of irreducibility.

Section \ref{Section_CBI_intCBI} is devoted to study the integral of a multi-type CBI process.
We show that a $(2d)$-dimensional stochastic process having coordinate processes a $d$-type CBI process
 and its integral process is a $(2d)$-type CBI process (where $d$ is a positive integer),
 see Proposition \ref{Pro_2dCBI}.
Using this result, we derive a formula for the Laplace transform of the integral of a multi-type CBI process,
 see Proposition \ref{Pro_int_CBI_Laplace}.
We emphasize that these two results are not new, only our proof technique is new (for more discussion,
 see Remark \ref{Rem_discussion} and the paragraph after it).
We also study the limit behaviour of a function which appears in the formula for the Laplace transform
 of the integral of a multi-type CBI process, this function is nothing else but a solution of a deterministic
 differential equation, see Proposition \ref{Pro_tv}.
This result can be considered as a multi-type counterpart of Proposition 2.2 in He and Li \cite{HeLi}.

In Section \ref{Section_dist_jump}, we investigate the distributional properties of jump times of multi-type CBI processes.
Given a multi-type CBI process $(\bX_u)_{u\geq 0}$, under some moment conditions, for
 any Borel set $A$ in $[0,\infty)^d\setminus\{\bzero\}$ having finite total L\'evy measure,
 we show that $\EE(J_t(A))<\infty$, $t>0$, and we derive an expression for
 $\PP(\tau_A>t \mid \bX_0=\bx)$, $t\geq 0$, $\bx\in\RR_+^d$, in terms of a solution of
 a deterministic differential equation, see Theorem \ref{Thm_jump_times_dCBI}.
Our result is a generalization of Proposition 3.1 and Theorem 3.1 in He and Li \cite{HeLi} to multi-type CBI processes.
We mention that part (iii) of our Theorem \ref{Thm_jump_times_dCBI} has just appeared as Example 12.2 in
 the new second edition of Li's book \cite{Li}.
The two research works have been carried out parallelly, so we decided to present our result and its proof as well.
For a more detailed discussion and a comparison of the two proofs, see the second part of the paragraph
 before Theorem \ref{Thm_jump_times_dCBI}.
In Corollary \ref{Cor_tau}, among others, we derive some sufficient conditions under which
 $\PP(\tau_A = \infty \mid \bX_0 = \bx)=1$ or
 $\PP(\tau_A<\infty \mid \bX_0 = \bx)=1$ holds, respectively.
We also present some conditions under which the probability $\PP(\tau_A = \infty \mid \bX_0 = \bx)$
 can be expressed in terms of the inverse of the branching mechanism of $(\bX_u)_{u\geq 0}$,
 introduced by Chaumont and Marolleau \cite[Theorem 2.1]{ChaMar2}.

In Section \ref{Section_max_jumps}, first we deal with the local and global supremum of the sizes of the jumps of multi-type CBI processes.
We derive expressions for the probability
  $\PP\big(\sup_{s\in(0,t]} \Vert \Delta \bX_s\Vert \leq r \mid \bX_0 = \bx\big)$,
 $r>0$, $t>0$, $\bx\in\RR_+^d$, and, if, in addition, the total L\'evy measure of $\RR_+^d\setminus\{\bzero\}$
 is finite, for $\PP\big(\sup_{s\in(0,t]} \Vert \Delta \bX_s\Vert = 0 \mid \bX_0 = \bx\big)$, $t>0$, $\bx\in\RR_+^d$, as well,
 see Proposition \ref{Pro_max_jump}.
This result is a generalization of Theorem 4.1 in He and Li \cite{HeLi} to multi-type CBI processes.
In Proposition \ref{Pro_sup_jump}, we derive sufficient conditions
 under which $\sup_{s\in(0,\infty)} \Vert \Delta\bX_s\Vert$ is a constant with probability one
 given any initial value $\bx\in[0,\infty)^d$.
This result is a multi-type counterpart of Corollary 4.1 in He and Li \cite{HeLi}.

In order to formulate the main result of the paper,
 for all $t>0$ and $\bx\in[0,\infty)^d$, let us introduce the probability measure $\pi_{t,\bx}$ on $([0,\infty)^d,\cB([0,\infty)^d))$ defined by
 \[
  \pi_{t,\bx}(A):=\PP\Bigg(  \sup_{s\in(0,t]} \Delta \bX_s\in A  \;\Big\vert\; \bX_0=\bx\Bigg), \qquad A\in\cB([0,\infty)^d),
 \]
 where $\sup_{s\in(0,t]} \Delta \bX_s$ denotes coordinate-wise supremum.
Further, let
 \[
 \cR_d:=\left\{ \left(\prod_{i=1}^d [0, w_i] \right) \setminus \{\bzero\}:  w_1,\ldots,w_d>0 \right\}
 \]
 be the collection of nondegenerate rectangles in $[0,\infty)^d$ anchored at $\bzero$, but not containing $\bzero$.
Under some moment conditions, we show that if $A\in\cR_d$ is such that its total L\'evy measure is zero,
 then $\pi_{t,\bx}(A)=0$ for all $t>0$ and $\bx\in[0,\infty)^d$.
Conversely, if $\pi_{t,\bx}(A)=0$ for some $t>0$, $\bx\in(0,\infty)^d$ and $A\in\cR_d$,
 then the total L\'evy measure of $A$ is zero.
Furthermore, in case that the immigration mechanism of $(\bX_u)_{u\geq 0}$
 has a drift $\Bbeta$ with strictly positive coordinates, then we can extend it to $\bx\in[0,\infty)^d$ as well, see Theorem \ref{Thm_sup_jump_equiv}.
This result can be considered as a multi-type counterpart of Theorem 4.2 in He and Li \cite{HeLi}, which is about single-type CBI processes.
In Remark \ref{Rem_HeLI_comparison},
 we point out the fact that in case of $d\geq 2$, part (i) of Theorem \ref{Thm_sup_jump_equiv} does not hold
 for a general Borel set $A$ in $[0,\infty)^d\setminus\{\bzero\}$,
 which also shows that, in general, the total L\'evy measure of $(\bX_t)_{t\in\RR_+}$
 is not equivalent to the probability measure $\pi_{t,\bx}$, where $\bx\in\RR_+^d$ and $t>0$.
In Remark \ref{Rem_nulla}, we highlight why the case $\bx=\Bbeta=\bzero$ is excluded
 in part (ii) of Theorem \ref{Thm_sup_jump_equiv}.
We also mention that the proof of part (ii) of our Theorem \ref{Thm_sup_jump_equiv}
 and that of Theorem 4.2 in He and Li \cite{HeLi} use quite different arguments.

Next, we summarize some of the difficulties that we encountered when passing from single-type CBI processes to multi-type ones.
First of all, we point out that, according to our knowledge, only few results are
 available for the analysis of the distributional properties of the jumps of multi-type CBI processes,
 see Chen and Li \cite[Theorem 5.3]{CheLi}, Horst and Xu \cite[Subsection 5.1]{HorXu} and Li \cite[Example 12.2]{Li}.
Compared to single-type CBI processes, the notion of irreducibility plays an important role in the multi-type case,
 in particular, we assume that the underlying multi-type CBI process is irreducible in part (iii) of Proposition \ref{Pro_tv},
 in Corollary \ref{Cor_inf_int_CBI}, in part (iii) of Corollary \ref{Cor_tau},
 in Lemma \ref{Lem_irreducibility}, and in Proposition \ref{Pro_sup_jump}.
Concerning our main Theorem \ref{Thm_sup_jump_equiv}, on the one hand, note that it is valid for
 a general multi-type CBI process satisfying the moment condition \eqref{moment_condition_m_new}
 (in particular, for reducible ones as well), on the other hand, we can handle only
 nondegenerate rectangles in $\RR_+^d$ anchored at $\bzero$ instead of a general Borel set of $\RR_+^d\setminus\{\bzero\}$.

Finally, we introduce the notations that will be used throughout the paper.
Let $\ZZ_+$, $\NN$, $\RR$, $\RR_+$, $\RR_{++}$ and $\CC$ denote the
 set of non-negative integers, positive integers, real numbers, non-negative real
 numbers, positive real numbers and complex numbers, respectively.
For each $d\in\NN$, let $\cU_d:= \RR_+^d \setminus \{\bzero\}$.
For $x , y \in \RR$, we will use the notations
 $x \land y := \min \{x, y\}$ and $x^+ := \max \{0, x\}$.
For $\bx=(x_1,\ldots,x_d)^\top$ and $\by=(y_1,\ldots,y_d)^\top\in\RR^d$,
 the inequality $\bx\leq \by$ means that $x_i\leq y_i$, $i\in\{1,\ldots,d\}$,
 the inequality $\bx<\by$ means that $\bx\leq \by$ and $\bx\ne\by$, and $\langle\bx, \by\rangle := \sum_{j=1}^d x_j y_j$
 denotes the Euclidean inner product.
Given a function $f:\RR^d \to\RR^p$ (where $d,p\in\NN$),
 we say that $f$ is increasing,
 if $f(\bx)\leq f(\by)$ holds for any $\bx,\by\in\RR^d$ with $\bx\leq \by$,
 and we say that $f$ is strictly increasing if $f(\bx)< f(\by)$ holds for any $\bx,\by\in\RR^d$ with $\bx<\by$.
For a function $g:\RR_+\times\RR_+^d\to \RR^d$, $\partial_1g$ denotes the partial derivative of $g$ with
 respect to its first variable (provided that it exists).
By $\|\bx\|$ and $\|\bA\|$, we denote the norm and the induced norm of $\bx \in \RR^d$ and
 $\bA \in \RR^{d\times d}$, respectively.
The null vector and the null matrix will be denoted by $\bzero$.
The open ball around $\bzero$ with radius $\vare>0$ in $\RR_+^d$ is denoted by
 $K_\vare:=\{\by\in\RR_+^d : \Vert \by\Vert<\vare\}$.
Moreover, $\bI_d \in \RR^{d\times d}$ denotes the identity matrix,
 and $\be_1^{(d)}$, \ldots, $\be_d^{(d)}$ denotes the natural bases in $\RR^d$.
The Borel $\sigma$-algebra on a subset $U$ of $\RR^d$ is denoted by $\cB(U)$,
 and recall that $\cB(U) = U\cap \cB(\RR^d)$.
By a Borel measure on a Borel set $S\in\cB(\RR^d)$, we mean a measure on $(S, \cB(S))$.
Every random variable will be defined on an appropriate probability space $(\Omega,\cF,\PP)$.
Throughout this paper, we make the conventions $\int_a^b := \int_{(a,b]}$
 and $\int_a^\infty := \int_{(a,\infty)}$ for any $a, b \in \RR$ with $a \leq b$.

\section{Preliminaries on multi-type CBI processes}
\label{section_CBI}

\begin{Def}\label{Def_essentially_non-negative}
A matrix \ $\bA = (a_{i,j})_{i,j\in\{1,\ldots,d\}} \in \RR^{d\times d}$ \ is
 called essentially non-negative if \ $a_{i,j} \in \RR_+$ \ whenever
 \ $i, j \in \{1,\ldots,d\}$ \ with \ $i \ne j$.
The set of essentially non-negative \ $d \times d$ \ matrices will be denoted by
 \ $\RR^{d\times d}_{(+)}$.
\end{Def}

\begin{Def}\label{Def_admissible}
A tuple \ $(d, \bc, \Bbeta, \bB, \nu, \bmu)$ \ is called a set of admissible
 parameters if
 \renewcommand{\labelenumi}{{\rm(\roman{enumi})}}
 \begin{enumerate}
  \item
   $d \in \NN$,
  \item
   $\bc = (c_i)_{i\in\{1,\ldots,d\}} \in \RR_+^d$,
  \item
   $\Bbeta = (\beta_i)_{i\in\{1,\ldots,d\}} \in \RR_+^d$,
  \item
   $\bB = (b_{i,j})_{i,j\in\{1,\ldots,d\}} \in \RR^{d \times d}_{(+)}$,
  \item
   $\nu$ \ is a Borel measure on \ $\cU_d$
    \ satisfying \ $\int_{\cU_d} (1 \land \|\br\|) \, \nu(\dd\br) < \infty$,
  \item
   $\bmu = (\mu_1, \ldots, \mu_d)$, \ where, for each
    \ $i \in \{1, \ldots, d\}$, \ $\mu_i$ \ is a Borel measure on
    \ $\cU_d$ \ satisfying
    \begin{align}\label{help_moment_mu}
      \int_{\cU_d}
       \biggl[\|\bz\| \land \|\bz\|^2
              + \sum_{j \in \{1, \ldots, d\} \setminus \{i\}} (1 \land z_j)\biggr]
       \mu_i(\dd\bz)
      < \infty .
    \end{align}
  \end{enumerate}
\end{Def}

Note that the measures \ $\nu$ \ and \ $\mu_i$, \ $i \in \{1, \ldots, d\}$, \ in Definition \ref{Def_admissible} are $\sigma$-finite,
 since the function \ $\cU_d\ni \br\mapsto 1 \land \Vert \br\Vert$ \ is strictly positive and measurable with a finite integral with respect to \ $\nu$,
 \ and, for each $i \in \{1, \ldots, d\}$, \ the function
 \ $\cU_d\ni \bz \mapsto \|\bz\| \land \|\bz\|^2 + \sum_{j \in \{1, \ldots, d\} \setminus \{i\}} (1 \land z_j)$ \
 is strictly positive and measurable with a finite integral with respect to \ $\mu_i$, \ see, e.g., Kallenberg \cite[Lemma 1.4]{K2}.

\begin{Thm}\label{CBI_exists}
Let \ $(d, \bc, \Bbeta, \bB, \nu, \bmu)$ \ be a set of admissible parameters.
Then there exists a unique conservative transition semigroup \ $(P_t)_{t\in\RR_+}$
 \ acting on the Banach space (endowed with the supremum norm) of real-valued
 bounded Borel measurable functions on the state space \ $\RR_+^d$ \ such that its
 Laplace transform has a representation
 \begin{equation*}
  \int_{\RR_+^d} \ee^{- \langle \blambda, \by \rangle} P_t(\bx, \dd \by)
  = \ee^{- \langle \bx, \bv(t, \blambda) \rangle
         - \int_0^t \psi(\bv(s, \blambda)) \, \dd s} , \qquad
  \bx \in \RR_+^d, \quad \blambda \in \RR_+^d , \quad t \in \RR_+ ,
 \end{equation*}
 where, for any \ $\blambda \in \RR_+^d$, \ the continuously differentiable
 function
 \ $\RR_+ \ni t \mapsto \bv(t, \blambda)
    = (v_1(t, \blambda), \ldots, v_d(t, \blambda))^\top \in \RR_+^d$
 \ is the unique locally bounded solution to the system of differential equations
 \begin{align*}
   \partial_1 v_i(t, \blambda) = - \varphi_i(\bv(t, \blambda)) ,  \qquad t\in\RR_+, \qquad
   v_i(0, \blambda) = \lambda_i , \qquad i \in \{1, \ldots, d\} ,
 \end{align*}
 with
 \[
   \varphi_i(\blambda)
   := c_i \lambda_i^2 -  \langle \bB \be_i^{(d)}, \blambda \rangle
      + \int_{\cU_d}
         \bigl( \ee^{- \langle \blambda, \bz \rangle} - 1
                + \lambda_i (1 \land z_i) \bigr)
         \, \mu_i(\dd \bz)
 \]
 for \ $\blambda \in \RR_+^d$, \ $i \in \{1, \ldots, d\}$, \ and
 \begin{align*}
   \psi(\blambda)
   := \langle \bbeta, \blambda \rangle
      + \int_{\cU_d}
         \bigl( 1 - \ee^{- \langle\blambda, \br\rangle} \bigr)
         \, \nu(\dd\br) , \qquad
   \blambda \in \RR_+^d .
 \end{align*}
\end{Thm}

Theorem \ref{CBI_exists} is a special case of Theorem 2.7 of Duffie et al.\ \cite{DufFilSch}
 with \ $m = d$, \ $n = 0$ \ and zero killing rate.
For the unique existence of a locally bounded solution to the system of differential equations
 in Theorem \ref{CBI_exists}, see Li \cite[page 48]{Li} or Duffie et al.\ \cite[Proposition 6.4]{DufFilSch}.

\begin{Def}\label{Def_CBI}
A conservative Markov process with state space \ $\RR_+^d$ \ and with transition
 semigroup \ $(P_t)_{t\in\RR_+}$ \ given in Theorem \ref{CBI_exists} is called a
 multi-type CBI process with parameters \ $(d, \bc, \Bbeta, \bB, \nu, \bmu)$.
\ The function
 \ $\RR_+^d \ni \blambda
    \mapsto \bvarphi(\blambda)=(\varphi_1(\blambda), \ldots, \varphi_d(\blambda))^\top \in \RR^d$
 \ is called its branching mechanism, and the function
 \ $\RR_+^d \ni \blambda \mapsto \psi(\blambda) \in \RR_+$ \ is called its
 immigration mechanism.
When there is no immigration, i.e., $\Bbeta = \bzero$ and $\nu = 0$,
 the process is simply called a multi-type CB process (a continuous state and continuous time branching process).
\end{Def}

Given a set of admissible parameters $(d, \bc, \Bbeta, \bB, \nu, \bmu)$, we get that
 $\nu$ and $\mu_i$, $i\in\{1,\ldots,d\}$, are L\'evy measures on $\cU_d$, since, by parts (v) and (vi) of Definition \ref{Def_admissible},
 we have
 \[
   \int_{\cU_d} (1 \land \|\br\|^2) \, \nu(\dd\br) \leq \int_{\cU_d}  (1 \land \|\br\|) \, \nu(\dd\br)< \infty,
 \]
 and
  \[
    \int_{\cU_d} (1\land \|\bz\|^2) \, \mu_i(\dd\bz)
       \leq \int_{\cU_d} (\|\bz\| \land \|\bz\|^2) \, \mu_i(\dd\bz)
      <\infty, \qquad i\in\{1,\ldots,d\}.
  \]
For this reason, we call $\nu+\sum_{i=1}^d \mu_i$ the total L\'evy measure corresponding to a multi-type CBI process
 with parameters $(d, \bc, \Bbeta, \bB, \nu, \bmu)$.

By Barczy et al.\ \cite[Remark 2.3 and (2.12)]{BarLiPap2}, for each $i\in\{1,\ldots,d\}$, the moment condition \eqref{help_moment_mu}
 is equivalent to
 \begin{align*}
      \int_{\cU_d}
       \biggl[\|\bz\| \land \|\bz\|^2
              + \sum_{j \in \{1, \ldots, d\} \setminus \{i\}} z_j\biggr]
       \mu_i(\dd\bz)
      < \infty ,
 \end{align*}
 which coincides with the moment condition in Example 2.5 in Li \cite[page 48]{Li} for multi-type CB processes,
 where multi-type CB processes are considered as special superprocesses.

By Li \cite[Theorem A.7]{Li}, $(\bX_t)_{t\in\RR_+}$ has c\`{a}dl\`{a}g realizations, and any such realization of the process has
 a c\`{a}dl\`{a}g modification $(\widetilde\bX_t)_{t\in\RR_+}$, and hence $\PP(\bX_t = \widetilde\bX_t)=1$, $t\in\RR_+$,
 and all the sample paths of $(\widetilde\bX_t)_{t\in\RR_+}$ are right continuous at every $t\in\RR_+$
 and possesses left limit at every $t\in\RR_{++}$.

Now, we present a property of \ $\psi$, \ which will be used in the proof of Corollary \ref{Cor_tau} as well.

\begin{Lem}\label{Lem_psi_positive}
If the immigration mechanism \ $\psi$ \ given in Definition \ref{Def_CBI} is not identically zero,
 then \ $\psi(\blambda)>0$ \ for all \ $\blambda\in(0,\infty)^d$.
\end{Lem}

\noindent{\bf Proof.}
Let us suppose that $\psi$ is not identically zero.
Then \ $\bbeta \ne \bzero$ \ or \ $\nu\ne 0$.
If \ $\bbeta\ne\bzero$, \ then, using that $\bbeta\in\RR_+^d$,
 there exists \ $i_0\in\{1,\ldots,d\}$ \ such that \ $\beta_{i_0}>0$, \ and hence, for all \ $\blambda\in(0,\infty)^d$,
 \ we have \ $\psi(\blambda)\geq \langle \bbeta, \blambda \rangle\geq \lambda_{i_0}\beta_{i_0}>0$.
\ If \ $\bbeta = \bzero$, \ then \ $\nu\ne 0$, \ and, similarly as before,
 for all \ $\blambda\in(0,\infty)^d$, \ we have \ $\langle \blambda, \br \rangle>0$, \ $\br\in \cU_d$, \
 and hence \ $1 - \ee^{-\langle \blambda, \br \rangle} > 0$, \ $\br\in \cU_d$.
\ This together with \ $\nu(\cU_d)>0$ \ (following from the fact that \ $\nu$ \ is not identically zero)
 yield that \ $\psi(\blambda) = \int_{\cU_d} \bigl( 1 - \ee^{- \langle\blambda, \br\rangle} \bigr) \, \nu(\dd\br) > 0$
 for all \ $\blambda\in(0,\infty)^d$.
\proofend

For a multi-type CBI process $(\bX_t)_{t\in\RR_+}$, $\bx\in\RR_+^d$, an event $A\in\sigma(\bX_t,t\in\RR_+)$ and an $\RR^p$-valued
 random variable $\bxi$ which is $\sigma(\bX_t,t\in\RR_+)$-measurable (where $p\in\NN$), let $\PP_{\bx}(A):=\PP(A\mid \bX_0=\bx)$
 and $\EE_{\bx}(\xi):=\EE(\xi \mid \bX_0=\bx)$, respectively.

Let \ $(\bX_t)_{t\in\RR_+}$ \ be a multi-type CBI process with parameters
 \ $(d, \bc, \Bbeta, \bB, \nu, \bmu)$ \ such that the moment condition
 \begin{equation}\label{moment_condition_m_new}
  \int_{\cU_d} \|\br\| \bbone_{\{\|\br\|\geq1\}} \, \nu(\dd\br) < \infty
 \end{equation}
 holds.
Then, by formula (3.4) in Barczy et al.\ \cite{BarLiPap2} (see also formula (79) in Li \cite{Li3}),
 \begin{equation}\label{EXcond}
  \EE_{\bx}(\bX_t)
  = \ee^{t\tbB} \bx + \int_0^t \ee^{u\tbB} \tBbeta \, \dd u ,
  \qquad \bx \in \RR_+^d , \quad t \in \RR_+ ,
 \end{equation}
 where
 \begin{gather}\label{help5}
  \tbB := (\tb_{i,j})_{i,j\in\{1,\ldots,d\}} , \qquad
  \tb_{i,j}
  := b_{i,j} + \int_{\cU_d} (z_i - \delta_{i,j})^+ \, \mu_j(\dd\bz) , \qquad
  \tBbeta := \Bbeta + \int_{\cU_d} \br \, \nu(\dd\br) ,
 \end{gather}
 with \ $\delta_{i,j}:=1$ \ if \ $i = j$, \ and \ $\delta_{i,j} := 0$ \ if
 \ $i \ne j$.
\ Note that, for all \ $\bx \in \RR^d_+$, \ the function
 \ $\RR_+ \ni t \mapsto \EE_{\bx}(\bX_t)$ \ is continuous, and
 \ $\tbB \in \RR^{d \times d}_{(+)}$ \ and \ $\tBbeta \in \RR_+^d$.
Indeed, part (v) of Definition \ref{Def_admissible} together with the moment condition \eqref{moment_condition_m_new}
 and Barczy et al.\ \cite[Remark 2.3 and formulas (2.11) and (2.12)]{BarLiPap2} yield that
 \[
   \int_{\cU_d} \|\br\| \, \nu(\dd\br) < \infty , \qquad
   \int_{\cU_d} (z_i - \delta_{i,j})^+ \, \mu_j(\dd \bz) < \infty , \quad
   i, j \in \{1, \ldots, d\} .
 \]
Here we point out that the fact that the matrix $\tbB$ belongs to $\RR^{d \times d}_{(+)}$ even if the moment condition
 \eqref{moment_condition_m_new} does not hold, however, the vector $\tBbeta$ belongs to $\RR_+^d$ if and only if
 the moment condition \eqref{moment_condition_m_new} holds.

Given a set of admissible parameters $(d, \bc, \Bbeta, \bB, \nu, \bmu)$ such that the moment condition
 \eqref{moment_condition_m_new} holds, let us consider the stochastic differential equation (SDE)
 \begin{align}\label{SDE_atirasa_dimd}
  \begin{split}
   \bX_t
   &=\bX_0
     + \int_0^t (\Bbeta + \tbB \bX_u) \, \dd u
     + \sum_{\ell=1}^d
        \int_0^t \sqrt{2 c_\ell \max \{0, X_{u,\ell}\}} \, \dd W_{u,\ell}
        \, \be_\ell^{(d)} \\
   &\quad
      + \sum_{\ell=1}^d
         \int_0^t \int_{\cU_d} \int_{\cU_1}
          \bz \bbone_{\{w\leq X_{u-,\ell}\}} \, \tN_\ell(\dd u, \dd\bz, \dd w)
      + \int_0^t \int_{\cU_d} \br \, M(\dd u, \dd\br)
  \end{split}
 \end{align}
 for $t \in\RR_+$, where
 $X_{t,\ell}$, $\ell \in \{1, \ldots, d\}$, denotes the $\ell^{\mathrm th}$ coordinate of $\bX_t$,
  $(W_{t,1})_{t\in\RR_+}$,$\ldots$,$(W_{t,d})_{t\in\RR_+}$ are standard Wiener processes,
 \ $N_\ell$, \ $\ell \in \{1, \ldots, d\}$, \ and \ $M$ \ are Poisson
 random measures on \ $\cU_1 \times \cU_d \times \cU_1$ \ and on
 \ $\cU_1 \times \cU_d$ \ with intensity measures
 \ $\dd u \, \mu_\ell(\dd\bz) \, \dd w$, \ $\ell \in \{1, \ldots, d\}$, \ and
 \ $\dd u \, \nu(\dd\br)$, \ respectively, and
 \ $\tN_\ell(\dd u, \dd\bz, \dd w)
    := N_\ell(\dd u, \dd\bz, \dd w) - \dd u \, \mu_\ell(\dd\bz) \, \dd w$,
 \ $\ell \in \{1, \ldots, d\}$.
\ We suppose that $\EE(\Vert \bX_0\Vert)<\infty$ and that \ $\bX_0$, \ $(W_{t,1})_{t\in\RR_+}$, \ldots, \ $(W_{t,d})_{t\in\RR_+}$, \ $N_1$,
 \ldots, \ $N_d$ \ and \ $M$ \ are mutually independent.
The SDE \eqref{SDE_atirasa_dimd} has a pathwise unique $\RR_+^d$-valued c\`{a}dl\`{a}g strong solution,
 and the solution is a multi-type CBI process with parameters \ $(d, \bc, \Bbeta, \bB, \nu, \bmu)$,
 \ see Theorem 4.6 and Section 5 in
 Barczy et al.~\cite{BarLiPap2}, where \eqref{SDE_atirasa_dimd} was proved only for
 \ $d \in \{1, 2\}$, \ but their method clearly works for all \ $d \in\NN$.
Consequently, given a c\`{a}dl\`{a}g CBI process $(\bX_t)_{t\in\RR_+}$ with parameters $(d, \bc, \Bbeta, \bB, \nu, \bmu)$
 such that $\EE(\Vert \bX_0\Vert)<\infty$ and the moment condition \eqref{moment_condition_m_new} hold,
 its law on the space of $\RR^d$-valued c\`{a}dl\`{a}g functions defined on $\RR_+$ coincides with
 the law of the pathwise unique c\`{a}dl\`{a}g strong solution of the SDE \eqref{SDE_atirasa_dimd}.
In the remaining part of the paper, when we refer to a CBI process $(\bX_t)_{t\in\RR_+}$ with parameters \ $(d, \bc, \Bbeta, \bB, \nu, \bmu)$
 such that $\EE(\Vert \bX_0\Vert)<\infty$ and the moment condition \eqref{moment_condition_m_new} hold, we consider it
 as a pathwise unique c\`{a}dl\`{a}g strong solution of the SDE \eqref{SDE_atirasa_dimd}.

Finally, we recall the notion of irreducibility for a matrix and for a multi-type CBI process.
A matrix \ $\bA \in \RR^{d\times d}$ \ is called reducible if there exist a
 permutation matrix \ $\bP \in \RR^{ d \times d}$ \ and an integer \ $p$ \ with
 \ $1 \leq p \leq d-1$ \ such that
 \[
  \bP^\top \bA \bP
   = \begin{pmatrix} \bA_1 & \bA_2 \\ \bzero & \bA_3 \end{pmatrix},
 \]
 where \ $\bA_1 \in \RR^{p\times p}$, \ $\bA_2 \in \RR^{p \times (d-p)}$, \ $\bA_3 \in \RR^{ (d-p) \times (d-p) }$, \ and \ $\bzero \in \RR^{(d-p)\times p}$ \ is
 a null matrix.
A matrix \ $\bA \in \RR^{d\times d}$ \ is called irreducible if it is not
 reducible, see, e.g., Horn and Johnson
 \cite[Definitions 6.2.21 and 6.2.22]{HorJoh}.
We do emphasize that no 1-by-1 matrix is reducible.

\begin{Def}\label{Def_irreducible}
Let \ $(\bX_t)_{t\in\RR_+}$ \ be a multi-type CBI process with parameters
 \ $(d, \bc, \Bbeta, \bB, \nu, \bmu)$.
Then \ $(\bX_t)_{t\in\RR_+}$ \ is called irreducible if \ $\tbB$ \ is irreducible.
\end{Def}

We point out the fact that for the definition of irreducibility of a multi-type CBI process in Definition \ref{Def_irreducible},
 we do not need the moment condition \eqref{moment_condition_m_new} on $\nu$, so the notion of irreducibility is a property only of
 the branching mechanism of a multi-type CBI process.
Note that every single-type CBI process is irreducible, and hence irreducibility comes into play for multi-type CBI processes with at least two types.
Recall that \ $\tbB \in \RR^{d\times d}_{(+)}$ \ is irreducible if and only if
 \ $\ee^{t\tbB} \in \RR^{d \times d}_{++}$ \ for all \ $t \in \RR_{++}$
 (see Barczy and Pap \cite[Lemma A.1]{BarPap}).
If $(\bX_t)_{t\in\RR_+}$ is an irreducible multi-type CBI process and $\bx\in\RR_+^d$ is such that at least one of its coordinates is positive
 (i.e., $\bx\ne\bzero$), then \eqref{EXcond} yields that $\EE_{\bx}(\bX_t)\in\RR_{++}^d$ for all $t\in\RR_{++}$.
Roughly speaking, irreducibility of a multi-type CBI process implies that each type continuously generates mass of
 all the types.

\section{A multi-type CBI process and its integral process}\label{Section_CBI_intCBI}

The next result states that a $(2d)$-dimensional stochastic process having coordinate processes a $d$-type CBI process
 and its integral process is a $(2d)$-type CBI process, where \ $d\in\NN$.

\begin{Pro}\label{Pro_2dCBI}
Let \ $(\bX_t)_{t\in\RR_+}$ \ be a multi-type CBI process with parameters
 \ $(d, \bc, \Bbeta, \bB, \nu, \bmu)$ \
 such that \ $\EE(\Vert \bX_0\Vert)<\infty$ \ and the moment condition \eqref{moment_condition_m_new} hold.
Let \ $\bY_0$ \ be an $\RR_+^d$-valued random variable such that \ $\EE(\Vert \bY_0\Vert)<\infty$ \
 and it is independent of \ $(W_{t,1})_{t\in\RR_+}$, \ldots, \ $(W_{t,d})_{t\in\RR_+}$, \ $N_1$, \ldots, \ $N_d$ \ and \ $M$
 \ appearing in the SDE \eqref{SDE_atirasa_dimd}.
Let \ $\bY_t:=\bY_0 + \int_0^t \bX_u\,\dd u$, \ $t\in\RR_+$.
\ Then \ $(\bX_t, \bY_t)_{t\in\RR_+}$ \ is a $(2d)$-type CBI process with branching mechanism
 \ $\bvarphi^*:\RR_+^{2d}\to\RR^{2d}$,
 \begin{align}\label{2dCBI_branching}
  \bvarphi^*(\blambda,\tblambda)
       := (\varphi_1(\blambda)-\tlambda_1, \ldots, \varphi_d(\blambda)-\tlambda_d,0,\ldots,0),
       \qquad (\blambda,\tblambda)\in\RR_+^{2d},
 \end{align}
 and with immigration mechanism \ $\psi^*:\RR_+^{2d}\to\RR_+$, \
 \begin{align}\label{2dCBI_immigration}
    \psi^*(\blambda,\tblambda):=\psi(\blambda),\qquad  (\blambda,\tblambda)\in\RR_+^{2d},
 \end{align}
 where \ $\bvarphi=(\varphi_1,\ldots,\varphi_d):\RR_+^d\to\RR^d$ \ and \ $\psi:\RR_+^d\to\RR_+$ \ is the branching mechanism
 and the immigration mechanism of \ $(\bX_t)_{t\in\RR_+}$, \ respectively, and \ $\tblambda:=(\tlambda_1,\ldots,\tlambda_d)$.
\end{Pro}

\noindent{\bf Proof.}
Using the SDE \eqref{SDE_atirasa_dimd} for \ $(\bX_t)_{t\in\RR_+}$, \ we have
 \begin{align*}
  \begin{bmatrix}
    \bX_t \\
    \bY_t \\
  \end{bmatrix}
  & = \begin{bmatrix}
      \bX_0 \\
      \bY_0 \\
     \end{bmatrix}
     + \int_0^t \left(
                      \begin{bmatrix}
                          \Bbeta \\
                          \bzero \\
                       \end{bmatrix}
                     + \begin{bmatrix}
                         \tbB & \bzero \\
                         \bI_d & \bzero \\
                       \end{bmatrix}
                       \begin{bmatrix}
                         \bX_u \\
                         \bY_u \\
                       \end{bmatrix}
                      \right) \, \dd u
     + \begin{bmatrix}
         \sum_{\ell=1}^d
        \int_0^t \sqrt{2 c_\ell \max \{0, X_{u,\ell}\}} \, \dd W_{u,\ell} \, \be_\ell^{(d)} \\
         \bzero \\
       \end{bmatrix}\\
  &\quad + \begin{bmatrix}
             \sum_{\ell=1}^d
             \int_0^t \int_{\cU_d} \int_{\cU_1}
              \bz \bbone_{\{w\leq X_{u-,\ell}\}} \, \tN_\ell(\dd u, \dd\bz, \dd w) \\
             \bzero \\
           \end{bmatrix}
         + \begin{bmatrix}
             \int_0^t \int_{\cU_d} \br \, M(\dd u, \dd\br) \\
             \bzero \\
           \end{bmatrix},\qquad t\in\RR_+.
 \end{align*}
Let $(\tW_{t,\ell})_{t\in\RR_+}$, $\ell\in\{1,\ldots,d\}$, be mutually independent standard Wiener processes such that they
 are independent of $(W_{t,\ell})_{t\in\RR_+}$, $N_\ell$, $\ell\in\{1,\ldots,d\}$, and $M$.
Let \ $N_\ell^*$, \ $\ell \in \{1, \ldots, d\}$, \ and \ $M^*$ \ be the embeddings of
 \ $N_\ell$, \ $\ell \in \{1, \ldots, d\}$, \ and \ $M$ \ into \ $\cU_1 \times (\cU_{d}\times\{\bzero\}) \times \cU_1$ \ and into
 \ $\cU_1\times (\cU_{d}\times\{\bzero\})$, respectively (where $\bzero\in\RR^d$).
Then, by the mapping theorem for Poisson point processes (see, e.g., Kingman \cite[Section 2.3]{Kin}),
 we get that \ $N_\ell^*$, \ $\ell \in \{1, \ldots, d\}$, \ and \ $M^*$ \ are Poisson random measures
 on \ $\cU_1 \times \cU_{2d} \times \cU_1$ \ and on \ $\cU_1\times \cU_{2d}$ \ with intensity measures
 \ $\dd u \, \mu_\ell^*(\dd\bz) \, \dd w$, \ $\ell \in \{1, \ldots, d\}$, \ and
 \ $\dd u \, \nu^*(\dd\br)$, \ respectively, given by
 \begin{align*}
  &\mu_\ell^*(B)= \int_{\cU_d} \bbone_{B}(\bz,\bzero)\,\mu_\ell(\dd \bz),\qquad B\in\cB(\cU_{2d}),\qquad \ell\in\{1,\ldots,d\},\\
  &\nu^*(B)= \int_{\cU_d} \bbone_{B}(\bz,\bzero)\,\nu(\dd \bz),\qquad B\in\cB(\cU_{2d}).
 \end{align*}
For each $\ell\in\{d+1,\ldots,2d\}$, let $N_\ell^*$ be the zero Poisson random measure on \ $\cU_1 \times \cU_{2d} \times \cU_1$
 with intensity measure \ $\dd u \, \mu_\ell^*(\dd\bz) \, \dd w$, \ where $\mu_\ell^*(B):= 0,$ $B\in\cB(\cU_{2d})$.
Since $(\tW_{t,\ell})_{t\in\RR_+}$, $(W_{t,\ell})_{t\in\RR_+}$, $N_\ell$, $\ell\in\{1,\ldots,d\}$, and $M$
 are mutually independent, we have that $(\tW_{t,\ell})_{t\in\RR_+}$, $(W_{t,\ell})_{t\in\RR_+}$, $\ell\in\{1,\ldots,d\}$, $N_\ell^*$, $\ell\in\{1,\ldots,2d\}$,
 and $M^*$ are mutually independent as well.
Then we have
 \begin{align*}
  \begin{bmatrix}
    \bX_t \\
    \bY_t \\
   \end{bmatrix}
  &= \begin{bmatrix}
      \bX_0 \\
      \bY_0 \\
     \end{bmatrix}
    + \int_0^t \left(
                      \begin{bmatrix}
                          \Bbeta \\
                          \bzero \\
                       \end{bmatrix}
                     + \begin{bmatrix}
                         \tbB & \bzero \\
                         \bI_d & \bzero \\
                       \end{bmatrix}
                       \begin{bmatrix}
                         \bX_u \\
                         \bY_u \\
                       \end{bmatrix}
                      \right) \, \dd u
     + \sum_{\ell=1}^d
        \int_0^t \sqrt{2 c_\ell \max \{0, X_{u,\ell}\}} \, \dd W_{u,\ell}
        \, \be_\ell^{(2d)} \\
   &\quad
     +\sum_{\ell=1}^d
        \int_0^t 0 \, \dd \tW_{u,\ell} \, \be_{d+\ell}^{(2d)}
     + \sum_{\ell=1}^d
         \int_0^t \int_{\cU_{2d}} \int_{\cU_1}
          \bz^* \bbone_{\{w\leq X_{u-,\ell}\}} \, \tN_\ell^*(\dd u, \dd\bz^*, \dd w)  \\
   &\quad
      + \sum_{\ell=1}^d
         \int_0^t \int_{\cU_{2d}} \int_{\cU_1}
          \bz^* \bbone_{\{w\leq Y_{u-,\ell}\}} \, \tN_{d+\ell}^*(\dd u, \dd\bz^*, \dd w)
      + \int_0^t \int_{\cU_{2d}} \br^* \, M^*(\dd u, \dd\br^*),\qquad t\in\RR_+.
 \end{align*}
Let
 \begin{align*}
      &\Bbeta^*:= \begin{bmatrix}
                          \Bbeta \\
                          \bzero \\
                       \end{bmatrix}
                      \in\RR_+^{2d},
      \qquad
      \bc^*:= \begin{bmatrix}
                  \bc \\
                  \bzero \\
              \end{bmatrix}
                 \in\RR_+^{2d},
      \qquad \bB^*:= \begin{bmatrix}
                  \bB & \bzero \\
                  \bI_d & \bzero \\
              \end{bmatrix}
             \in\RR^{(2d)\times(2d)}_{(+)}.
  \end{align*}
Note that for any Borel measurable function \ $f:\cU_{2d}\to\RR_+$ \ being integrable with respect to \ $\nu^*$,
 \ we have that
 \begin{align}\label{help_int1}
  \int_{\cU_{2d}} f(\br^*)\,\nu^*(\dd\br^*)
       = \int_{\cU_d} f(\br,\bzero)\,\nu(\dd\br),
 \end{align}
 as a consequence of the construction of the Lebesgue integral.
Similarly, for any Borel measurable function \ $f:\cU_{2d}\to\RR_+$ \ being integrable with respect to \ $\mu_\ell^*$, \ where \ $\ell\in\{1,\ldots,d\}$,
 \ we have that
 \begin{align}\label{help_int2}
  \int_{\cU_{2d}} f(\bz^*)\,\mu_\ell^*(\dd\bz^*)
       = \int_{\cU_d} f(\bz,\bzero)\,\mu_\ell(\dd\bz).
 \end{align}
Then, by \eqref{help5}, we have
 \[
  \widetilde{\bB^*} = \bB^*
            + \begin{bmatrix}
                  \left(\int_{\cU_d} (z_i - \delta_{i,j})^+\,\mu_j(\dd\bz) \right)_{i,j=1}^d & \bzero \\
                  \bzero & \bzero \\
              \end{bmatrix}
        = \begin{bmatrix}
                  \tbB & \bzero \\
                  \bI_d & \bzero \\
              \end{bmatrix},
 \]
 since, by \eqref{help_int2}, for all $i\in\{d+1,\ldots,2d\}$ and $j\in\{1,\ldots,d\}$, we have
 \[
  \int_{\cU_{2d}} (z_i^* - \delta_{i,j})^+\,\mu_j^*(\dd\bz^*)
    =  \int_{\cU_d} (0 - 0)^+\,\mu_j(\dd\bz)
    =0.
 \]
Furthermore, by \eqref{moment_condition_m_new}, \eqref{help_int1}, \eqref{help_int2} and the moment conditions given in parts (v), (vi) of Definition \ref{Def_admissible},
 we have that
 \begin{align*}
   & \int_{\cU_{2d}} (1 \land \|\br^*\|) \, \nu^*(\dd\br^*)
       \leq \int_{\cU_{2d}} \|\br^*\| \, \nu^*(\dd\br^*)
       = \int_{\cU_d} \|\br\| \, \nu(\dd\br)
       <\infty,\\
   &  \int_{\cU_{2d}}
       \biggl[\|\bz^*\| \land \|\bz^*\|^2
              + \sum_{j \in \{1, \ldots, 2d\} \setminus \{i\}} (1 \land z_j^*)\biggr]
       \mu_i^*(\dd\bz^*) \\
   &\qquad   =  \int_{\cU_d}
       \biggl[\|\bz\| \land \|\bz\|^2
              + \sum_{j \in \{1, \ldots, d\} \setminus \{i\}} (1 \land z_j)\biggr]
       \mu_i(\dd\bz) <\infty, \qquad i\in\{1,\ldots,d\}\\
   & \int_{\cU_{2d}}
       \biggl[\|\bz^*\| \land \|\bz^*\|^2
              + \sum_{j \in \{1, \ldots, 2d\} \setminus \{i\}} (1 \land z_j^*)\biggr]
       \mu_i^*(\dd\bz^*) = 0 <\infty, \qquad i\in\{d+1,\ldots,2d\}.
 \end{align*}
In particular, the moment condition \eqref{moment_condition_m_new} holds for the measure \ $\nu^*$.
\ Thus, by Theorem 4.6 in Barczy et al.\ \cite{BarLiPap2}, \ $(\bX_t,\bY_t)_{t\in\RR_+}$ \ is a CBI process
 with parameters \ $(2d, \bc^*, \Bbeta^*, \bB^*, \nu^*, \bmu^*)$, \ where \ $\bmu^*:=(\mu_1^*, \ldots,\mu_{2d}^*)$.
For all \ $(\blambda,\tblambda)\in\RR_+^{2d}$, we have
 \begin{align*}
  \left\langle \bB^* \be_i^{(2d)}, \begin{bmatrix}
                                 \blambda \\
                                 \tblambda \\
                               \end{bmatrix} \right\rangle
    = \begin{cases}
          \langle \bB \be_i^{(d)},\blambda \rangle + \tlambda_i  & \text{if \ $i\in\{1,\ldots,d\}$,}\\
          0   & \text{if \ $i\in\{d+1,\ldots,2d\}$,}
       \end{cases}
 \end{align*}
 and for all \ $(\blambda,\tblambda)\in\RR_+^{2d}$, using \eqref{help_int1} with \ $f(\br^*):=1 - \ee^{- \langle (\blambda, \tblambda), \br^*\rangle}$, \ $\br^*\in\cU_{2d}$,
 and \eqref{help_int2} with \ $f(\bz^*):= \ee^{- \langle (\blambda, \tblambda), \bz^* \rangle} - 1 + \lambda_i (1 \land z_i^*)$, \ $\bz^*\in\cU_{2d}$,
 \ where \ $i\in\{1,\ldots,d\}$, \ we get
  \[
    \int_{\cU_{2d}}
         \bigl( 1 - \ee^{- \langle (\blambda, \tblambda), \br^*\rangle} \bigr)
         \, \nu^*(\dd\br^*)
     = \int_{\cU_d}
         \bigl( 1 - \ee^{- \langle\blambda, \br\rangle} \bigr)
         \, \nu(\dd\br)
  \]
  and
  \[
  \int_{\cU_{2d}}
         \bigl( \ee^{- \langle (\blambda, \tblambda), \bz^* \rangle} - 1 + \lambda_i (1 \land z_i^*) \bigr)
         \, \mu_i^*(\dd \bz^*)
   =
   \int_{\cU_d}
         \bigl( \ee^{- \langle \blambda, \bz \rangle} - 1
                + \lambda_i (1 \land z_i) \bigr)
         \, \mu_i(\dd \bz), \qquad i\in\{1,\ldots,d\}.
  \]
Consequently, we have that the CBI process \ $(\bX_t, \bY_t)_{t\in\RR_+}$ \ has branching and immigration mechanisms given in
 \eqref{2dCBI_branching} and \eqref{2dCBI_immigration}, respectively, as desired.
\proofend

\begin{Rem}\label{Rem_discussion}
We note that the fact that \ $(\bX_t, \bY_t)_{t\in\RR_+}$ \ is a $(2d)$-dimensional CBI process
 also follows from Filipovi\'c et al.\ \cite[paragraph before Theorem 4.3]{FilMaySch} or
 Theorem 4.10 in Keller-Ressel \cite{Kel}, where this property is stated for general regular affine processes.
Filipovi\'c et al.\ \cite{FilMaySch} did not give any proof, Keller-Ressel \cite{Kel} gave a proof,
 which is completely different from ours, since he calculated conditional moment generating function of \ $(\bX_t, \bY_t)_{t\in\RR_+}$,
 while our proof is based on the SDE \eqref{SDE_atirasa_dimd} for \ $(\bX_t)_{t\in\RR_+}$.
The formulae derived for the branching and immigration mechanisms of \ $(\bX_t,\bY_t)_{t\in\RR_+}$ \
 in Proposition \ref{Pro_2dCBI} are also in accordance with the ones in Theorem 4.10 in Keller-Ressel \cite{Kel}.
\proofend
\end{Rem}

Next, we present a result on the Laplace transform of the integral of a multi-type CBI process, see Proposition \ref{Pro_int_CBI_Laplace}.
In fact, this result is a special case of Theorem 9.22 in Li \cite{Li} for immigration superprocesses
 or of Corollary 4.11 in Keller-Ressel \cite{Kel} for analytic affine processes, but, for completeness, we give an independent proof.
Proposition \ref{Pro_int_CBI_Laplace} with \ $d=1$ \ gives back Proposition 2.1 in He and Li \cite{HeLi}
 and Theorem 8.12 with $\lambda=0$ in Li \cite{Li3} as well.
The proof of Proposition 2.1 in He and Li \cite{HeLi} and that of Theorem 8.12 in Li \cite{Li3} are different from
 the proof of our Proposition \ref{Pro_int_CBI_Laplace}.
We directly use that \ $(\bX_t,\int_0^t \bX_u\,\dd u)_{t\in\RR_+}$ \ is a \ $(2d)$-type CBI process starting from \ $(\bX_0,\bzero)$
 (see Proposition \ref{Pro_2dCBI}), while Li \cite[Theorem 8.12]{Li3} used It\^{o}'s formula for appropriate functions of
 $(t,\bX_t,\int_0^t \bX_u\,\dd u)$, \ where $t\in\RR_+$.
\ The proof of Proposition 2.1 in He and Li \cite{HeLi} is an application of Theorem 9.22 in Li \cite{Li} for immigration superprocesses.

\begin{Pro}\label{Pro_int_CBI_Laplace}
Let \ $(\bX_t)_{t\in\RR_+}$ \ be a multi-type CBI process with parameters
 \ $(d, \bc, \Bbeta, \bB, \nu, \bmu)$ \ such that the moment condition  \eqref{moment_condition_m_new} holds.
Then
 \[
   \EE_\bx\left( \exp\left\{ - \left\langle \tblambda, \int_0^t \bX_u \,\dd u  \right\rangle \right\} \right)
     = \exp\left\{ -\langle \bx,  \tbv(t,\tblambda)\rangle
                   - \int_0^t \psi(\tbv(s,\tblambda) ) \,\dd s \right\}
 \]
 for \ $t\in\RR_+$, \ $\tblambda\in \RR_+^d$ \ and \ $\bx\in\RR_+^d$, \ where,
 for all \ $\tblambda = (\tlambda_1,\ldots,\tlambda_d)^\top\in\RR_+^d$, \ the continuously differentiable function
 \ $\RR_+\ni t\mapsto \tbv(t,\tblambda) =: ( \tv_1(t,\tblambda),\ldots, \tv_d(t,\tblambda) )^\top\in\RR_+^d$ \
 is the unique locally bounded solution to the system of differential equations
 \begin{align}\label{help1}
  \partial_1\tv_i(t, \tblambda) = \tlambda_i - \varphi_i(\tbv(t, \tblambda)) , \qquad
   \tv_i(0, \tblambda) = 0 , \qquad i \in \{1, \ldots, d\} .
 \end{align}
\end{Pro}

\noindent{\bf Proof.}
By Proposition \ref{Pro_2dCBI}, for all $\bx\in\RR_+^d$, under the probability measure $\PP_{\bx}$,
 we have that \ $(\bX_t,\int_0^t \bX_u\,\dd u)_{t\in\RR_+}$ \ is a \ $(2d)$-type CBI process  starting from \ $(\bx,\bzero)$.
Hence Theorem \ref{CBI_exists} and \eqref{2dCBI_immigration} yield that
 \begin{align*}
  &\EE_{\bx}\left( \exp\left\{ - \left\langle \tblambda, \int_0^t \bX_u \,\dd u  \right\rangle \right\} \right)
     = \EE_{\bx}\left( \exp\left\{ - \sum_{\ell=1}^d \tlambda_\ell \int_0^t X_{u,\ell} \,\dd u \right\} \right)\\
  &\qquad  = \exp\left\{ -\left\langle \begin{bmatrix}
                                   \bx \\
                                   \bzero \\
                                 \end{bmatrix},
                               \bv^*(t,\bzero,\tblambda)\right\rangle
                   - \int_0^t \psi^*(\bv^*(s,\bzero,\tblambda) ) \,\dd s \right\}\\
  &\qquad = \exp\left\{ -\left\langle  \bx,  \begin{bmatrix}
                                          v^*_1(t,\bzero,\tblambda) \\
                                          \vdots \\
                                          v^*_d(t,\bzero,\tblambda) \\
                                        \end{bmatrix}
                               \right\rangle
                   - \int_0^t \psi(v^*_1(s,\bzero,\tblambda),\ldots,v^*_d(s,\bzero,\tblambda)) \,\dd s \right\}
 \end{align*}
 for \ $\bx\in\RR_+^d$, \ $\tblambda\in\RR_+^d$ \ and \ $t\in\RR_+$, \ where,
 for any \ $(\blambda,\tblambda)\in\RR_+^{2d}$, \ the continuously differentiable function
 \[
   \RR_+\ni t\mapsto \bv^*(t,\blambda,\tblambda)
               := \big( v^*_1(t,\blambda,\tblambda),\ldots,v^*_{2d}(t,\blambda,\tblambda)
               \in\RR_+^{2d}
 \]
 is the unique locally bounded solution to the system of differential equations
 \[
   \partial_1 v^*_i(t, \blambda, \tblambda) = - \varphi_i^*(\bv^*(t, \blambda,\tblambda)), \qquad i\in\{1,\ldots,2d\},
 \]
 with initial conditions
 \begin{align*}
    v_i^*(0,\blambda, \tblambda) = \lambda_i , \qquad
    v_{d+i}^*(0,\blambda, \tblambda) = \tlambda_i , \qquad i \in \{1, \ldots, d\},
 \end{align*}
 where \ $\varphi_i^*$, \ $i=1,\ldots,2d$, \ are defined in \eqref{2dCBI_branching}.
Note that, for all \ $(\blambda,\tblambda)\in\RR_+^{2d}$ \ and \ $i\in\{1,\ldots,d\}$, \ we have that
 \[
   \partial_1 v^*_{d+i}(t, \blambda, \tblambda) = 0 \qquad \text{with \ $v_{d+i}^*(0,\blambda, \tblambda) = \tlambda_i$.}
 \]
Hence, for all \ $(\blambda,\tblambda)\in\RR_+^{2d}$ \ and \ $i\in\{1,\ldots,d\}$, \ we get
 \ $v^*_{d+i}(t, \blambda, \tblambda) = \tlambda_i$, \ $t\in\RR_+$.
\ Consequently, for all \ $(\blambda,\tblambda)\in\RR_+^{2d}$ \ and \ $i\in\{1,\ldots,d\}$, we have that
 \begin{align*}
 \partial_1 v^*_i(t, \blambda, \tblambda)
    = v^*_{d+i}(t, \blambda, \tblambda) - \varphi_i(\bv^*(t,\blambda, \tblambda))
    =\tlambda_i - \varphi_i(\bv^*(t,\blambda, \tblambda))
 \end{align*}
 with $v_i^*(0,\blambda, \tblambda) = \lambda_i$.
Hence, for all \ $\tblambda\in\RR_+^d$ \ and \ $i\in\{1,\ldots,d\}$, \ we get that
 \[
    \partial_1 v^*_i(t, \bzero, \tblambda)
       =\tlambda_i - \varphi_i(\bv^*(t,\bzero, \tblambda))
      \qquad \text{with \ $v_i^*(0,\bzero, \tblambda) = 0$.}
 \]
Since the system of differential equations \eqref{help1} has a unique locally bounded solution (see Li \cite[page 48]{Li}
 or Duffie et al.\ \cite[Proposition 6.4]{DufFilSch}),
 we have that \ $v_i^*(t,\bzero,\tblambda)=\tv_i(t,\tblambda)$, \ $t\in\RR_+$, \ $\tblambda\in\RR_+^d$ \ and \ $i\in\{1,\ldots,d\}$.
\ This yields the statement.
\proofend

Next, we specialize Proposition \ref{Pro_int_CBI_Laplace} to multi-type CB processes (see Definition \ref{Def_CBI}).

\begin{Cor}\label{Cor_int_CBI_Laplace}
Let \ $(\bZ_t)_{t\in\RR_+}$ \ be a multi-type CB process with parameters
 \ $(d, \bc, \bzero, \bB, 0, \bmu)$.
Then
 \begin{align}\label{help10}
  \EE_\bz\left( \exp\left\{ - \left\langle \tblambda, \int_0^t \bZ_u \,\dd u  \right\rangle \right\}  \right)
     = \ee^{-\langle \bz,  \tbv(t,\tblambda)\rangle}
 \end{align}
 for \ $t\in\RR_+$, \ $\tblambda\in \RR_+^d$ \ and \ $\bz\in\RR_+^d$, \ where
 \ $\RR_+\ni t\mapsto \tbv(t,\tblambda) = ( \tv_1(t,\tblambda),\ldots, \tv_d(t,\tblambda) )^\top\in\RR_+^d$ \
 is the unique locally bounded solution to the system of differential equations \eqref{help1}.
\end{Cor}

Now, we recall a result on the existence of an inverse of the branching mechanism \ $\bvarphi$ \
 of a multi-type CBI process with parameters \ $(d, \bc, \Bbeta, \bB, \nu, \bmu)$ \
 due to Chaumont and Marolleau \cite[Theorem 2.1]{ChaMar2}.
Let
 \[
   D_{\bvarphi}:=\big\{ \blambda\in\RR_+^d : \varphi_j(\blambda) >0, \;\; j=1,\ldots,d \big\}.
 \]
By part (2) of Theorem 2.1 in Chaumont and Marolleau \cite{ChaMar2}, if \ $D_{\bvarphi}$ \ is not empty,
 then there exists a mapping
 \ $\bphi=(\phi_1,\ldots,\phi_d):\RR_+^d \to \RR_+^d$ \ such that
 \ $\bphi(\blambda)\in D_{\bvarphi}$ \ for all \ $\blambda\in(0,\infty)^d$,
 \ and the restriction \ $\bphi\vert_{(0,\infty)^d} : (0,\infty)^d \to D_{\bvarphi}$ \ of the mapping \ $\bphi$ \ onto \ $(0,\infty)^d$ \
 is a diffeomorphism (a differentiable bijection with a differentiable inverse) such that its inverse is \ $\bvarphi\vert_{D_{\bvarphi}}:D_{\bvarphi}\to(0,\infty)^d$ \ satisfying
 \begin{align}\label{varphi_phi}
  \bvarphi(\bphi(\blambda)) = \blambda, \qquad \blambda\in(0,\infty)^d,
    \qquad \text{and}\qquad   \bphi(\bvarphi(\blambda)) = \blambda,\qquad \blambda\in D_{\bvarphi}.
 \end{align}
For simplicity, we will call \ $\bphi$ \ the inverse of the branching mechanism \ $\bvarphi$.
\ Here, for completeness, we note that one can indeed apply Theorem 2.1 in Chaumont and Marolleau \cite{ChaMar2}, since
  \ $\varphi_i$, \ $i=1,\ldots,d$ \ (given in Theorem \ref{CBI_exists}), can be written in the form of (5) in Chaumont and Marolleau \cite{ChaMar2},
 and the moment condition \eqref{help_moment_mu} implies the moment condition after formula (5) on page 3 in Chaumont and Marolleau \cite{ChaMar2}.
To give more details, for all \ $\blambda\in\RR_+^d$ \ and \ $i\in\{1,\ldots,d\}$, \ using the moment condition \eqref{help_moment_mu},
 we have
 \begin{align*}
  \varphi_i(\blambda)
   & = c_i \lambda_i^2 -  \langle \bB \be_i^{(d)}, \blambda \rangle
        + \int_{\cU_d} \Big( \lambda_i (1 \land z_i) - \langle \blambda,\bz\rangle \bone_{\{\Vert \bz\Vert <1\}} \Big) \, \mu_i(\dd \bz) \\
   &\phantom{=\;}
        + \int_{\cU_d}
          \bigl( \ee^{- \langle \blambda, \bz \rangle} - 1
                + \langle \blambda,\bz\rangle \bone_{\{\Vert \bz\Vert <1\}} \bigr)
          \, \mu_i(\dd \bz),
 \end{align*}
 since
  \begin{align*}
   & \int_{\cU_d} \Big\vert  \lambda_i (1 \land z_i) - \langle \blambda,\bz\rangle \bone_{\{\Vert \bz\Vert <1\}} \Big\vert \, \mu_i(\dd \bz)  \\
   & = \int_{\cU_d} \Big\vert  \lambda_i (1 \land z_i) \bone_{\{ \Vert \bz\Vert \geq 1\}}
           - \sum_{j\in\{1,\ldots,d\}\setminus\{i\}} \lambda_j (1\wedge z_j) \bone_{\{\Vert \bz\Vert <1\}} \Big\vert \, \mu_i(\dd \bz) \\
   & \leq \lambda_i \int_{\cU_d}  \Vert \bz\Vert  \bone_{\{ \Vert \bz\Vert \geq 1\}} \, \mu_i(\dd \bz)
           + \sum_{j\in\{1,\ldots,d\}\setminus\{i\}} \lambda_j \int_{\cU_d}  (1\wedge z_j) \bone_{\{\Vert \bz\Vert <1\}} \, \mu_i(\dd \bz)
   <\infty.
  \end{align*}
Hence for all \ $\blambda\in\RR_+^d$ \ and \ $i\in\{1,\ldots,d\}$, \ we get
 \begin{align}\label{help53}
  \begin{split}
  \varphi_i(\blambda) & = c_i \lambda_i^2 - \langle \bB \be_i^{(d)}, \blambda \rangle
        + \int_{\cU_d}
          \bigl( \ee^{- \langle \blambda, \bz \rangle} - 1
                + \langle \blambda,\bz\rangle \bone_{\{\Vert \bz\Vert <1\}} \bigr)
          \, \mu_i(\dd \bz)  \\
  &\phantom{=\;}
        + \int_{\cU_d} \Big( \lambda_i (1 \land z_i) \bone_{\{\Vert \bz\Vert \geq 1\}}
                              - \sum_{j\in\{1,\ldots,d\}\setminus\{i\}} \lambda_j (1\wedge z_j) \bone_{\{\Vert \bz\Vert <1\}} \Big) \, \mu_i(\dd \bz) \\
 & =  c_i \lambda_i^2 - \langle (\bB + \bS ) \be_i^{(d)}, \blambda \rangle
        + \int_{\cU_d}
          \bigl( \ee^{- \langle \blambda, \bz \rangle} - 1
                + \langle \blambda,\bz\rangle \bone_{\{\Vert \bz\Vert <1\}} \bigr)
          \, \mu_i(\dd \bz) ,
 \end{split}
 \end{align}
 where \ $\bS=(S_{j,i})_{j,i=1}^d$ \ is a $d\times d$ matrix having entries
  \begin{align*}
     S_{j,i}
        :=\begin{cases}
              \int_{\cU_d}(1\wedge z_j) \bone_{\{\Vert \bz\Vert <1\}}\,\mu_i(\dd\bz)
                 & \text{if \ $j\ne i$,}\\[1mm]
              -\int_{\cU_d}(1\wedge z_i) \bone_{\{\Vert \bz\Vert \geq 1\}}\,\mu_i(\dd\bz)
                  & \text{if \ $j= i$,}
          \end{cases}
          \qquad j,i\in\{1,\ldots,d\}.
  \end{align*}
Here \ $\bB, \, \bS\in\RR^{d \times d}_{(+)}$, \ yielding that \ $\bB + \bS\in\RR^{d \times d}_{(+)}$, \
 and hence, for each $i\in\{1,\ldots,d\}$, \ $\varphi_i$ \ has the form of (5) in Chaumont and Marolleau \cite{ChaMar2}, as desired.

The next proposition can be considered as a multi-type counterpart of Proposition 2.2 in He and Li \cite{HeLi}.

\begin{Pro}\label{Pro_tv}
Let \ $(d, \bc, \Bbeta, \bB, \nu, \bmu)$ \ be a set of admissible parameters.
 For all \ $\tblambda\in \RR_+^d$, \ let us consider
 the continuously differentiable function
    \begin{align}\label{tv_function}
     \RR_+\ni t \mapsto \tbv(t,\tblambda)= ( \tv_1(t,\tblambda),\ldots, \tv_d(t,\tblambda) )^\top\in\RR_+^d,
    \end{align}
    which is the unique locally bounded solution to the system of differential equations \eqref{help1}.
 \begin{enumerate}
  \item[(i)] For all \ $\tblambda\in \RR_+^d$, the function \eqref{tv_function} is increasing, and, consequently, the limit
   \[
      \tbv(\infty,\tblambda) = (\tv_1(\infty,\tblambda),\ldots,\tv_d(\infty,\tblambda))^\top := \lim_{t\to\infty} \tbv(t,\tblambda)\in[0,\infty]^d \qquad\text{exists}.
   \]
 \item[(ii)] Let \ $(\bZ_t)_{t\in\RR_+}$ \ be a multi-type CB process with parameters \ $(d, \bc, \bzero, \bB, 0, \bmu)$.
      If \ $\tblambda\in(0,\infty)^d$ \ and \ $i\in\{1,\ldots,d\}$ \ are such that \ $\tv_i(\infty,\tblambda)=\infty$, \ then
             \begin{align*}
               \PP_{\be^{(d)}_i}\left( \exists \; j\in\{1,\ldots,d\} : \int_0^\infty Z_{u,j} \,\dd u = \infty \right)
                =\PP_{\be^{(d)}_i}\left( \int_0^\infty \bZ_u\,\dd u \notin \RR_+^d \right)=1.
             \end{align*}
 \item[(iii)] If, in addition, \ $\tbB\in\RR_{(+)}^{d\times d}$ \ (given in \eqref{help5}) is irreducible, then, for all \ $\tblambda\in\RR_+^d$ \ with \ $\tblambda\ne\bzero$,
              \ we get \ $\tbv(\infty,\tblambda)\in (0,\infty]^d$.
 \item[(iv)]
  If, in addition, \ $D_\bvarphi\ne\emptyset$, \ then for all \ $\tblambda\in(0,\infty)^d$ \ satisfying \ $\tbv(\infty,\tblambda)\in(0,\infty)^d$, \ we have that
    \ $\tbv(\infty,\tblambda)=\bphi(\tblambda)$, \ where the diffeomorphism \ $\bphi\vert_{(0,\infty)^d} : (0,\infty)^d \to D_{\bvarphi}$ \
   is the inverse of \ $\bvarphi\vert_{D_\bvarphi}:D_\bvarphi\to(0,\infty)^d$ \ (see \eqref{varphi_phi}).
 \end{enumerate}
\end{Pro}

\noindent{\bf Proof.}
(i):
First, note that the function \ $\tbv$ \ depends only on the parameters \ $\bc$, \ $\bB$ \ and \ $\bmu$,
 \ but not on \ $\Bbeta$ \ and \ $\nu$, \ since the functions \ $\varphi_1,\ldots,\varphi_d$ \ appearing
 in the system of differential equations \eqref{help1} depend only on \ $\bc$, \ $\bB$ \ and \ $\bmu$.
\ Hence, in order to prove part (i), we may and do assume that  \ $\Bbeta=\bzero$ \ and \ $\nu=0$.
\ Let \ $(\bZ_t)_{t\in\RR_+}$ \ be a multi-type CB process with parameters \ $(d, \bc, \bzero, \bB, 0, \bmu)$.
For all $\bz\in\RR_+^d$ and \ $\tblambda\in\RR_+^d$, the left hand side of \eqref{help10} as a function of $t$ is decreasing
 (since $\bZ_u\in\RR_+^d$, $u\in\RR_+$), and hence the same holds for the right hand side of \eqref{help10} as well.
Consequently, for all $\bz\in\RR_+^d$ and \ $\tblambda\in\RR_+^d$,
 the mapping $\RR_+\ni t \mapsto \langle \bz,  \tbv(t,\tblambda)\rangle$ is increasing.
By choosing \ $\bz=\be^{(d)}_i$, $i=1,\ldots,d$, we have that for all \ $\tblambda\in\RR_+^d$,
 \ the mappings $\RR_+\ni t \mapsto \tv_i(t,\tblambda)$, $i=1,\ldots,d$, are increasing.
This implies that the function \eqref{tv_function} is increasing for all \ $\tblambda\in\RR_+^d$, as desired.
Consequently, for all $\tblambda\in\RR_+^d$,  we have that \ $\partial_1\tv_i(t,\tblambda)\geq 0$, $t\in\RR_+$, \ $i = 1,\ldots,d$, \ and the limit
 \[
    \tbv(\infty,\tblambda):=\lim_{t\to\infty}\tbv(t,\tblambda)\in[0,\infty]^d \qquad \text{exists.}
 \]

(ii):
Let \ $\tblambda\in(0,\infty)^d$ \ and \ $i\in\{1,\ldots,d\}$ \ be such that \ $\tv_i(\infty,\tblambda)=\infty$.
\ Then we have \ $\ee^{-\langle \be^{(d)}_i,  \tbv(t,\tblambda)\rangle}  = \ee^{-\tv_i(t,\tblambda)}\to 0$ \ as \ $t\to\infty$.
\ Hence, using \eqref{help10} and the dominated convergence theorem, we have that
 \begin{align*}
   0 & = \lim_{t\to\infty} \EE_{\be^{(d)}_i}\left( \exp\left\{ - \left\langle \tblambda, \int_0^t \bZ_u \,\dd u  \right\rangle \right\}  \right) \\
     & = \EE_{\be^{(d)}_i}\left( \exp\left\{ - \left\langle \tblambda, \int_0^\infty \bZ_u \,\dd u  \right\rangle \right\}
                                   \bone_{\left\{ \int_0^\infty \bZ_u \,\dd u \in \RR_+^d\right\}} \right),
 \end{align*}
 yielding that
 \[
  \PP_{\be^{(d)}_i}\left( \exp\left\{ - \left\langle \tblambda, \int_0^\infty \bZ_u \,\dd u  \right\rangle \right\}
                                   \bone_{\left\{  \int_0^\infty \bZ_u \,\dd u  \in\RR_+^d\right\}} = 0 \right)=1.
 \]
Consequently, we have
 \[
  \PP_{\be^{(d)}_i}\left( \bone_{\left\{  \int_0^\infty \bZ_u \,\dd u  \in\RR_+^d\right\}} = 0  \right)=1,
 \]
 which implies (ii).

(iii):
Similarly as it was explained at the beginning of the proof of part (i),
 we may and do assume that  \ $\Bbeta=\bzero$ \ and \ $\nu=0$.
\ Let \ $(\bZ_t)_{t\in\RR_+}$ \ be a multi-type CB process with parameters \ $(d, \bc, \bzero, \bB, 0, \bmu)$.
Since \ $\tbB$ \ is irreducible, by Definition \ref{Def_irreducible}, we have \ $(\bZ_t)_{t\in\RR_+}$ \ is irreducible.
Let \ $\tblambda=(\tlambda_1,\ldots,\tlambda_d)\in\RR_+^d$ \ be such that \ $\tblambda\ne \bzero$.
\ We need to check that
 \begin{align}\label{help15}
   \tv_j(\infty,\tblambda)>0 \qquad \text{for each \ $j\in\{1,\ldots,d\}$.}
 \end{align}
On the contrary, let us assume that there exists \ $i\in\{1,\ldots,d\}$ \ such that \ $\tv_i(\infty,\tblambda)=0$.
\ Then, since the function \ $\RR_+\ni t\mapsto \tv_i(t,\tblambda)$ \ is non-negative and increasing (see part (i)),
 we get \ $\tv_i(t,\tblambda)=0$, \ $t\in\RR_+$.
\ Using \eqref{help10} with the choice \ $\bz:=\be^{(d)}_i$, \ it implies that
 \[
  \EE_{\be^{(d)}_i}\left( \exp\left\{ - \left\langle \tblambda, \int_0^t \bZ_u \,\dd u  \right\rangle \right\} \right)
       = 1, \qquad t\in\RR_+.
 \]
Consequently, we have that
 \[
   \PP_{\be^{(d)}_i}\left(  \left\langle \tblambda, \int_0^t \bZ_u \,\dd u  \right\rangle =0 \right)=1, \qquad t\in\RR_+,
 \]
 i.e.,
 \begin{align}\label{help51}
   \PP_{\be^{(d)}_i}\left(  \int_0^t \langle \tblambda, \bZ_u \rangle \,\dd u = 0 \right)=1, \qquad t\in\RR_+.
 \end{align}
Next, we check that \eqref{help51} yields that
 \begin{align}\label{help14}
     \PP_{\be^{(d)}_i}\Big( \tlambda_j Z^{(j)}_t  = 0  \Big)=1, \qquad t\in\RR_+, \;\; j\in\{1,\ldots,d\}.
 \end{align}
It will easily follow from the following auxiliary lemma.

{\sl Auxiliary lemma:} If $t>0$ and $f:[0,t]\to \RR$ is a c\`{a}dl\`{a}g function such that $f(u)=0$ Lebesgue a.e.\ $u\in[0,t]$,
  then $f(u) =0$ for all $u\in[0,t)$.

 {\sl Proof of Auxiliary lemma:}
By the assumption, there exists a Lebesgue measurable set $S\subseteq [0,t]$ having Lebesgue measure $0$ such that
 $f(t)=0$ for $u\in[0,t]\setminus S$.
Let $u_0\in [0,t)$ be arbitrary.
Since $S$ has Lebesgue measure $0$, there is no right neighbourhood of $u_0$ contained in $S$.
Hence for all $\vare>0$, we have $[u_0,u_0+\vare)\nsubseteq S$, yielding that for each $n\in\NN$,
 there exists $u_n\in[0,t]\setminus S$ such that $u_n\in[u_0,u_0+\frac{1}{n})$.
Then $u_n\to u_0$ as $n\to\infty$, and $f(u_n)=0$, $n\in\NN$.
Since $f$ is right continuous, we get $f(u_n)\to f(u_0)$ as $n\to\infty$, yielding that $f(u_0)=0$, as desired.

We start to check \eqref{help14}.
Using that \ $\langle \tblambda, \bZ_u \rangle\in\RR_+$, $u\in\RR_+$, \
 if $\omega\in\Omega$ and $t\in\RR_+$ are such that $\int_0^t \langle \tblambda, \bZ_u(\omega) \rangle\,\dd u  = 0$, then
 $\langle \tblambda, \bZ_u(\omega) \rangle =0$ Lebesgue a.e.\ $u\in[0,t]$.
Since $(\bZ_u(\omega))_{u\in\RR_+}$ is c\`{a}dl\`{a}g, we have that
 $(\langle \tblambda, \bZ_u(\omega) \rangle)_{u\in\RR_+}$ is c\`{a}dl\`{a}g as well, and hence,
 by the Auxiliary lemma above, we get that
 $\langle \tblambda, \bZ_u(\omega) \rangle =0$ for all $u\in[0,t)$.
Using that \ $\langle \tblambda, \bZ_u(\omega)\rangle = 0$ \ holds if and only if \ $\tlambda_jZ^{(j)}_u(\omega)=0$, \ $j\in\{1,\ldots,d\}$,
 \ we get that if $\omega\in\Omega$ and $t\in\RR_+$ are such that $\int_0^t \langle \tblambda, \bZ_u(\omega) \rangle\,\dd u  = 0$, then
 $\tlambda_jZ^{(j)}_u(\omega)=0$ for all \ $u\in[0,t)$, $j\in\{1,\ldots,d\}$.
Taking into account \eqref{help51}, it implies \eqref{help14}.

In particular, \eqref{help14} with \ $j:=i$ \ and \ $t:=0$ \ yields that
 \ $\PP_{\be^{(d)}_i}( \tlambda_i Z^{(i)}_0 = 0)=1$,
 \ and hence \ $\tlambda_i=0$.
\ Further, \eqref{help14} also yields that
 \begin{align*}
     \PP_{\be^{(d)}_i}\Big( Z^{(j)}_t  = 0 \Big)=1, \qquad t\in\RR_+,
 \end{align*}
 for each \ $j\in\{1,\ldots,d\}\setminus\{i\}$ \ with \ $\tlambda_j>0$.
\ In particular, it implies that \ $\EE_{\be^{(d)}_i}( Z^{(j)}_1)=0$ \ for each \ $j\in\{1,\ldots,d\}\setminus\{i\}$ \ with \ $\tlambda_j>0$.
By \eqref{EXcond}, we have that \ $(\be_j^{(d)})^\top \ee^\tbB \be_i^{(d)} = \EE_{\be_i^{(d)}}( Z^{(j)}_1)$, \
 and hence \ $(\ee^\tbB)_{j,i}=0$ \ for each \ $j\in\{1,\ldots,d\}\setminus\{i\}$ \ with \ $\tlambda_j>0$.
\ Since \ $\tbB$ \ is irreducible, we have \ $\ee^{\tbB} \in \RR^{d \times d}_{++}$ \ (see, e.g., Barczy and Pap \cite[Lemma A.1]{BarPap}),
 which yields that there does not exist \ $j\in\{1,\ldots,d\}\setminus\{i\}$ \ with \ $\tlambda_j>0$.
Taking into account that \ $\tlambda_i=0$, \ it implies that \ $\tblambda=\bzero$, \ leading us to a contradiction.
That is, we get \eqref{help15}, as desired.

(iv):
Let us suppose that \ $D_\varphi\ne\emptyset$, \ and let $\tblambda\in (0,\infty)^d$ be such that $\tbv(\infty,\tblambda)\in(0,\infty)^d$.
\ Let \ $i\in\{1,\ldots,d\}$ \ be fixed.
Recall that  $\varphi_i$ is continuous, which can be proved using the dominated convergence theorem taking into account \eqref{help53},
 the moment assumption \eqref{help_moment_mu} and the inequalities
 \ $\vert \ee^{-x} - 1 \vert\leq x$, $x\in\RR_+$, \ and \ $0\leq \ee^{-x} - 1 + x \leq \frac{x^2}{2}$, $x\in\RR_+$.
Hence we have that the limit $\lim_{t\to\infty} \varphi_i(\tbv(t,\tblambda))=\varphi_i(\tbv(\infty,\tblambda))$ exists.
Further, using \eqref{help1} and that \ $\partial_1\tv_i(t,\tblambda)\geq 0$, $t\in\RR_+$ \ (following from part (i)),  we get
 \[
   \varphi_i(\tbv(t,\tblambda)) = \tlambda_i - \partial_1 \tv_i(t,\tblambda) \leq \tlambda_i, \qquad t\in\RR_+,
 \]
 which yields that $\varphi_i(\tbv(\infty,\tblambda))\in(-\infty,\tlambda_i]$.
Applying again \eqref{help1}, we have
 \begin{align}\label{help52}
   \lim_{t\to\infty} \partial_1 \tv_i(t,\tblambda)
     = \tlambda_i - \lim_{t\to\infty} \varphi_i(\tbv(t,\tblambda))
     = \tlambda_i - \varphi_i(\tbv(\infty,\tblambda)) \in[0,\infty).
 \end{align}
Using that $\RR_+\ni t\mapsto \tv_i(t,\tblambda)$ is increasing (see part (i)), $\tv_i(\infty,\tblambda)\in(0,\infty)$
 and $\lim_{t\to\infty} \partial_1 \tv_i(t,\tblambda)$ exists (it belongs to $[0,\infty)$),
 an elementary calculus shows that $\lim_{t\to\infty} \partial_1 \tv_i(t,\tblambda)=0$.
Indeed, for all $t\in\RR_+$, by the mean value theorem, there exists \ $\xi_t\in[t,t+1]$ \ such that
 \[
     \tv_i(t+1,\tblambda)  - \tv_i(t,\tblambda)  =   \partial_1 \tv_i(\xi_t,\tblambda).
 \]
By taking the limit of both sides as $t\to\infty$, we get
 \[
    \tv_i(\infty,\tblambda) - \tv_i(\infty,\tblambda) = \lim_{t\to\infty} \partial_1 \tv_i(t,\tblambda),
 \]
 implying that $\lim_{t\to\infty} \partial_1 \tv_i(t,\tblambda)=0$, as desired.
Hence, by \eqref{help52}, for all \ $\tblambda\in(0,\infty)^d$ \ satisfying \ $\tbv(\infty,\tblambda)\in(0,\infty)^d$, we have
 \[
  \varphi_i(\tbv(\infty,\tblambda)) = \tlambda_i, \qquad i=1,\ldots,d,
   \qquad \text{that is,} \qquad \bvarphi(\tbv(\infty,\tblambda)) = \tblambda.
 \]
Using that the inverse of the diffeomorphism \ $\bphi\vert_{(0,\infty)^d} : (0,\infty)^d \to D_{\bvarphi}$ \
 is \ $\bvarphi\vert_{D_\bvarphi}:D_\bvarphi\to(0,\infty)^d$ \ (see \eqref{varphi_phi}), we have
 $\bvarphi(\bphi(\tblambda))=\tblambda$, \ and hence
 \ $\tbv(\infty,\tblambda)=\bphi(\tblambda)$ \ for all \ $\tblambda\in (0,\infty)^d$ \ satisfying
 \ $\tbv(\infty,\tblambda)\in(0,\infty)^d$, as desired.
\proofend

Note that in the proof of part (iii) of Proposition \ref{Pro_tv}, we need that
 \ $\ee^{\tbB} \in \RR^{d \times d}_{++}$, which holds under our assumption that $\tbB\in\RR_{(+)}^{d\times d}$ is irreducible.
Further, we call the attention that the irreducibility of $\tbB\in\RR_{(+)}^{d\times d}$ is not only sufficient,
 but also necessary in order that $\ee^{\tbB} \in \RR^{d \times d}_{++}$ hold, see Lemma A.1 in Barczy and Pap \cite{BarPap}.

The following Corollary \ref{Cor_inf_int_CBI} can be considered as a multi-type counterpart of Corollary 2.1 in He and Li \cite{HeLi},
 which is also contained as Corollary 5.21 in Li \cite{Li}.

\begin{Cor}\label{Cor_inf_int_CBI}
Let \ $(\bZ_t)_{t\in\RR_+}$ \ be an irreducible multi-type CB process with parameters
 \ $(d, \bc, \bzero, \bB, 0, \bmu)$ \ such that \ $D_\bvarphi\ne\emptyset$.
Then
 \begin{align*}
  \EE_{\bz}\left( \exp\left\{ - \left\langle \tblambda, \int_0^\infty \bZ_u \,\dd u  \right\rangle \right\}
              \bone_{\left\{ \int_0^\infty \bZ_u \,\dd u \in\RR_+^d\right\}}  \right)
     = \ee^{-\langle \bz,  \tbv(\infty,\tblambda) \rangle}
     = \ee^{-\langle \bz,  \bphi(\tblambda) \rangle}
 \end{align*}
 for all \ $\bz\in\RR_+^d$ \ and \ $\tblambda\in (0,\infty)^d$ \ satisfying \ $v_i(\infty, \tblambda)<\infty$, \ $i\in\{1,\ldots,d\}$.
\end{Cor}

\noindent{\bf Proof.}
It follows from dominated convergence theorem, Corollary \ref{Cor_int_CBI_Laplace} and Proposition \ref{Pro_tv}.
\proofend

In the next remark, we discuss the connections between Proposition 1 in Chaumont and Marolleau \cite{ChaMar2} and our Corollary \ref{Cor_inf_int_CBI}.

\begin{Rem}
We observe that Corollary \ref{Cor_inf_int_CBI} can also be derived by the theory of spectrally positive additive L\'evy fields.
Let \ $(\bZ_t)_{t\in\RR_+}$ \ be an irreducible multi-type CB process with parameters
 \ $(d, \bc, \bzero, \bB, 0, \bmu)$ \ such that $D_\bvarphi\ne\emptyset$.
According to Proposition 2.1 in Chaumont and Marolleau \cite{ChaMar2},  if $\bZ_0 =\bz\in\RR_+^d$,
 then \ $\int_0^\infty \bZ_u \,\dd u=\bT_{\bz}$ \ holds \ $\PP$-almost surely,
 where \ $\bT_{\bz}$ \ is the (multivariate) first hitting time of the level $-\bz$
 by the spectrally positive additive L\'evy field corresponding to the branching mechanism of $(\bZ_t)_{t\in\RR_+}$,
 see Proposition 1 in Chaumont and Marolleau \cite{ChaMar2}.
Then using Theorem 2.1 in Chaumont and Marolleau \cite{ChaMar2} and the convention $\ee^{-\infty}:=0$,
 we get Corollary \ref{Cor_inf_int_CBI}.
\proofend
\end{Rem}

\section{Distributional properties of jump times}
\label{Section_dist_jump}

Recall that \ $\cB(\cU_d)$ \ denotes the set of Borel subsets of \ $\cU_d$.
\ For all \ $t\in\RR_{++}$ \ and \ $A \in\cB(\cU_d)$, \ let
\begin{align*}
  \tau_A
    := \inf\{u \in \RR_{++}
              : \Delta \bX_u \in A\},
  \qquad
  J_t(A)
   := \card(\{u \in (0, t] : \Delta \bX_u  \in A\}),
 \end{align*}
 with the convention \ $\inf(\emptyset) := \infty$, \  where
 \ $\Delta \bX_u := \bX_u - \bX_{u-}$, \ $u \in \RR_{++}$, \
 and \ $\card(H)$ \ denotes the cardinality of a set \ $H$.

In the forthcoming results, given a set \ $A\in\cB(\cU_d)$, \ the condition that \ $\nu(A)+\sum_{\ell=1}^d \mu_\ell(A)<\infty$
 \ will come into play.
In the next remark, we give a sufficient condition under which it holds.

\begin{Rem}\label{Rem_A_condition}
First, recall the notation  $K_\vare = \{\by\in\RR_+^d : \Vert \by\Vert<\vare\}$.
If \ $A\in\cB(\cU_d)$ \ is such that there exists \ $\vare\in(0,1)$ \ with \ $A\subseteq \RR_+^d\setminus K_\vare$
 \ (or equivalently, \ $A\in\cB(\cU_d)$ \ is such that \ $\bzero$ \ is not contained in the closure of $A$), \ then
 \ $\nu(A)+\sum_{\ell=1}^d \mu_\ell(A)<\infty$.
\ Indeed, if \ $A\in \cB(\cU_d)$ \ is such that \ $A\subseteq \RR_+^d\setminus K_\vare$ \ with some $\vare\in(0,1)$,
 \ then by parts (v) and (vi) of Definition \ref{Def_admissible}, we have
 \begin{align*}
    \nu(A) \leq \nu(\RR_+^d\setminus K_\vare)
           \leq \frac{1}{\vare} \int_{\{ \br\in\cU_d\,:\, \vare\leq \Vert \br\Vert \leq 1 \}} \Vert\br\Vert \,\nu(\dd\br)
             + \int_{\{ \br\in\cU_d\,:\, \Vert \br\Vert >1 \}}1 \,\nu(\dd\br)
           <\infty,
 \end{align*}
 and, for each \ $\ell\in\{1,\ldots,d\}$,
 \begin{align*}
    \mu_\ell(A) \leq \mu_\ell(\RR_+^d\setminus K_\vare)
           \leq \frac{1}{\vare^2} \int_{\{ \bz\in\cU_d\,:\, \vare\leq \Vert \bz\Vert \leq 1 \}} \Vert\bz\Vert^2 \,\mu_\ell(\dd\bz)
             + \int_{\{ \bz\in\cU_d\,:\, \Vert \bz\Vert >1 \}} \Vert \bz\Vert \,\mu_\ell(\dd\bz)
           <\infty,
 \end{align*}
 as desired.
\proofend
\end{Rem}

The next result is a generalization of Proposition 3.1 and Theorem 3.1 in He and Li \cite{HeLi}.
 We mention that Theorem 3.1 in He and Li \cite{HeLi} is also contained as Theorem 10.13 in Li \cite{Li}.
For completeness, we note that He and Li \cite{HeLi} do not assume the moment condition \eqref{moment_condition_m_new}
 for deriving their results, however we will assume this condition in order to derive our forthcoming results
 in Theorem \ref{Thm_jump_times_dCBI}.
The moment condition \eqref{moment_condition_m_new} is needed for our approach
 in order to be able to use Theorem 4.6 in Barczy et al.\ \cite{BarLiPap2} about the SDE representations of multi-type CBI processes,
 and in order to ensure the finiteness of the expectation of the norm of a multi-type CBI process at a given time point.
Further, as we already mentioned in the Introduction, the vector $\tBbeta$ belongs $\RR_+^d$ if and only if the moment condition
 \eqref{moment_condition_m_new} holds.
Under the condition \eqref{moment_condition_m_new}, we can prove not only that $J_t(A)$ is finite almost surely
 for all \ $t\in\RR_+$, \ but also that
 $J_t(A)$ has a finite expectation, see part (i) of Theorem \ref{Thm_jump_times_dCBI}.
We also mention that part (iii) of our Theorem \ref{Thm_jump_times_dCBI} has just appeared as Example 12.2 in
 the new second edition of Li's book \cite{Li}.
The two research works have been carried out parallelly, so we decided to present our result and its proof as well.
Concerning the proofs of Example 12.2 in Li \cite{Li} and our Theorem 4.2,
 Li \cite{Li} derived his Example 12.2 as a special case of Theorem 12.22 in Li \cite{Li}, which
 is a measure-valued generalization of Theorem 3.1 in He and Li \cite{HeLi}.
The proof of Theorem 12.22 in Li \cite{Li} is based on the fact that
 a c\`{a}dl\`{a}g immigration superprocess can be represented as a pathwise unique
 c\`{a}dl\`{a}g strong solution of a stochastic integral equation (see Theorem 12.14 in Li \cite{Li}).
In the proof of our Theorem \ref{Thm_jump_times_dCBI}, we do not apply Theorems 12.14 and 12.22 in Li \cite{Li},
 instead, we directly use the SDE \eqref{SDE_atirasa_dimd} for a multi-type CBI process,
 which may be more easily accessible than that of immigration superprocess.

Given a set of admissible parameters \ $(d, \bc, \Bbeta, \bB, \nu, \bmu)$ \ and a Borel set
 $A\in\cB(\cU_d)$ such that $\nu(A)+\sum_{\ell=1}^d \mu_\ell(A)<\infty$, let us introduce
 \ $\bB^{(A)}\in\RR^{d\times d}_{(+)}$, \ $\nu^{(\cU_d\setminus A)}$, \ and $\bmu^{(\cU_d\setminus A)}$ \ by
  \begin{align}\label{mod_parameters}
   \begin{split}
             &(\bB^{(A)})_{i,j}
                 := b_{i,j} + \delta_{i,j}\int_A \big( (z_i - 1)^+ - z_i \big)\,\mu_i(\dd\bz),\qquad i,j\in\{1,\ldots,d\},\\
             &\nu^{(\cU_d\setminus A)}:=\nu\vert_{\cU_d\setminus A},\\
             &\bmu^{(\cU_d\setminus A)}:=\Big(\mu^{(\cU_d\setminus A)}_1,\ldots,\mu^{(\cU_d\setminus A)}_d\Big)  \quad \text{with}\quad
               \mu^{(\cU_d\setminus A)}_\ell:=\mu_\ell\vert_{\cU_d\setminus A}, \quad \ell\in\{1,\ldots,d\}.
   \end{split}
  \end{align}
Then one can easily see that \ $\big(d, \bc, \Bbeta, \bB^{(A)}, \nu^{(\cU_d\setminus A)}, \bmu^{(\cU_d\setminus A)}\big)$ \
 is a set of admissible parameters.
Indeed, for each $i\in\{1,\ldots,d\}$, we have $\int_A \vert (z_i - 1)^+ - z_i \vert \,\mu_i(\dd\bz)
 = \int_A  z_i \bbone_{\{z_i\leq1\}}\,\mu_i(\dd\bz) + \int_A  \bbone_{\{z_i>1\}}\,\mu_i(\dd\bz) \leq 2\mu_i(A)<\infty$.
Further, given a multi-type CBI process $(\bX_t)_{t\in\RR_+}$ with parameters \ $(d, \bc, \Bbeta, \bB, \nu, \bmu)$ \
 and a Borel set $A\in\cB(\cU_d)$ such that $\nu(A)+\sum_{\ell=1}^d \mu_\ell(A)<\infty$,
 let \ $(\bX^{(A)}_t)_{t\in\RR_+}$ \ be a multi-type CBI process with parameters
 \ $\big(d, \bc, \Bbeta, \bB^{(A)}, \nu^{(\cU_d\setminus A)}, \bmu^{(\cU_d\setminus A)}\big)$ \ such that \ $\bX_0^{(A)} = \bX_0$.
Intuitively, the process $(\bX^{(A)}_t)_{t\in\RR_+}$ is obtained by removing from $(\bX_t)_{t\in\RR_+}$ all the
 masses produced by the jumps with size vectors in the set $A$.
This argument is made precise in mathematical terms in the proof of our next Theorem \ref{Thm_jump_times_dCBI}.
Note that if the moment condition \eqref{moment_condition_m_new} holds for $(\bX_t)_{t\in\RR_+}$, then it also holds for
 $(\bX^{(A)}_t)_{t\in\RR_+}$, since
 \begin{align}\label{moment_condition_m_new_A}
  \int_{\cU_d} \|\br\| \bbone_{\{\|\br\|\geq1\}} \, \nu^{(\cU_d\setminus A)}(\dd\br)
     =  \int_{\cU_d\setminus A} \|\br\| \bbone_{\{\|\br\|\geq1\}} \, \nu(\dd\br) <\infty.
 \end{align}
For the branching and immigration mechanisms, and an SDE for $(\bX^{(A)}_t)_{t\in\RR_+}$, see
 Theorem \ref{Thm_jump_times_dCBI} and \eqref{SDE_atirasa_dimd_R}, respectively.

\begin{Thm}\label{Thm_jump_times_dCBI}
Let \ $(\bX_t)_{t\in\RR_+}$ \ be a multi-type CBI process with parameters \ $(d, \bc, \Bbeta, \bB, \nu, \bmu)$ \ such that
 \ $\EE(\|\bX_0\|) < \infty$ \ and the moment condition \eqref{moment_condition_m_new} hold.
Then for all \ $A\in\cB(\cU_d)$ \ such that
 \ $\nu(A)+\sum_{\ell=1}^d \mu_\ell(A)<\infty$, \ we have
  \begin{itemize}
    \item[(i)] $\EE(J_t(A))<\infty$, \ $t \in \RR_{++}$, \ which, in particular, implies that
                \ $\PP(J_t(A) < \infty) = 1$, \ $t\in\RR_{++}$;
    \item[(ii)] for all $t\in\RR_+$,
         \begin{align}\label{formula1}
           \PP(\tau_A>t \mid \bX_0 )
                =  \ee^{-\nu(A)t} \EE\left( \exp\left\{ -\sum_{\ell=1}^d \mu_\ell(A) \int_0^t X_{u,\ell}^{(A)}\,\dd u \right\} \,\bigg\vert\, \bX_0 \right),
         \end{align}
         where $(\bX^{(A)}_t)_{t\in\RR_+}$ is a multi-type CBI process with parameters
         $\big(d, \bc, \Bbeta, \bB^{(A)}, \nu^{(\cU_d\setminus A)}, \bmu^{(\cU_d\setminus A)}\big)$
         such that $\bX_0^{(A)} = \bX_0$,
         where \ $\bB^{(A)}$, \ $\nu^{(\cU_d\setminus A)}$, \ and $\bmu^{(\cU_d\setminus A)}$ \
         are given in \eqref{mod_parameters},
         and \ $X^{(A)}_{t,\ell}$ \ denotes the \ $\ell^{\mathrm th}$ \ coordinate of \ $\bX^{(A)}_t$
         \ for any \ $t\in\RR_+$ \ and \ $\ell \in \{1, \ldots, d\}$;
    \item[(iii)] for all \ $t\in\RR_+$ \ and \ $\bx=(x_1,\ldots,x_d)^\top\in\RR_+^d$,
       \begin{align}\label{formula3}
         \PP_{\bx}(\tau_A>t)
          = \exp\left\{ -\nu(A)t - \sum_{\ell=1}^d x_\ell\,\tv_\ell^{(A)}(t,\bmu(A))
                           - \int_0^t \psi^{(A)}\big(\tbv^{(A)}(s,\bmu(A))\big)\,\dd s \right\},
        \end{align}
        where
          \begin{itemize}
           \item[$\bullet$] $\bmu(A):=(\mu_1(A),\ldots,\mu_d(A))$,
           \item[$\bullet$] the function
                    \begin{align}\label{formula_psi_A}
                       \RR_+^d \ni \blambda \mapsto
                       \psi^{(A)}(\blambda)
                        : = \langle \bbeta, \blambda \rangle
                             + \int_{\cU_d\setminus A}
                             \bigl( 1 - \ee^{- \langle\blambda, \br\rangle} \bigr)
                             \, \nu(\dd\br)\in \RR_+
                      \end{align}
           is the immigration mechanism of \ $(\bX^{(A)}_t)_{t\in\RR_+}$,
           \item[$\bullet$] the continuously differentiable function
                 \begin{align}\label{tvA}
                    \RR_+\ni t \mapsto \tbv^{(A)}(t,\bmu(A)):=(\tv^{(A)}_1(t,\bmu(A)),\ldots,\tv^{(A)}_d(t,\bmu(A)))^\top\in\RR_+^d
                  \end{align}
                  is the unique locally bounded solution to the system of differential equations
                 \begin{align}\label{help_DE_tv}
                  \begin{split}
                    &\partial_1 \tv_i^{(A)}(t, \bmu(A)) = \mu_i(A) - \varphi_i^{(A)}( \tbv^{(A)}(t, \bmu(A)) ), \qquad i \in \{1, \ldots, d\} ,\\
                    & \tv_i^{(A)}(0, \bmu(A)) = 0 , \qquad i \in \{1, \ldots, d\},
                  \end{split}
                 \end{align}
                 where the function \ $\RR_+^d\ni \blambda \mapsto \bvarphi^{(A)}(\blambda):=(\varphi^{(A)}_1(\blambda),\ldots,\varphi^{(A)}_d(\blambda))^\top\in\RR^d$
                 \ with
                 \begin{align}\label{help54}
                  \begin{split}
                    \varphi^{(A)}_\ell(\blambda)
                   & := \varphi_\ell(\blambda) + \int_A  (1  - \ee^{-\langle \blambda, \bz \rangle}) \, \mu_\ell(\dd \bz)
                 \end{split}
              \end{align}
              for \ $\blambda\in\RR_+^d$, \ $\ell\in\{1,\ldots,d\}$,
              \ is the branching mechanism of \ $(\bX^{(A)}_t)_{t\in\RR_+}$.
          \end{itemize}
     \end{itemize}
\end{Thm}

Concerning the notations in Theorem \ref{Thm_jump_times_dCBI}, we note that \ $(A)$ \ in the superscript
 of a formula (e.g. in that of \ $\varphi^{(A)}_\ell$) \ means only that the corresponding expression
 depends on \ $A$ \ (however, it may depend on \ $\cU_d\setminus A$ \ as well).
Nonetheless, the restrictions of the measures \ $\nu$ \ and \ $\bmu$ \ onto \ $\cU_d\setminus A$
 \ are denoted by \ $\nu^{(\cU_d\setminus A)}$ \ and \ $\bmu^{(\cU_d\setminus A)}$, \ respectively,
 in order to avoid some possible confusion.
In the proof of Theorem \ref{Thm_jump_times_dCBI}, we also follow this
 convention when a (random) measure is restricted to a subset.

\noindent
\textbf{Proof of Theorem \ref{Thm_jump_times_dCBI}.}
First, note that, by \eqref{EXcond}, we have
 \[
     \Vert \EE_{\bx}(\bX_t) \Vert
         \leq \left\Vert \ee^{t\tbB} \right\Vert  \Vert \bx \Vert + \int_0^t \left\Vert \ee^{u\tbB}\right\Vert \Vert \tBbeta \Vert \, \dd u ,
           \qquad \bx \in \RR_+^d , \quad t \in \RR_+ ,
 \]
 yielding that
 \[
  \Vert \EE(\bX_t ) \Vert
    = \Vert \EE(\EE( \bX_t \mid \bX_0)) \Vert
    \leq \EE( \Vert \EE(\bX_t\mid \bX_0)\Vert)
    \leq  \left\Vert \ee^{t\tbB} \right\Vert \EE(\Vert \bX_0 \Vert)
          + \Vert \tBbeta \Vert \int_0^t \left\Vert \ee^{u\tbB}\right\Vert \, \dd u
 \]
 for \ $t\in\RR_+$.
\ Consequently,
 \begin{align}\label{help_int_moment}
   \int_0^t \Vert \EE( \bX_u ) \Vert \, \dd u
      \leq \EE(\Vert \bX_0 \Vert) \int_0^t \left\Vert \ee^{u\tbB}\right \Vert \,\dd u
           + \Vert \tBbeta \Vert \int_0^t \int_0^u \left\Vert \ee^{v\tbB}\right\Vert \,\dd v \, \dd u
      <\infty,  \qquad t\in\RR_+,
 \end{align}
 since the function \ $\RR_+\ni u \mapsto \left\Vert \ee^{u\tbB} \right\Vert$ \ is continuous.

Let \ $A\in\cB(\cU_d)$ \ be fixed such that \ $\nu(A)+\sum_{\ell=1}^d \mu_\ell(A)<\infty$.
\ The proof is divided into Steps 1-7.

{\sl Step 1.}
Roughly speaking, the jumps of \ $(\bX_t)_{t\in \RR_+}$ \ are associated with
 the jumps of the mutually independent Poisson point processes corresponding to the mutually independent Poisson random measures
 \ $N_\ell$, \ $\ell \in\{1, \ldots, d\}$, \ and \ $M$.
\ In what follows, we make it precise.
We can rewrite the SDE \eqref{SDE_atirasa_dimd} in the form
 \begin{align}\label{help3}
  \begin{split}
  \bX_t
  &=\bX_0 + \int_0^t (\Bbeta + \tbB \bX_u) \, \dd u
    - \sum_{\ell=1}^d
       \int_0^t \int_{\cU_d} \int_{\cU_1}
        \bz \bbone_A(\bz) \bbone_{\{w\leq X_{u,\ell}\}}
        \, \dd u \, \mu_\ell(\dd\bz) \, \dd w \\
  &\quad
    + \sum_{\ell=1}^d
       \int_0^t \sqrt{2 c_\ell \max \{0, X_{u,\ell}\}} \, \dd W_{u,\ell}
       \, \be_\ell^{(d)}
    + \int_0^t \int_{\cU_d}
       \br \bbone_{\cU_d\setminus A}(\br) \, M(\dd u, \dd\br) \\
  &\quad
    + \sum_{\ell=1}^d
       \int_0^t \int_{\cU_d} \int_{\cU_1}
        \bz \bbone_{\cU_d\setminus A}(\bz)
        \bbone_{\{w\leq X_{u-,\ell}\}} \, \tN_\ell(\dd u, \dd\bz, \dd w) \\
  &\quad
    + \sum_{\ell=1}^d
       \int_0^t \int_{\cU_d} \int_{\cU_1}
        \bz \bbone_A(\bz) \bbone_{\{w\leq X_{u-,\ell}\}} \, N_\ell(\dd u, \dd\bz, \dd w) \\
  &\quad
    + \int_0^t \int_{\cU_d} \br \bbone_A(\br) \, M(\dd u, \dd\br) , \qquad t\in\RR_+,
  \end{split}
 \end{align}
 since, by page 62 in Ikeda and Watanabe \cite{IkeWat}, part (v) in Definition \ref{Def_admissible} and \eqref{moment_condition_m_new},
 we have
 \begin{align*}
   \EE\left(\int_0^t \int_{\cU_d}
     \|\br\| \, M(\dd u, \dd\br)  \right)
  & = \int_0^t \int_{\cU_d}
       \|\br\|  \, \dd u\,\nu(\dd\br)
   = t \int_{\cU_d} \|\br\| \, \nu(\dd\br)
   < \infty,
 \end{align*}
 and, by page 62 in Ikeda and Watanabe \cite{IkeWat}, part (vi) in Definition \ref{Def_admissible} and \eqref{help_int_moment},
 for all \ $\ell\in\{1,\ldots,d\}$, \ we get
 \begin{align}\label{help6}
 \begin{split}
  &\EE\biggl(\int_0^t \int_{\cU_d} \int_{\cU_1}
             \|\bz\| \bbone_{A}(\bz)
             \bbone_{\{w\leq X_{u,\ell}\}}
             \, \dd u \, \mu_\ell(\dd\bz) \, \dd w\biggr)\\
  &= \int_0^t \int_{\cU_d}
         \|\bz\| \bbone_{A}(\bz) \EE(X_{u,\ell})
         \, \dd u \, \mu_\ell(\dd\bz) \\
  &\leq  \int_0^t \Vert \EE(\bX_u) \Vert \, \dd u
        \int_{\cU_d}  \|\bz\| \bbone_{A}(\bz)  \, \mu_\ell(\dd\bz)\\
  & = \int_0^t \Vert \EE( \bX_u ) \Vert \, \dd u
        \left(\int_{\cU_d}  \|\bz\| \bbone_{A}(\bz) \bbone_{\{ \Vert \bz\Vert \leq 1\} }   \, \mu_\ell(\dd\bz)
             + \int_{\cU_d}  \|\bz\| \bbone_{A}(\bz) \bbone_{\{ \Vert \bz\Vert > 1\}} \, \mu_\ell(\dd\bz)\right) \\
  &\leq \int_0^t \Vert \EE(\bX_u ) \Vert \, \dd u
        \left( \mu_\ell(A) + \int_{\cU_d}  \|\bz\| \bbone_{\{ \Vert \bz\Vert > 1\}} \, \mu_\ell(\dd\bz)  \right)
     < \infty.
 \end{split}
 \end{align}

{\sl Step 2.}
As a consequence of the SDE \eqref{help3} for $(\bX_t)_{t\in\RR_+}$,
 the jumps of \ $(\bX_t)_{t\in\RR_+}$ \ are related to the last four terms of the right hand side of \eqref{help3}.
We check that
 \begin{align}\label{help4}
  J_t(A)
  = \sum_{\ell=1}^d
      \int_0^t \int_{\cU_d} \int_{\cU_1} \bbone_A(\bz)
       \bbone_{\{w\leq X_{u-,\ell}\}}
       \, N_\ell(\dd u, \dd\bz, \dd w)
    + \int_0^t \int_{\cU_d} \bbone_A(\br) \, M(\dd u, \dd\br)
 \end{align}
 almost surely for all $t\in\RR_{++}$.
To prove \eqref{help4}, it is enough to verify that the last four terms on the right hand side of \eqref{help3}
 cannot have jumps simultaneously almost surely, and that the jumps with size vectors in $A$ correspond to the last two terms
 on the right hand side of \eqref{help3}.
These follow by the following steps (a)--(e):
\begin{itemize}
 \item[(a)] Note that the mutually independent Poisson point processes corresponding to \ $N_\ell$, \ $\ell \in\{1, \ldots, d\}$, \ and
      \ $M$, \ do not jump simultaneously almost surely, since the corresponding characteristic measures $\mu_\ell$, \ $\ell \in\{1, \ldots, d\}$,
      \ and \ $\nu$ \ are \ $\sigma$-finite (see the paragraph after Definition \ref{Def_admissible}),
       following, e.g., from the Disjointness Lemma and the paragraph after its proof in Section 2.2 in Kingman \cite{Kin}.
 \item[(b)] Observe that, under the moment condition \eqref{moment_condition_m_new}, by page 62 in Ikeda and Watanabe \cite{IkeWat},
       we have that for all \ $\omega\in\Omega$ \ and \ $t\in\RR_{++}$,
       \begin{align*}
         &\left(\int_0^t \int_{\cU_d} \br \bbone_{\cU_d\setminus A}(\br) \, M(\dd u, \dd\br) \right)(\omega)
             = \sum_{s\leq t, \; s\in D_{p_M(\omega)}} p_M(\omega,s)\bbone_{\cU_d\setminus A}(p_M(\omega,s)),\\
        &\left(\int_0^t \int_{\cU_d} \br \bbone_A(\br) \, M(\dd u, \dd\br)\right)(\omega)
             = \sum_{s\leq t, \; s\in D_{p_M(\omega)}} p_M(\omega,s)\bbone_A(p_M(\omega,s)),
      \end{align*}
      where \ $p_M: \Omega \to \Pi_{\cU_d}$ \ denotes a Poisson point process on \ $\cU_d$ \
      with characteristic measure \ $\nu$, \ $\Pi_{\cU_d}$ is the collection of point functions
      on $\cU_d$,  for all $\omega\in\Omega$, the countable set $D_{p_M(\omega)}\subseteq (0,\infty)$ is the domain
      of the point function $p_M(\omega)$, \ and the value of $p_M(\omega)$ at $s\in D_{p_M(\omega)}$ is simply denoted by
      $p_M(\omega, s)$ \ (see Ikeda and Watanabe \cite[Chapter 1, Section 9]{IkeWat}).
     Consequently, the set of jump times of $\big(\int_0^t \int_{\cU_d} \br \, M(\dd u, \dd\br)\big)_{t\in\RR_+}$
     coincides with $D_{p_M}$,
      \ $\big(\int_0^t \int_{\cU_d} \br \bbone_{\cU_d\setminus A}(\br) \, M(\dd u, \dd\br)\big)_{t\in\RR_+}$ \
     does not have any jump with size vectors in \ $A$, \ and all the size vectors of the jumps of \ $\big(\int_0^t \int_{\cU_d} \br \bbone_A(\br) \, M(\dd u, \dd\br)\big)_{t\in\RR_+}$ \
     are in \ $A$.
 \item[(c)] Similarly to step (b), using \eqref{help6}, by page 62 in Ikeda and Watanabe \cite{IkeWat}, for each \ $\ell\in\{1,\ldots,d\}$,
      we have that for all $\omega\in\Omega$ and $t\in\RR_{++}$,
      \begin{align*}
       \Bigg(\int_0^t \int_{\cU_d} \int_{\cU_1}
       &\bz \bbone_A(\bz) \bbone_{\{w\leq X_{u-,\ell}\}} \, N_\ell(\dd u, \dd\bz, \dd w)\Bigg)(\omega)\\
       &\qquad = \sum_{s\leq t, \; s\in D_{p_{N_\ell}(\omega)}} p_{N_\ell}^{(1)}(\omega,s)
                 \bbone_A\big(p_{N_\ell}^{(1)}(\omega,s)\big) \bbone_{\big\{ p_{N_\ell}^{(2)}(\omega,s) \leq X_{s-,\ell}\big\} },
       \end{align*}
       where \ $p_{N_\ell}=(p_{N_\ell}^{(1)}, p_{N_\ell}^{(2)}): \Omega \to \Pi_{\cU_d\times \cU_1}$ \
      denotes a Poisson point process on \ $\cU_d\times \cU_1$ \ with characteristic measure \ $\mu_\ell\otimes \dd w$,
      $\Pi_{\cU_d\times \cU_1}$ is the collection of point functions
      on $\cU_d\times \cU_1$, for all $\omega\in\Omega$, the countable set $D_{N_\ell(\omega)}\subseteq (0,\infty)$ is the domain
      of the point function $p_{N_\ell}(\omega)$, \ and the value of $p_{N_\ell}(\omega)$ at $s\in D_{p_{N_\ell}(\omega)}$ is simply denoted by
      $p_{N_\ell}(\omega, s)$ \ (see Ikeda and Watanabe \cite[Chapter 1, Section 9]{IkeWat}).
      Consequently, the jump times of $\big(\int_0^t \int_{\cU_d} \int_{\cU_1}
      \bz \bbone_A(\bz) \bbone_{\{w\leq X_{u-,\ell}\}} \, N_\ell(\dd u, \dd\bz, \dd w)\big)_{t\in\RR_+}$
      are in $D_{p_{N_\ell}}$, and every jump's size vector is in $A$.
 \item[(d)] It remains to check that, for each $\ell\in\{1,\ldots,d\}$,  the jump times of
       \[
         \left(\int_0^t \int_{\cU_d} \int_{\cU_1} \bz \bbone_{\cU_d\setminus A}(\bz)
       \bbone_{\{w\leq X_{u-,\ell}\}} \, \tN_\ell(\dd u, \dd\bz, \dd w)\right)_{t\in\RR_+}
       \]
       are in $D_{p_{N_\ell}}$, and every jump's size vector is in $\cU_d\setminus A$.
       Note that
      \begin{align}\label{help7}
       \begin{split}
        &\int_0^t \int_{\cU_d} \int_{\cU_1} \bz \bbone_{\cU_d\setminus A}(\bz)
         \bbone_{\{w\leq X_{u-,\ell}\}} \, \tN_\ell(\dd u, \dd\bz, \dd w)\\
        & = \int_0^t \int_{\cU_d} \int_{\cU_1} \bz \bbone_{\cU_d\setminus A}(\bz)
                \bbone_{\{w\leq X_{u-,\ell}\}}  \bbone_{\{\Vert \bz\Vert\leq 1\}} \, \tN_\ell(\dd u, \dd\bz, \dd w)\\
        &\quad + \int_0^t \int_{\cU_d} \int_{\cU_1} \bz \bbone_{\cU_d\setminus A}(\bz)
              \bbone_{\{w\leq X_{u-,\ell}\}}  \bbone_{\{\Vert \bz\Vert>1\}} \, \tN_\ell(\dd u, \dd\bz, \dd w),
              \qquad t\in\RR_+.
       \end{split}
      \end{align}
     By \eqref{help_int_moment} and part (vi) of Definition \ref{Def_admissible}, we have
     \begin{align*}
     &\EE\biggl(\int_0^t \int_{\cU_d} \int_{\cU_1}
             \|\bz\| \bbone_{\cU_d\setminus A}(\bz)
             \bbone_{\{w\leq X_{u,\ell}\}} \bbone_{\{\Vert \bz\Vert>1\}}
             \, \dd u \, \mu_\ell(\dd\bz) \, \dd w\biggr) \\
     &\leq \int_0^t \Vert \EE( \bX_u ) \Vert \, \dd u
            \int_{\cU_d}  \|\bz\| \bbone_{\{ \Vert \bz\Vert > 1\}} \, \mu_\ell(\dd\bz)
        < \infty,
     \end{align*}
     and hence the second term on right hand side of \eqref{help7} can be written as a sum of two terms:
     an integral with respect to \ $N_\ell$ \ and an integral with respect to \ $\dd u\,\mu_\ell(\dd \bz)\,\dd w$,
     \ and, similarly, as above, one can see that none of the stochastic processes corresponding to
     these two parts have jumps with size vectors in $A$.
     The first term on the right hand side of \eqref{help7} is handled in the next step (e).

\item[(e)] Finally, we check that all the jumps's size vectors of
      \begin{align}\label{help9}
        \left(\int_0^t \int_{\cU_d} \int_{\cU_1} \bz \bbone_{\cU_d\setminus A}(\bz)
                \bbone_{\{w\leq X_{u-,\ell}\}}  \bbone_{\{\Vert \bz\Vert\leq 1\}} \, \tN_\ell(\dd u, \dd\bz, \dd w)\right)_{t\in\RR_+}
      \end{align}
      are in \ $\cU_d\setminus A$.
      Similarly to the proof of part (iv) of Theorem 2.10 in Kyprianou \cite{Kyp},
      since the compensator measure \ $\dd u \, \mu_\ell(\dd\bz) \, \dd w$ \ has no atoms,
      every jump time of the process \eqref{help9} corresponds to a unique point in $D_{p_{N_\ell}}$.
      By \eqref{help_int_moment} and part (vi) of Definition \ref{Def_admissible}, we have
      \begin{align*}
     &\EE\left(\int_0^t \int_{\cU_d} \int_{\cU_1}
             \|\bz\|^2 \bbone_{\cU_d\setminus A}(\bz)
             \bbone_{\{w\leq X_{u,\ell}\}} \bbone_{\{\Vert \bz\Vert\leq 1\}}
             \, \dd u \, \mu_\ell(\dd\bz) \, \dd w \right) \\
     &\leq \int_0^t \Vert \EE( \bX_u) \Vert \, \dd u
            \int_{\cU_d}  \|\bz\|^2 \bbone_{\{ \Vert \bz\Vert \leq 1\}} \, \mu_\ell(\dd\bz)
        < \infty,
     \end{align*}
     and hence the process \eqref{help9} can be uniformly approximated on compact intervals almost surely by a sequence of appropriate square integrable martingales
     of which the jumps's size vectors are in \ $\cU_d\setminus A$, \ see, e.g., Applebaum \cite[part (2) of Theorem 4.3.4]{App}.
     Consequently, using also the form of the approximating martingales in question,
     we have that all the jumps of the process \eqref{help9} are in \ $\cU_d\setminus A$.
     Indeed, it cannot happen that the size vectors of jumps of the process \eqref{help9} are on the boundary of \ $\cU_d\setminus A$,
     since the approximating martingales in question are constructed as stochastic integrals with respect to the
     same compensated Poisson random measure \ $\tN_\ell$.
 \end{itemize}

{\sl Step 3.}
By taking the expectation of both sides of \eqref{help4} and using page 62 in Ikeda and Watanabe \cite{IkeWat} and \eqref{help_int_moment},  we have
 \begin{align*}
  \EE(J_t(A))
   = \sum_{\ell=1}^d
      \int_0^t \int_{\cU_d}
       \bbone_A(\bz) \EE(X_{u,\ell})
       \, \dd u \, \dd\mu_\ell(\bz)
     + \int_0^t \int_{\cU_d}
        \bbone_A(\br)
        \, \dd u \, \nu(\dd\br)
   < \infty ,\qquad t\in\RR_{++}.
 \end{align*}
Indeed, by \eqref{help_int_moment}, for each \ $\ell \in \{1, \ldots, d\}$ \ and \ $t\in\RR_+$, \ we have
 \begin{align*}
  \int_0^t \!\int_{\cU_d} \bbone_A(\bz) \EE(X_{u,\ell}) \, \dd u \, \dd\mu_\ell(\bz)
  & \leq \int_0^t \!\int_{\cU_d} \bbone_A(\bz) \Vert \EE(  \bX_u ) \Vert \, \dd u \, \dd\mu_\ell(\bz)\\
  & = \mu_\ell(A)  \int_0^t \Vert \EE(  \bX_u ) \Vert \, \dd u < \infty ,
 \end{align*}
 and
 \[
   \int_0^t \int_{\cU_d} \bbone_A(\bz) \, \dd u \, \nu(\dd\br)
   = t \nu(A) < \infty .
 \]
Consequently, for all \ $t \in \RR_{++}$, \ we have \ $\EE(J_t(A))<\infty$ \ and hence \ $\PP(J_t(A) < \infty) = 1$, \ i.e., we get part (i).
It also yields that \ $\PP(\Delta \bX_{\tau_A} \in A \mid \tau_A<\infty)=1$.

{\sl Step 4.}
Motivated by the SDE \eqref{help3}, we will consider another SDE (see \eqref{SDE_atirasa_dimd_R}) in a way that
 we remove all the parts from the SDE \eqref{help3} which correspond to the jumps with size vectors in the set \ $A$.
\ We will check that this new SDE \eqref{SDE_atirasa_dimd_R} admits a pathwise unique strong solution,
 which is a multi-type CBI process, and we calculate that its branching mechanism and immigration mechanism take the form \eqref{formula_psi_A} and \eqref{help54}, respectively.
Let \ $\tbD^{(A)} \in \RR_+^{d\times d}$ \ be the matrix with \ $\ell^{\mathrm th}$
 column given by \ $\int_{\cU_d} \bz \bbone_A(\bz) \, \mu_\ell(\dd\bz)$, \ $\ell \in \{1, \ldots, d\}$.
The entries of \ $\tbD^{(A)}$ \ are indeed in \ $\RR_+$, \ since $\int_{\cU_d} \Vert \bz\Vert \bbone_A(\bz) \, \mu_\ell(\dd\bz)<\infty$,
 $\ell\in\{1,\ldots,d\}$, following from \eqref{help6}.
First, note that
 \begin{align*}
  &\sum_{\ell=1}^d
    \int_0^t \int_{\cU_d} \int_{\cU_1}
     \bz \bbone_A(\bz) \bbone_{\{w\leq X_{u,\ell}\}}
     \, \dd u \, \mu_\ell(\dd\bz) \, \dd w
    = \sum_{\ell=1}^d
      \int_0^t \int_{\cU_d}
       \bz \bbone_A(\bz) X_{u,\ell} \, \dd u\,\mu_\ell(\dd\bz) \\
  &\qquad = \int_0^t \left(\sum_{\ell=1}^d  X_{u,\ell}
            \int_{\cU_d} \bz \bbone_A(\bz) \, \mu_\ell(\dd\bz)\right) \dd u
           = \int_0^t \tbD^{(A)} \bX_u \, \dd u, \qquad  t\in\RR_+.
 \end{align*}
Let us consider the SDE
 \begin{align}\label{SDE_atirasa_dimd_R}
  \begin{split}
  \bX_t^{(A)}
  &= \bX_0^{(A)} + \int_0^t \big(\Bbeta + (\tbB - \tbD^{(A)}) \bX_u^{(A)} \big) \, \dd u\\
  &\quad + \sum_{\ell=1}^d
        \int_0^t \sqrt{2 c_\ell \max \{0, X_{u,\ell}^{(A)}\}} \, \dd W_{u,\ell}
        \, \be_\ell^{(d)} \\
  &\quad
    + \sum_{\ell=1}^d
       \int_0^t \int_{\cU_d} \int_{\cU_1}
        \bz \bbone_{\cU_d\setminus A}(\bz)
        \bbone_{\{w\leq X_{u-,\ell}^{(A)}\}} \, \tN_\ell(\dd u, \dd\bz, \dd w) \\
  &\quad
    + \int_0^t \int_{\cU_d} \br \bbone_{\cU_d\setminus A}(\br) \, M(\dd u, \dd\br) , \qquad t \in \RR_+ ,
  \end{split}
 \end{align}
 with \ $\bX_0^{(A)} := \bX_0$.
\ For all \ $t \in \RR_+$ \ and \ $\ell\in\{1,\ldots,d\}$, \ by the definition of stochastic integrals with respect to
 (compensated) Poisson random measures, we have
 \begin{align*}
  &\int_0^t \int_{\cU_d} \int_{\cU_1}
        \bz \bbone_{\cU_d\setminus A}(\bz)
        \bbone_{\{w\leq X_{u-,\ell}^{(A)}\}} \, \tN_\ell(\dd u, \dd\bz, \dd w)\\
  &\qquad = \int_0^t \int_{\cU_d\setminus A} \int_{\cU_1}
       \bz \bbone_{\{w\leq X_{u-,\ell}^{(A)}\}}
       \, \tN_\ell^{(\cU_d\setminus A)}(\dd u, \dd\bz, \dd w)
 \end{align*}
 and
 \begin{align*}
  \int_0^t \int_{\cU_d} \br \bbone_{\cU_d\setminus A}(\br) \, M(\dd u, \dd\br)
    = \int_0^t \int_{\cU_d\setminus A} \br\, M^{(\cU_d\setminus A)}(\dd u, \dd\br) ,
 \end{align*}
 where \ $N_1^{(\cU_d\setminus A)}, \ldots, N_d^{(\cU_d\setminus A)}$ \ denote the restrictions of
 \ $N_1, \ldots, N_d$ \ to \ $\cU_1 \times (\cU_d\setminus A) \times \cU_1$ \ and
 \ $M^{(\cU_d\setminus A)}$ \ denotes the restriction of \ $M$ \ to \ $\cU_1 \times (\cU_d\setminus A)$,
 and $\tN_\ell^{(\cU_d\setminus A)}(\dd u, \dd\bz, \dd w)
    := N_\ell^{(\cU_d\setminus A)}(\dd u, \dd\bz, \dd w) - \dd u \, \bbone_{\cU_d\setminus A}(\bz)\mu_\ell(\dd\bz) \, \dd w$,
 \ $\ell\in\{1,\ldots,d\}$.
\ Here \ $N_1^{(\cU_d\setminus A)}, \ldots, N_d^{(\cU_d\setminus A)}$ \ and \ $M^{(\cU_d\setminus A)}$ \ are mutually
 independent Poisson random measures on \ $\cU_1 \times (\cU_d\setminus A) \times \cU_1$ \ and on \ $\cU_1 \times (\cU_d\setminus A)$ \
 with intensity measures \ $\dd u \, \bbone_{\cU_d\setminus A}(\bz)\mu_\ell(\dd\bz) \, \dd w$, \ $\ell \in \{1, \ldots, d\}$, \ and
 \ $\dd u \, \bbone_{\cU_d\setminus A}(\br)\nu(\dd\br)$, \ respectively (see, e.g., Kyprianou \cite[Corollary 2.5]{Kyp}).
Hence the SDE \eqref{SDE_atirasa_dimd_R} takes the form
 \begin{align}\label{SDE_atirasa2_dimd_R}
  \begin{split}
  \bX_t^{(A)}
  &= \bX_0^{(A)} + \int_0^t \big(\Bbeta + (\tbB - \tbD^{(A)}) \bX_u^{(A)} \big) \, \dd u\\
  &\quad  + \sum_{\ell=1}^d
             \int_0^t \sqrt{2 c_\ell \max \{0, X_{u,\ell}^{(A)}\}} \, \dd W_{u,\ell}
              \, \be_\ell^{(d)} \\
  &\quad
    + \sum_{\ell=1}^d \int_0^t \int_{\cU_d\setminus A} \int_{\cU_1}
       \bz \bbone_{\{w\leq X_{u-,\ell}^{(A)}\}}
       \, \tN_\ell^{(\cU_d\setminus A)}(\dd u, \dd\bz, \dd w)\\
   &\quad + \int_0^t \int_{\cU_d\setminus A} \br\, M^{(\cU_d\setminus A)}(\dd u, \dd\br)
  \end{split}
 \end{align}
 for \ $t \in \RR_+$.
Here \ $\tbB^{(A)} := \tbB - \tbD^{(A)} \in \RR_{(+)}^{d\times d}$,
 \ since, with the notation \ $\tbB^{(A)} = (\tb_{i,j}^{(A)})_{i,j=1}^d$,
 \ for each \ $i, j \in \{1, \ldots, d\}$ \ with \ $i \ne j$, \ we have
 \begin{align}\label{help21}
   \tb_{i,j}^{(A)} = \tb_{i,j} - (\tbD^{(A)})_{i,j}
   = b_{i,j} + \int_{\cU_d} z_i (1 - \bbone_{A}(\bz)) \, \mu_j(\dd\bz)
   \in\RR_+,
 \end{align}
 where we used \eqref{help5} and that \ $b_{i,j}\in\RR_+$ \ for each \ $i,j\in\{1,\ldots,d\}$ \ with \ $i\ne j$ \
 (due to \ $\bB\in \RR_{(+)}^{d\times d}$).
Note that the matrices \ $\bB^{(A)}$ \ and \ $\tbB^{(A)}$ \ satisfy the relations
 \begin{align}\label{help58}
 \tb_{i,j}^{(A)}
     = b_{i,j}^{(A)} + \int_{\cU_d} (z_i - \delta_{i,j})^+ \, \mu_j^{(\cU_d\setminus A)}(\dd\bz),
     \qquad i,j\in\{1,\ldots,d\},
 \end{align}
 which are required (see \eqref{help5} for the corresponding relations for \ $\bB$ \ and \ $\tbB$).
\ Indeed, using \eqref{help5}, we have
 \begin{align*}
   &\tb_{i,j}^{(A)} - \int_{\cU_d}(z_i - \delta_{i,j})^+ \mu_j^{(\cU_d\setminus A)}(\dd\bz)\\
    & = \tb_{i,j} - (\tbD^{(A)})_{i,j} - \int_{\cU_d} (z_i - \delta_{i,j})^+ \,\mu_j^{(\cU_d\setminus A)}(\dd \bz)\\
    & = b_{i,j} + \int_{\cU_d}(z_i - \delta_{i,j})^+ \mu_j(\dd\bz)
         - \int_A z_i \,\mu_j(\dd\bz)
         - \int_{\cU_d\setminus A} (z_i - \delta_{i,j})^+ \mu_j(\dd \bz) \\
    & = b_{i,j} + \int_A \big( (z_i - \delta_{i,j})^+ - z_i \big)\,\mu_j(\dd\bz)\\
    & = b_{i,j} + \delta_{i,j}\int_A \big( (z_i - 1)^+ - z_i \big)\,\mu_i(\dd\bz)
      = b_{i,j}^{(A)},\qquad i,j\in\{1,\ldots,d\},
 \end{align*}
 where \ $\bB^{(A)} = (b_{i,j}^{(A)})_{i,j=1}^d\in\RR_{(+)}^{d\times d}$ \ is introduced
  in \eqref{mod_parameters}.
Let us introduce the diagonal matrix \ $\bD^{(A)}=(d^{(A)}_{i,j})_{i,j=1}^d$ \ given by
 \[
  d^{(A)}_{i,j}:=\delta_{i,j}\int_A \big( (z_i - 1)^+ - z_i \big)\,\mu_i(\dd\bz), \qquad i,j\in\{1,\ldots,d\}.
 \]
Then, by \eqref{mod_parameters}, we have \ $b_{i,j}^{(A)}  = b_{i,j} + d^{(A)}_{i,j}$, $i,j\in\{1,\ldots,d\}$.

Recall that \ $(d, \bc, \Bbeta, \bB^{(A)}, \nu^{(\cU_d\setminus A)}, \bmu^{(\cU_d\setminus A)})$ \ is a set of admissible parameters,
 $\EE(\Vert \bX_0^{(A)}\Vert)<\infty$,
 and the measure \ $\nu^{(\cU_d\setminus A)}$ \ satisfies the moment condition \eqref{moment_condition_m_new} (see \eqref{moment_condition_m_new_A}).
Consequently, by Barczy et al.~\cite[Theorem 4.6 and Section 5]{BarLiPap2}, the SDE \eqref{SDE_atirasa2_dimd_R} admits
 a pathwise unique c\`{a}dl\`{a}g strong solution \ $(\bX^{(A)}_t)_{t\in\RR_+}$, \ which is a multi-type CBI process
 with parameters \ $(d, \bc, \Bbeta, \bB^{(A)}, \nu^{(\cU_d\setminus A)}, \bmu^{(\cU_d\setminus A)})$.
\ The branching mechanism \ $\RR_+^d \ni \blambda \mapsto (\varphi^{(A)}_1(\blambda), \ldots, \varphi^{(A)}_d(\blambda))^\top \in \RR^d$ \
 and the immigration mechanism  \ $\RR_+^d \ni \blambda \mapsto \psi^{(A)}(\blambda) \in \RR_+$ \
 of \ $(\bX^{(A)}_t)_{t\in\RR_+}$ \  takes the form
 \begin{align*}
   \varphi^{(A)}_\ell(\blambda)
   & = c_\ell \lambda_\ell^2 -  \langle \bB^{(A)} \be_\ell^{(d)}, \blambda \rangle
      + \int_{\cU_d}
         \bigl( \ee^{- \langle \blambda, \bz \rangle} - 1
                + \lambda_\ell (1 \land z_\ell) \bigr) \mu_\ell^{(\cU_d\setminus A)}(\dd \bz)\\
   & = c_\ell \lambda_\ell^2 -  \langle (\bB + \bD^{(A)}) \be_\ell^{(d)}, \blambda \rangle
      + \int_{\cU_d\setminus A}
         \bigl( \ee^{- \langle \blambda, \bz \rangle} - 1
                + \lambda_\ell (1 \land z_\ell) \bigr)
          \, \mu_\ell(\dd \bz)\\
   & = c_\ell \lambda_\ell^2 -  \langle \bB \be_\ell^{(d)}, \blambda \rangle
        - \lambda_\ell  \int_A ((z_\ell-1)^+ - z_\ell)\,\mu_\ell(\dd\bz) \\
   &\quad + \int_{\cU_d\setminus A}
         \bigl( \ee^{- \langle \blambda, \bz \rangle} - 1
                + \lambda_\ell (1 \land z_\ell) \bigr)
          \, \mu_\ell(\dd \bz)\\
  &= \varphi_\ell(\blambda) + \int_A  (1  - \ee^{-\langle \blambda, \bz \rangle}) \, \mu_\ell(\dd \bz)
 \end{align*}
 for \ $\blambda \in \RR_+^d$ \ and \ $\ell\in\{1,\ldots,d\}$ \ (where, at the last step, we used that \ $z_\ell\land 1 = -((z_\ell-1)^+ - z_\ell)$),
 and
 \begin{align*}
   \psi^{(A)}(\blambda)
    = \langle \bbeta, \blambda \rangle
       + \int_{\cU_d}
           \bigl( 1 - \ee^{- \langle\blambda, \br\rangle} \bigr)
              \, \nu^{(\cU_d\setminus A)}(\dd\br)
     = \langle \bbeta, \blambda \rangle
       + \int_{\cU_d\setminus A}
           \bigl( 1 - \ee^{- \langle\blambda, \br\rangle} \bigr)
              \, \nu(\dd\br)
 \end{align*}
 for \ $\blambda \in \RR_+^d$, \ respectively.
For completeness, we note that the integral $\int_A  (1  - \ee^{-\langle \blambda, \bz \rangle}) \, \mu_\ell(\dd \bz)$
 in the expression for $ \varphi^{(A)}_\ell(\blambda)$ is indeed finite, since,
 using $1-\ee^{-x}\leq x$, $x\in\RR$, and \eqref{help6}, we have that
 \begin{align*}
    \int_A  \vert 1  - \ee^{-\langle \blambda, \bz \rangle} \vert \, \mu_\ell(\dd \bz)
       \leq \int_A  \langle \blambda, \bz \rangle \, \mu_\ell(\dd \bz)
       \leq \int_A  \Vert \blambda \Vert \Vert \bz\Vert  \, \mu_\ell(\dd \bz)
        <\infty.
 \end{align*}
In particular, we get formulae \eqref{formula_psi_A} and \eqref{help54} in part (iii) (of the present Theorem \ref{Thm_jump_times_dCBI}).
Intuitively, the process \ $(\bX^{(A)}_t)_{t\in\RR_+}$ \ is obtained by removing from
 \ $(\bX_t)_{t\in\RR_+}$ \ all the masses produced by the jumps with size vectors in the set \ $A$.

{\sl Step 5.}
We check that, for all \ $t \in \RR_{++}$, \ we have
 \begin{align}\label{formula2}
 \begin{split}
  \{\tau_A > t\}
  = \biggl\{\sum_{\ell=1}^d
              \int_0^t \int_{A} \int_{\cU_1}
              \bbone_{\{w\leq X_{u-,\ell}^{(A)}\}}
              \, N_\ell^{(A)}(\dd u, \dd\bz, \dd w)
  + \int_0^t \int_{A} 1 \, M^{(A)}(\dd u, \dd\br) = 0\biggr\}
  \end{split}
 \end{align}
 up to a \ $\PP$-null set, where \ $N_1^{(A)}, \ldots, N_d^{(A)}$ \ are the restrictions of \ $N_1, \ldots, N_d$ \ to
 \ $\cU_1 \times A \times \cU_1$, \ and \ $M^{(A)}$ \ is the restriction of \ $M$ \ to \ $\cU_1 \times A$.
\ Here \ $N_1^{(A)}, \ldots, N_d^{(A)}$ \ and \ $M^{(A)}$ \ are mutually independent Poisson random measures on
 \ $\cU_1 \times A \times \cU_1$ \ and on \ $\cU_1 \times A$ \ with intensity measures
 \ $\dd u \, \bbone_A(\bz)\mu_\ell(\dd\bz) \, \dd w$, \ $\ell \in \{1, \ldots, d\}$, \ and
 \ $\dd u \, \bbone_A(\br)\nu(\dd\br)$, \ respectively (see, e.g., Kyprianou \cite[Corollary 2.5]{Kyp}).
Note that \ $\{\tau_A > 0\} \subseteq \{ J_s(A) = 0, \; s\in(0,\tau_A)\}$, \ and hence, using \eqref{help4},
 we get the following inclusion up to a $\PP$-null set:
 \begin{align*}
  \{\tau_A > 0 \}
     \subseteq & \biggl\{\sum_{\ell=1}^d
                \int_0^s \int_{\cU_d} \int_{\cU_1}
                \bbone_A(\bz)
                \bbone_{\{w\leq X_{u-,\ell}\}}
                \, N_\ell(\dd u, \dd\bz, \dd w) = 0, \; s\in[0,\tau_A) \biggr\}\\
               & \cap \biggl\{ \int_0^s \int_{\cU_d} \bbone_A(\br) \, M(\dd u, \dd\br) = 0,  \; s\in[0,\tau_A) \biggr\}.
 \end{align*}
In view of \eqref{help3} and \eqref{SDE_atirasa_dimd_R}, this implies that \ $\PP(\bX_s = \bX_s^{(A)}, s\in[0,\tau_A))=1$, \ which
 yields that \ $\PP(\bX_{s-} = \bX_{s-}^{(A)}, s\in(0,\tau_A])=1$.
\ Hence, using again \eqref{help4}, for all $t\in\RR_{++}$, we obtain the following inclusion up to a \ $\PP$-null set:
 \begin{align*}
  &\{\tau_A>t\}= \{\tau_A>t\}\cap \{J_t(A)=0\}\\
  & = \{\tau_A>t\}
     \cap \biggl\{\sum_{\ell=1}^d
              \int_0^t \int_{\cU_d} \int_{\cU_1}
              \bbone_A(\bz)
              \bbone_{\{w\leq X_{u-,\ell}\}}
              \, N_\ell(\dd u, \dd\bz, \dd w)
           + \int_0^t \int_{\cU_d} \bbone_A(\br) \, M(\dd u, \dd\br) = 0\biggr\} \\
   & = \{\tau_A>t\}
     \cap \biggl\{\sum_{\ell=1}^d
              \int_0^t \int_{A} \int_{\cU_1}
              \bbone_{\{w\leq X_{u-,\ell}\}}
              \, N_\ell^{(A)}(\dd u, \dd\bz, \dd w)
           + \int_0^t \int_{A} 1 \, M^{(A)}(\dd u, \dd\br) = 0\biggr\} \\
   &\subseteq
     \biggl\{\sum_{\ell=1}^d
              \int_0^t \int_{A} \int_{\cU_1}
              \bbone_{\{w\leq X^{(A)}_{u-,\ell}\}}
              \, N_\ell^{(A)}(\dd u, \dd\bz, \dd w)
           + \int_0^t \int_{A} 1 \, M^{(A)}(\dd u, \dd\br) = 0\biggr\}.
 \end{align*}

In what follows, we show that the reversed inclusion up to a $\PP$-null set holds as well.
Similarly as above, using part (i) (which was already proved, see Step 3),
 \eqref{help4} and \ $\PP(\bX_{s-} = \bX_{s-}^{(A)}, s\in(0,\tau_A])=1$, \
 for all $t\in\RR_{++}$, \ we have the following inclusions up to a \ $\PP$-null set:
 \begin{align*}
  \{\tau_A\leq t\} & \subseteq \{\tau_A\leq t\} \cap \{J_{\tau_A}(A)>0\} \\
  &= \{\tau_A\leq t\}
     \cap \biggl\{\sum_{\ell=1}^d
              \int_0^{\tau_A} \int_{\cU_d} \int_{\cU_1}
               \bbone_A(\bz)
              \bbone_{\{w\leq X_{u-,\ell}\}}
              \, N_\ell(\dd u, \dd\bz, \dd w)\\
  &\phantom{\subseteq \{\tau_A\leq t\} \cap \biggl\{}
           + \int_0^{\tau_A} \int_{\cU_d}  \bbone_A(\br) \, M(\dd u, \dd\br) > 0\biggr\} \\
   & = \{\tau_A\leq t\}
     \cap \biggl\{\sum_{\ell=1}^d
              \int_0^{\tau_A} \int_{A} \int_{\cU_1}
              \bbone_{\{w\leq X_{u-,\ell}\}}
              \, N^{(A)}_\ell(\dd u, \dd\bz, \dd w)\\
   &\phantom{= \{\tau_A\leq t\} \cap \biggl\{}
           + \int_0^{\tau_A} \int_{A} 1  \, M^{(A)}(\dd u, \dd\br) > 0\biggr\} =
 \end{align*}
 \begin{align*}
   &= \{\tau_A\leq t\}
     \cap \biggl\{\sum_{\ell=1}^d
              \int_0^{\tau_A}  \int_{A} \int_{\cU_1}
              \bbone_{\{w\leq X^{(A)}_{u-,\ell}\}}
              \, N^{(A)}_\ell(\dd u, \dd\bz, \dd w)\\
   &\phantom{=\{\tau_A\leq t\} \cap \biggl\{}
           + \int_0^{\tau_A}  \int_{A} 1 \, M^{(A)}(\dd u, \dd\br) > 0\biggr\}\\
   & \subseteq
     \biggl\{\sum_{\ell=1}^d
              \int_0^t \int_{A} \int_{\cU_1}
              \bbone_{\{w\leq X^{(A)}_{u-,\ell}\}}
              \, N^{(A)}_\ell(\dd u, \dd\bz, \dd w)
           + \int_0^t \int_{A} 1 \, M^{(A)}(\dd u, \dd\br) > 0\biggr\}.
 \end{align*}
Hence, by taking complement, we have
 \[
     \biggl\{\sum_{\ell=1}^d
              \int_0^t \int_{A} \int_{\cU_1}
              \bbone_{\{w\leq X^{(A)}_{u-,\ell}\}}
              \, N^{(A)}_\ell(\dd u, \dd\bz, \dd w)
           + \int_0^t \int_{A} 1 \, M^{(A)}(\dd u, \dd\br) = 0\biggr\}
      \subseteq \{\tau_A> t\}.
 \]
Consequently, we obtain \eqref{formula2}.

{\sl Step 6.}
Using \eqref{formula2} and \ $\bX_0=\bX^{(A)}_0$, \ we have for all $t\in\RR_{++}$,
 \begin{align*}
  &\PP(\tau_A > t \mid \bX_0 ) = \PP(\tau_A > t \mid \bX^{(A)}_0 )\\
  &= \PP\biggl(\sum_{\ell=1}^d
                \int_0^t \int_{A} \int_{\cU_1}
                 \bbone_{\{w\leq X_{u-,\ell}^{(A)}\}}
                 \, N_\ell^{(A)}(\dd u, \dd\bz, \dd w)
               + \int_0^t \int_{A}
                  1 \, M^{(A)}(\dd u, \dd\br) = 0 \,\Big\vert\, \bX^{(A)}_0 \biggr).
 \end{align*}
Since the process \ $(\bX_t^{(A)})_{t\in\RR_+}$ \ is a pathwise unique strong
 solution of the SDE \eqref{SDE_atirasa2_dimd_R} (see Step 4),
 it is progressively measurable with respect to the filtration generated by
 \[
    \Big\{ \bX_0^{(A)}, (W_{t,1})_{t\in\RR_+},\ldots,(W_{t,d})_{t\in\RR_+},
          N_1^{(\cU_d\setminus A)},\ldots,N_d^{(\cU_d\setminus A)}, M^{(\cU_d\setminus A)} \Big\},
 \]
 and all the $\sigma$-algebras belonging to this filtration are independent of \ $N_1^{(A)},\ldots,N_d^{(A)}$ \ and \ $M^{(A)}$.
This yields that the process \ $(\bX_t^{(A)})_{t\in\RR_+}$ \ is independent of \ $N_1^{(A)},\ldots, N_d^{(A)}$ \ and \ $M^{(A)}$.
\ Consequently, using also the independence of \ $N^{(A)}_1,\ldots, N^{(A)}_d$ \ and \ $M^{(A)}$,
 and that $\bX_0 = \bX_0^{(A)}$, \ by the tower rule for conditional expectation, for all \ $t \in \RR_{++}$, \ we obtain
 \begin{align*}
  &\PP(\tau_A > t \mid \bX_0)\\
  &= \EE\biggl( \PP\biggl(\biggl\{\int_0^t \int_{A}
                 1 \, M^{(A)}(\dd u, \dd\br) = 0 \biggr\} \\
  &\phantom{=\EE\biggl(\PP\biggl(\;}
                  \cap \bigcap_{\ell=1}^d \biggl\{ \int_0^t \int_{A} \int_{\cU_1}
                  \bbone_{\{w\leq X_{u-,\ell}^{(A)}\}}
                  \, N_\ell^{(A)}(\dd u, \dd\bz, \dd w) = 0 \biggr\}
                  \;\bigg \vert\;
                  (\bX_u^{(A)})_{u\in[0,t]}
                 \biggr) \,\bigg\vert\,  \bX_0^{(A)}\biggr)\\
  &= \EE\biggl( \PP\biggl( \int_0^t \int_{A} 1 \, M^{(A)}(\dd u, \dd\br) = 0 \biggr)\\
  &\phantom{=\EE\biggl[\,}\times \prod_{\ell=1}^d
             \PP\biggl( \int_0^t \int_{A} \int_{\cU_1}
                  \bbone_{\{w\leq X_{u-,\ell}^{(A)}\}}
                  \, N_\ell^{(A)}(\dd u, \dd\bz, \dd w) = 0
                  \;\bigg \vert\;
                  (\bX_u^{(A)})_{u\in[0,t]}
                  \biggr) \,\bigg\vert\, \bX_0^{(A)} \biggr) =
 \end{align*}
 \begin{align*}
  &= \PP\biggl( \int_0^t \int_{A} 1 \, M^{(A)}(\dd u, \dd\br) = 0 \biggr) \\
  &\phantom{=\,}
     \times\EE\biggl( \prod_{\ell=1}^d
             \PP\biggl( \int_0^t \int_{A} \int_{\cU_1}
                  \bbone_{\{w\leq X_{u-,\ell}^{(A)}\}}
                  \, N_\ell^{(A)}(\dd u, \dd\bz, \dd w) = 0
                  \;\bigg \vert\;
                  (\bX_u^{(A)})_{u\in[0,t]}
                  \biggr) \,\bigg\vert\, \bX_0^{(A)} \biggr) \\
  &= \ee^{-\nu(A)t} \EE\left( \exp\left\{ -\sum_{\ell=1}^d \mu_\ell(A) \int_0^t X_{u,\ell}^{(A)}\,\dd u \right\}
                 \,\bigg\vert\, \bX_0^{(A)} \right)\\
  &= \ee^{-\nu(A)t} \EE\left( \exp\left\{ -\sum_{\ell=1}^d \mu_\ell(A) \int_0^t X_{u,\ell}^{(A)}\,\dd u \right\}
                 \,\bigg\vert\, \bX_0 \right),
 \end{align*}
 yielding \eqref{formula1} in case of $t \in \RR_{++}$,  i.e., part (ii) in case of $t \in \RR_{++}$.
\ Indeed, the last but one equality can be checked as follows.
For all \ $t\in\RR_+$, \  we have \ $\int_0^t \int_{A} 1 \, M^{(A)}(\dd u, \dd\br)$ \ has a Poisson distribution
 with parameter \ $\int_0^t \int_{A} 1 \,\dd u\,\nu(\dd\br)=t\nu(A)$, \
 and for all \ $t\in\RR_+$ \ and \ $\ell\in\{1,\ldots,d\}$, \ the conditional distribution of
 \[
   \int_0^t \int_{A} \int_{\cU_1}
                  \bbone_{\{w\leq X_{u-,\ell}^{(A)}\}}
                  \, N_\ell^{(A)}(\dd u, \dd\bz, \dd w)
 \]
 given \ $(\bX_u^{(A)})_{u\in[0,t]}$ \ is a Poisson distribution with parameter
 \[
   \int_0^t \int_{A} \int_{\cU_1}
                  \bbone_{\{w\leq X_{u-,\ell}^{(A)}\}}
                  \,\dd u\,\mu_\ell(\dd\bz)\, \dd w
   = \mu_\ell(A) \int_0^t X^{(A)}_{u,\ell}\,\dd u.
 \]
Consequently, we have
 \begin{align*}
   \PP\left(  \int_0^t \int_{A}
                1 \, M^{(A)}(\dd u, \dd\br) =0 \right)
      = \ee^{-\nu(A)t},\qquad t\in\RR_+,
 \end{align*}
 and for each \ $\ell\in\{1,\ldots,d\}$, \
 \begin{align*}
   &\PP\left( \int_0^t \int_{A} \int_{\cU_1}
                  \bbone_{\{w\leq X_{u-,\ell}^{(A)}\}}
                  \, N_\ell^{(A)}(\dd u, \dd\bz, \dd w) =0
            \;\bigg\vert\; (\bX_u^{(A)})_{u\in[0,t]} \right)\\
   &\qquad =  \exp\left\{- \mu_\ell(A) \int_0^t X_{u,\ell}^{(A)}\,\dd u \right\},
       \qquad  t\in\RR_+.
 \end{align*}
In case of $t=0$, \eqref{formula1} follows from the facts that \eqref{formula1}
 holds for all $t\in\RR_{++}$ and both sides of \eqref{formula1} as functions of \ $t$
 \ are right continuous at \ $0$. \
Indeed, the latter statement in case of the right hand side of \eqref{formula1} follows from the dominated convergence theorem
 for conditional expectations, and, in case of the left hand side of \eqref{formula1},
 from
 \begin{align*}
  \lim_{t\downarrow 0} \PP(\tau_A > t \mid \bX_0)
    & = \lim_{t\downarrow 0} \EE(\bone_{\{\tau_A > t \}} \mid \bX_0)
      = \EE\big(\lim_{t\downarrow 0} \bone_{\{\tau_A > t \}} \mid \bX_0\big) \\
    & = \EE( \bone_{\{\tau_A > 0 \}} \mid \bX_0)
      = \PP( \tau_A > 0  \mid \bX_0).
 \end{align*}

{\sl Step 7.}
Formula \eqref{formula3} directly follows by \eqref{formula1} and Proposition \ref{Pro_int_CBI_Laplace}
 taking into account that the branching and immigration mechanisms of $(\bX^{(A)}_t)_{t\in\RR_+}$ \ have been calculated in Step 4,
  i.e., we get part (iii) as well.
\proofend

Next, we formulate a corollary of Theorem \ref{Thm_jump_times_dCBI} for multi-type CB processes.

\begin{Cor}\label{Cor_CB_tauA}
Let \ $(\bZ_t)_{t\in\RR_+}$ \ be a multi-type CB process with parameters \ $(d, \bc, \bzero, \bB, 0, \bmu)$.
Then for all \ $t \in \RR_+$ \ and \ $A\in\cB(\cU_d)$ \ with \ $\sum_{\ell=1}^d \mu_\ell(A)<\infty$, \ we have
  \begin{align}\label{formula3_CB}
     \PP_\bz(\tau_A>t)
        = \exp\left\{ - \sum_{\ell=1}^d z_\ell\,\tv_\ell^{(A)}(t,\bmu(A)) \right\},
     \qquad  \bz=(z_1,\ldots,z_d)^\top\in\RR_+^d,
  \end{align}
 where $\bmu(A)=(\mu_1(A),\ldots,\mu_d(A))$ and the function $\RR_+\ni t \mapsto \tbv^{(A)}(t,\bmu(A))\in\RR_+^d$
 is given in \eqref{tvA}.
\end{Cor}

\noindent{\bf Proof.}
This readily follows from part (iii) of Theorem \ref{Thm_jump_times_dCBI}, since, by the assumption, \ $\psi\equiv 0$, \
 which yields that $\nu(A)=0$ and $\psi^{(A)}\equiv 0$.
\proofend

The following corollary can be considered as a multi-type counterpart of Corollaries 3.1 and 3.2 in He and Li \cite{HeLi}.

\begin{Cor}\label{Cor_tau}
Let \ $(\bX_t)_{t\in\RR_+}$ \ be a multi-type CBI process with parameters \ $(d, \bc, \Bbeta, \bB, \nu, \bmu)$ \
  such that the moment condition \eqref{moment_condition_m_new} holds.
\begin{enumerate}
 \item[(i)] If \ $A,B\in\cB(\cU_d)$ \ are such that \ $A\subseteq B$ \ and \ $\sum_{i=1}^d \mu_i(B)<\infty$, \ then
            \[
            \tbv^{(A)}(t,\bmu(A)) \leq \tbv^{(B)}(t,\bmu(B)),\qquad t\in\RR_+,
            \]
            where the functions \ $\RR_+\ni t \mapsto \tbv^{(A)}(t,\bmu(A))$ \ and \ $\RR_+\ni t \mapsto \tbv^{(B)}(t,\bmu(B))$ \
            are the unique locally bounded solutions of the system of differential equations \eqref{help_DE_tv} with the
            choices $A$ and $B$, respectively.
 \item[(ii)] If \ $A\in\cB(\cU_d)$ \ is such that \ $\nu(A)+\sum_{i=1}^d \mu_i(A)=0$, \ then
             \[
               \PP_\bx(\tau_A=\infty) = 1 \qquad  \mbox{ and } \qquad \PP_\bx(J_t(A)=0)=1, \qquad \bx\in\RR_+^d, \; t\in\RR_{++}.
             \]
 \item[(iii)] If one of the following two conditions hold:
               \begin{itemize}
                  \item[(a)] $A\in\cB(\cU_d)$ \ is such that \ $\nu(A)\in(0,\infty]$,
                  \item[(b)] $A\in\cB(\cU_d)$ \ is such that \ $\nu(A)=0$ \ and \ $\sum_{i=1}^d \mu_i(A)\in(0,\infty)$,
                             \ $(\bX^{(A)}_t)_{t\in\RR_+}$ \ is irreducible and \ $\psi$ \ is not identically zero,
               \end{itemize}
              then
              \[
                \PP_\bx(\tau_A<\infty) = 1, \qquad \bx\in\RR_+^d.
              \]
 \item[(iv)] If \ $D_{\bvarphi^{(A)}}\ne \emptyset$, \ $\psi\equiv 0$ \ and \ $A\in\cB(\cU_d)$ \ is such that \ $0<\mu_i(A)<\infty$, \ $i\in\{1,\ldots,d\}$, \ then
              \[
                 \PP_\bx(\tau_A=\infty)
                      = \ee^{-\big\langle \bx,\bphi^{(A)}(\bmu(A)) \big\rangle},\qquad \bx\in\RR_+^d,
              \]
              provided that \ $\tbv^{(A)}(\infty,\bmu(A))\in(0,\infty)^d$, \ where
              \ $\bphi^{(A)}$ \ is the inverse of the branching mechanism
              \ $\bvarphi^{(A)}=(\varphi^{(A)}_1,\ldots,\varphi^{(A)}_d)$
              \ (for more details on the existence of \ $\bphi^{(A)}$, \ see the discussion
              after Corollary \ref{Cor_int_CBI_Laplace}).
\end{enumerate}
\end{Cor}

\noindent{\bf Proof.}
(i):
First, note that the function \ $\tbv^{(A)}$ \ depends on the parameters \ $\bc$, \ $\bB$ \ and \ $\bmu$,
 \ but not on \ $\Bbeta$ \ and \ $\nu$, \ since the functions \ $\varphi_1^{(A)},\ldots,\varphi_d^{(A)}$ \ appearing
 in the system of differential equations \eqref{help_DE_tv} depend on \ $\bc$, \ $\bB$ \ and \ $\bmu$,
 but not on \ $\Bbeta$ \ and \ $\nu$.
\ Hence, in order to prove part (i), we may and do assume that  \ $\Bbeta=\bzero$ \ and \ $\nu=0$.
\ Let \ $(\bZ_t)_{t\in\RR_+}$ \ be a multi-type CB process with parameters \ $(d, \bc, \bzero, \bB, 0, \bmu)$.
Since \ $A\subseteq B$, \ we have \ $\sum_{i=1}^d \mu_i(A)<\infty$ \ and \ $\tau_A\geq \tau_B$, \ yielding that
 \[
    \PP_\bz(\tau_B>t)  \leq \PP_\bz(\tau_A>t), \qquad t\in\RR_+,\;\; \bz\in\RR_+^d.
 \]
Using \eqref{formula3_CB}, it implies that
 \[
   \exp\left\{ - \sum_{\ell=1}^d z_\ell\,\tv_\ell^{(B)}(t,\bmu(B)) \right\}
      \leq \exp\left\{ - \sum_{\ell=1}^d z_\ell\,\tv_\ell^{(A)}(t,\bmu(A)) \right\}
 \]
 for $t\in\RR_+$ and $\bz=(z_1,\ldots,z_d)^\top\in\RR_+^d$.
By choosing \ $\bz:=\be_i^{(d)}$, $i=1,\ldots,d$, \ we have \ $\tv_i^{(B)}(t,\bmu(B))\geq \tv_i^{(A)}(t,\bmu(A))$, \ $t\in\RR_+$, $i=1,\ldots,d$,
 \ implying that \ $\tbv^{(B)}(t,\bmu(B))\geq \tbv^{(A)}(t,\bmu(A))$, \ $t\in\RR_+$, \ as desired.

(ii): Let \ $A\in\cB(\cU_d)$ \ be such that \ $\nu(A)+\sum_{i=1}^d \mu_i(A)=0$.
\ Since \ $\bmu(A)=\bzero$ \ and \ $\varphi^{(A)}_\ell(\bzero) = 0$, $\ell\in\{1,\ldots,d\}$,
 by the uniqueness of a locally bounded solution of the system of differential equations \eqref{help_DE_tv},
 we get \ $\tbv^{(A)}(t,\bmu(A))=\tbv^{(A)}(t,\bzero)=\bzero$, $t\in\RR_+$.
\ Consequently, using again \ $\nu(A)=0$, \ $\bmu(A)=\bzero$ \ and \ $\psi^{(A)}(\bzero)=0$, \ \eqref{formula3} yields that
 \[
   \PP_\bx(\tau_A>t) = 1, \qquad \bx\in\RR_+^d, \;\; t\in\RR_+.
 \]
Therefore, using that $\{\tau_A>t\}\subseteq \{J_t(A)=0\}$, $t\in\RR_{++}$, we have $\PP_\bx(J_t(A)=0)=1$, $\bx\in\RR_+^d$, $t\in\RR_{++}$.
Further, we have
 \[
    \PP_\bx(\tau_A = \infty)
        = \lim_{t\to\infty} \PP_\bx(\tau_A>t) = 1,\qquad \bx\in\RR_+^d.
 \]

(iii): First, we check that we may and do assume that \ $0<\nu(A)+\sum_{i=1}^d \mu_i(A)<\infty$.
\ In case of (b), it holds trivially by our assumption.
In case of (a), if \ $\nu(A)\in(0,\infty]$ \ and \ $\nu(A)+\sum_{i=1}^d \mu_i(A)=\infty$,
 then, by Remark \ref{Rem_A_condition}, we have $\nu( A\setminus K_\epsilon)+\sum_{i=1}^d \mu_i(A\setminus K_\epsilon)<\infty$
 for all $\vare\in(0,1)$.
Using that $A\setminus K_\vare \uparrow A$ as $\vare\downarrow 0$, the continuity from below of the measure $\nu$ implies that
 $\nu(A\setminus K_\vare)\uparrow \nu(A)$ as $\vare\downarrow 0$.
Consequently, since $\nu(A)>0$, we get the existence of an $\epsilon\in(0,1)$ such that
  $0<\nu( A\setminus K_\epsilon)+\sum_{i=1}^d \mu_i(A\setminus K_\epsilon)<\infty$.
\ We also have \ $\tau_{A\setminus K_\epsilon}\geq \tau_A$, \ so if \ $\PP_\bx(\tau_{A\setminus K_\epsilon}<\infty) = 1$, \ $\bx\in\RR_+^d$,
 \ then \ $\PP_\bx(\tau_A<\infty) = 1$, \ $\bx\in\RR_+^d$, \ holds as well.
Therefore, it is enough to consider $A\in\cB(\cU_d)$ \ such that \ $0<\nu(A)+\sum_{i=1}^d \mu_i(A)<\infty$.

In what follows, let \ $A\in\cB(\cU_d)$ \ be such that \ $0<\nu(A)+\sum_{i=1}^d \mu_i(A)<\infty$.

In case of (a), we have \ $t\nu(A)\to\infty$ \ as \ $t\to\infty$, and, since \ $\sum_{i=1}^d \mu_i(A)<\infty$, \
 we also get \ $\bmu(A)\in\RR_+^d$.
\ Then using that \ $\tbv^{(A)}(t,\bmu(A))\in\RR_+^d$ \ and \ $\psi^{(A)}(\blambda)\in\RR_+$, \ $\blambda\in\RR_+^d$, \
 as a consequence of the continuity of probability and \eqref{formula3}, we can conclude that
 \[
   \PP_\bx(\tau_A=\infty) = \lim_{t\to\infty} \PP_\bx(\tau_A>t) \leq \lim_{t\to\infty} \ee^{-\nu(A)t}=0,
    \qquad \bx\in\RR_+^d,
 \]
 yielding that \ $\PP_\bx(\tau_A=\infty) =0$, \ $\bx\in\RR_+^d$, \ as desired.

In case of (b), we have \ $\nu(A)=0$, \ yielding that \ $\psi^{(A)}=\psi$, \ and, by our assumption, we have
  that \ $\sum_{i=1}^d \mu_i(A) \in (0,\infty)$.
\ It implies that \ $\bmu(A)\in\RR_+^d$ \ with \ $\bmu(A)\ne \bzero$.
\ By part (i) of Proposition \ref{Pro_tv}, we have that the function
 \ $\RR_+\ni t\mapsto \tbv^{(A)}(t,\bmu(A))\in\RR_+^d$ \ is increasing, and hence
 \[
   \tbv^{(A)}(\infty,\bmu(A)):=\lim_{t\to\infty} \tbv^{(A)}(t,\bmu(A)) \in[0,\infty]^d \qquad \text{exists}.
 \]
By the assumption, $(\bX^{(A)}_t)_{t\in\RR_+}$ is irreducible (i.e., \ $\tbB^{(A)}$, \ given in \eqref{help58}, is irreducible),
 and hence, by part (iii) of Proposition \ref{Pro_tv}, we have that \ $\tbv^{(A)}(\infty,\bmu(A))\in (0,\infty]^d$.
\ Since \ $\psi^{(A)} = \psi$ \ and \ $\psi$ \ is not identically zero, by Lemma \ref{Lem_psi_positive}, we get that
 \begin{align}\label{help13}
   \lim_{t\to\infty} \int_0^t \psi^{(A)}\big(\tbv^{(A)}(s,\bmu(A))\big)\,\dd s = \infty.
 \end{align}
Indeed, using that the limit of the increasing function \ $\RR_+\ni t\mapsto \tbv^{(A)}(t,\bmu(A))\in\RR_+^d$ \
 as \ $t\to\infty$ \ is \ $\tbv^{(A)}(\infty,\bmu(A))\in (0,\infty]^d$, \ we get that
 there exists \ $t_0\in\RR_+$ \ such that
 \[
   \tbv^{(A)}(t,\bmu(A)) \geq \tbv^{(A)}(t_0,\bmu(A)) > \bzero \qquad \text{for all \ $t\geq t_0$.}
 \]
Using that \ $\psi^{(A)}=\psi$ \ is monotone increasing, it implies that
 \[
   \psi^{(A)}(\tbv^{(A)}(t,\bmu(A)))\geq \psi^{(A)}\big(\tbv^{(A)}(t_0,\bmu(A))\big), \qquad t\geq t_0.
 \]
Further, since \ $\tbv^{(A)}(t_0,\bmu(A))\in(0,\infty)^d$, \ and \ $\psi^{(A)}=\psi$ \ is not identically zero,
 by Lemma \ref{Lem_psi_positive}, we have that $\psi^{(A)}\big(\tbv^{(A)}(t_0,\bmu(A))\big) = \psi\big(\tbv^{(A)}(t_0,\bmu(A))\big)>0$,
 which yields \eqref{help13}.

In view of \eqref{formula3} and \eqref{help13}, using the continuity of probability and that \ $\tbv^{(A)}(t,\bmu(A))\in\RR_+^d$, \ $t\in\RR_+$, \
 we have that
 \begin{align*}
  \PP_\bx(\tau_A=\infty)
           = \lim_{t\to\infty} \PP_\bx(\tau_A>t)
           \leq \limsup_{t\to\infty}  \exp\left\{ - \int_0^t \psi^{(A)}\big(\tbv^{(A)}(s,\bmu(A))\big)\,\dd s \right\}
           =0, \qquad \bx\in\RR_+^d,
 \end{align*}
 yielding that \ $\PP_\bx(\tau_A=\infty) =0$, \ $\bx\in\RR_+^d$, \ as desired.

(iv): First, note that $\bmu(A)\in(0,\infty)^d$.
\ By the continuity of probability, \eqref{formula3_CB} and parts (i) and (iv) of Proposition \ref{Pro_tv}, we have that
 \begin{align*}
   \PP_\bx(\tau_A=\infty)
    & =  \lim_{t\to\infty} \PP_\bx(\tau_A>t)
      = \exp\left\{ - \sum_{\ell=1}^d x_\ell \, \tv^{(A)}_\ell(\infty,\bmu(A)) \right\}\\
    & = \ee^{- \big\langle \bx, \tbv^{(A)}(\infty,\bmu(A)) \big\rangle}
      = \ee^{-\big\langle \bx,\bphi^{(A)}(\bmu(A)) \big\rangle},\qquad \bx = (x_1,\ldots,x_\ell)\in\RR_+^d,
 \end{align*}
 provided that \ $\tbv^{(A)}(\infty,\bmu(A))\in(0,\infty)^d$, \ as desired.
\proofend

\section{Local and global supremum for jumps}
\label{Section_max_jumps}

Recall that for all \ $r>0$, \ $K_r=\{\by\in\RR_+^d : \Vert \by\Vert<r\}$ denotes the open ball
 in $\RR_+^d$ around $\bzero$ with radius \ $r$.
For all \ $r>0$, \ let \ $\overline{K}_r$ \ be the closure of $K_r$, and let \ $K_r^c:=\cU_d\setminus K_r$ \ and
 \ $\overline{K}_r^c :=\cU_d\setminus \overline{K}_r$.
\ The next result is a generalization of Theorem 4.1 in He and Li \cite{HeLi} to multi-type CBI processes.

\begin{Pro}\label{Pro_max_jump}
Let \ $(\bX_t)_{t\in\RR_+}$ \ be a multi-type CBI process with parameters \ $(d, \bc, \Bbeta, \bB, \nu, \bmu)$ \
 such that the moment condition \eqref{moment_condition_m_new} holds.
 \begin{itemize}
 \item[(i)] Then for all \ $r>0$, \ $t>0$ \ and \ $\bx=(x_1,\ldots,x_d)^\top\in\RR_+^d$, \ we have \ $\PP(J_t(\overline{K}_r^c) < \infty)=1$, \ and
 \begin{align*}
    &\PP_\bx\Big(\sup_{s\in(0,t]} \Vert \Delta \bX_s\Vert \leq r\Big)\\
    &\qquad  = \exp\left\{ -\nu(\overline{K}_r^c)t - \sum_{\ell=1}^d x_\ell\,\tv_\ell^{(\overline{K}_r^c)}(t,\bmu(\overline{K}_r^c ))
                           - \int_0^t \psi^{(\overline{K}_r^c)}\big(\tbv^{(\overline{K}_r^c )}(s,\bmu(\overline{K}_r^c ))\big)\,\dd s \right\},
 \end{align*}
 where the continuously differentiable function
 \begin{align*}
   \RR_+\ni t \mapsto \tbv^{(\overline{K}_r^c)}(t,\bmu(\overline{K}_r^c ))
        :=(\tv^{(\overline{K}_r^c )}_1(t,\bmu(\overline{K}_r^c )),\ldots,\tv^{(\overline{K}_r^c )}_d(t,\bmu(\overline{K}_r^c)))^\top\in\RR_+^d
 \end{align*}
 is the unique locally bounded solution to the system of differential equations \eqref{help_DE_tv} with the choice $A:=\overline{K}_r^c$.
 \item[(ii)] If, in addition, \ $\nu(\cU_d) + \sum_{i=1}^d \mu_i(\cU_d)<\infty$, then for all $t>0$ and
             $\bx=(x_1,\ldots,x_d)^\top\in\RR_+^d$, we have \ $\PP(J_t(\cU_d) < \infty)=1$ \ and
             \begin{align*}
                &\PP_\bx\Big(\sup_{s\in(0,t]} \Vert \Delta \bX_s\Vert = 0\Big)\\
                &\qquad  = \exp\left\{ -\nu(\cU_d)t - \sum_{\ell=1}^d x_\ell\,\tv_\ell^{(\cU_d)}(t,\bmu(\cU_d))
                            - \int_0^t \Big\langle \Bbeta, \tbv^{(\cU_d)}(s,\bmu(\cU_d)) \Big\rangle \,\dd s \right\}.
            \end{align*}
 \end{itemize}
\end{Pro}

\noindent{\bf Proof.}
(i). Let \ $r>0$, \ $t>0$ \ and \ $\bx\in\RR_+^d$ \ be fixed arbitrarily.
Note that the closure of $\overline{K}_r^c$ coincides with $K_r^c = \cU_d\setminus K_r$,
 and since $\bzero$ is not contained in $K_r^c$, by Remark \ref{Rem_A_condition},
 we have that \ $\nu(\overline{K}_r^c) + \sum_{i=1}^d \mu_i(\overline{K}_r^c)<\infty$.
Then, by part (i) of Theorem \ref{Thm_jump_times_dCBI}, we have \ $\PP(J_t(\overline{K}_r^c) < \infty)=1$.
Next, we check that
 \begin{align}\label{help22}
   \PP_\bx\Big(\sup_{s\in(0,t]} \Vert \Delta \bX_s\Vert \leq r \Big)
     = \PP_\bx\big(\tau_{\overline{K}_r^c} > t \big).
 \end{align}
First, note that for all $\bz\in\RR^d$, \ the inclusion \ $\bz\in \overline{K}_r^c$ \ holds if and only if $\Vert \bz\Vert>r$.
If $\sup_{s\in(0,t]} \Vert \Delta \bX_s\Vert \leq r$ holds, then \ $\Vert \Delta \bX_s\Vert \leq r$, \ $s\in(0,t]$, \
 yielding that $\Delta \bX_s\notin \overline{K}_r^c$, $s\in(0,t]$,  and hence, due to the definition of $\tau_{\overline{K}_r^c}$,
 we have that $\tau_{\overline{K}_r^c} \geq t$.
Here $\tau_{\overline{K}_r^c} = t$ can hold only with probability zero,
 since, if $\tau_{\overline{K}_r^c} = t$, then, using that \ $\Vert \Delta \bX_t\Vert\leq r$, \
  we have that there exists a strictly decreasing sequence \ $(t_n)_{n\in\NN}$ \ such that $t_n\downarrow t$ as $n\to\infty$ and
 $\Vert \Delta \bX_{t_n}\Vert > r$, $n\in\NN$, and, in particular, we have $J_{t+1}(\overline{K}_r^c) = \infty$,
 which has probability zero.
This implies that the left hand side of \eqref{help22} is less than or equal to its right hand side.
Conversely, if $\tau_{\overline{K}_r^c} > t$, then $\Delta \bX_s \notin \overline{K}_r^c$, $s\in(0,t]$, and hence
 $\Vert \Delta \bX_s\Vert\leq r$, $s\in(0,t]$, \ yielding that
 \ $\sup_{s\in(0,t]} \Vert \Delta \bX_s\Vert \leq r$.
This implies that the right hand side of \eqref{help22} is less than or equal to its left hand side.
Thus we get \eqref{help22}.
Consequently, \eqref{help22} and part (iii) of Theorem \ref{Thm_jump_times_dCBI} yield part (i).

(ii).  The same proof works as in case of (i).
Namely, similarly as \eqref{help22} with the convention \ $\overline{K}_0^c:= \cU_d$,
 one can derive that
 \begin{align*}
  \PP_\bx\Big(\sup_{s\in(0,t]} \Vert \Delta \bX_s\Vert =0\Big)
     = \PP_\bx\big(\tau_{\cU_d} > t \big),
     \qquad t>0,\quad \bx\in\RR_+^d.
 \end{align*}
Consequently, part (iii) of Theorem \ref{Thm_jump_times_dCBI} together with the fact that
 \ $\psi^{(\cU_d)}(\blambda) = \langle \bbeta,\blambda \rangle$, \ $\blambda\in\RR_+^d$,
 yield part (ii).
\proofend

\begin{Rem}
If $(\bX_t)_{t\in\RR_+}$ is a multi-type CBI process with parameters $(d, \bc, \Bbeta, \bB, \nu, \bzero)$
 (i.e., $\bmu=\bzero$) such that the moment condition \eqref{moment_condition_m_new} holds,
 then part (i) of Proposition \ref{Pro_max_jump} yields that
 \begin{align}\label{help_Levy_supjump}
    \PP_\bx\Big(\sup_{s\in(0,t]} \Vert \Delta \bX_s\Vert > r \Big)
      = 1 - \ee^{-\nu(\overline{K}_r^c)t} = 1 - \PP\big(M( (0,t]\times \overline{K}_r^c) = 0 \big)
 \end{align}
for all $r>0$, $t>0$ and $\bx\in\RR_+^d$.
Indeed, $\bmu(\overline{K}_r^c) = \bzero$ yields that $\tbv^{(\overline{K}_r^c )}(s,\bmu(\overline{K}_r^c )) = \tbv^{(\overline{K}_r^c )}(s,\bzero) = \bzero$, $s\in\RR_+$
 (due to the uniqueness of a locally bounded solution of the system of differential equations  \eqref{help_DE_tv})
 and, by $\psi^{(\overline{K}_r^c)}(\bzero)=0$, we have the first equality in \eqref{help_Levy_supjump}.
The second equality in \eqref{help_Levy_supjump} follows from the fact that $M$ is a Poisson random measure on $\cU_1\times\cU_d$
 with intensity measure $\dd u\,\nu(\dd\br)$ (appearing in the SDE \eqref{SDE_atirasa_dimd}).
Recall that, given a real-valued L\'evy process \ $(\xi_t)_{t\in\RR_+}$ \ such that its L\'evy measure \ $\Pi$ \ satisfies
 \ $\Pi((-\infty,0])=0$, \ then
 \[
   \PP_x\Big( \sup_{s\in(0,t]} \vert \Delta\xi_s \vert > r\Big)
       = 1 - \ee^{-\Pi(\RR_+ \setminus [0,r]) t},\qquad t>0, \quad x>0,
 \]
 see, e.g., Kyprianou \cite[Exercise 2.7]{Kyp}.
This formula is in accordance with \eqref{help_Levy_supjump} with $d=1$.
Indeed, in case of $d=1$ and $\mu=0$, by the SDE \eqref{SDE_atirasa_dimd}, we have that
 \[
   \PP\Big( \sup_{s\in(0,t]} \vert \Delta X_t \vert >r \Big)
    = \PP\Big( \sup_{s\in(0,t]} \Delta Y_t >r\Big),
   \qquad t>0, \quad r>0,
 \]
 where $(Y_t)_{t\in\RR_+}$ is the L\'evy process given by
 \[
   Y_t:=\int_0^t \int_{\cU_1} r \, M(\dd u,\dd r), \qquad t\in\RR_+,
 \]
 which has a L\'evy measure $\nu$.

In particular, if $(\bX_t)_{t\in\RR_+}$ is a multi-type CBI process with parameters $(d, \bc, \Bbeta, \bB, 0, \bzero)$
 (i.e., $\nu=0$ and $\bmu=\bzero$), then \eqref{help_Levy_supjump} yields that
 \begin{align*}
    \PP_\bx\Big(\sup_{s\in(0,t]} \Vert \Delta \bX_s\Vert > r \Big) = 0, \qquad r>0,\quad t>0, \quad \bx\in\RR_+^d.
 \end{align*}
Consequently, we obtain that \ $\PP_\bx\big(\sup_{s\in(0,t]} \Vert \Delta \bX_s\Vert = 0 \big)=1$ \ for all \ $t>0$ \
 and \ $\bx\in\RR_+^d$.
\proofend
\end{Rem}

Given a measure \ $\kappa$ \ on \ $(\cU_d,\cB(\cU_d))$, \ let
 \begin{align}\label{Def_k_sup}
    \kappa_{\mathrm{sup}}
      := \begin{cases}
            \sup\Big\{ r>0 : \kappa(\overline{K}_r^c) > 0 \Big\} & \text{if \ $\kappa\ne 0$,}\\
            0 & \text{if \ $\kappa = 0$.}
         \end{cases}
 \end{align}
One may call \ $\kappa_{\mathrm{sup}}$ \ the supremum of the support of $\kappa$ (the smallest closed set whose complement has measure $0$\
 under $\kappa$).
\ Note that if \ $\kappa\ne 0$, \ then \ $\kappa_{\mathrm{sup}}\in (0,\infty]$.
\ Indeed, if \ $\kappa\ne 0$, \ then there exists a Borel set $U\in\cB(\cU_d)$ such that $\kappa(U)>0$.
\ Since \ $\overline{K}_r^c\cap U\uparrow U$ \ as \ $r\downarrow 0$, \ by the continuity from below of the measure \ $\kappa$,
 \ we get that \ $\kappa(\overline{K}_r^c\cap U)\uparrow \kappa(U)$ \ as $r\downarrow 0$.
As a consequence, there exists an $r_0>0$ such that $\kappa(\overline{K}_{r_0}^c\cap U)>0$, \ yielding that $\kappa(\overline{K}_{r_0}^c)>0$,
 and hence \ $\kappa_{\mathrm{sup}}\geq r_0>0$.

Remark also that \ $\kappa_{\mathrm{sup}} \leq \eta_{\mathrm{sup}}$ \ for measures $\kappa$ and $\eta$ on $(\cU_d,\cB(\cU_d))$ satisfying
 $\kappa \leq \eta$ (following from $\{r>0 : \kappa(\overline{K}_r^c) > 0\} \subseteq \{r>0 : \eta(\overline{K}_r^c) > 0\}$).

Next, we will give a counterpart of Corollary 4.1 in He and Li \cite{HeLi} for multi-type CBI processes.
For this, we need an auxiliary lemma in which we give a set of sufficient conditions under which
 the multi-type CBI process \ $(\bX^{(\overline{K}_r^c)}_t)_{t\in\RR_+}$ \ (for its definition, see the paragraph before Theorem \ref{Thm_jump_times_dCBI}
 by choosing \ $A:=\overline K_r^c$) \ is irreducible for sufficiently large \ $r>0$.
\ Note that in case of \ $d=1$, \ the single-type CBI process \ $(X^{(\overline{K}_r^c)}_t)_{t\in\RR_+} = (X^{((r,\infty))}_t)_{t\in\RR_+} $
 is irreducible for all $r>0$,\ since all the single-type CBI processes are irreducible.
Therefore, in the next lemma, we only consider the case $d\geq 2$, $d\in\NN$.

\begin{Lem}\label{Lem_irreducibility}
Let \ $(\bX_t)_{t\in\RR_+}$ \ be an irreducible multi-type CBI process with parameters \ $(d, \bc, \Bbeta, \bB, \nu, \bmu)$ \ such that
 \ $d\geq 2$, \ $\EE(\|\bX_0\|) < \infty$ \ and the moment condition \eqref{moment_condition_m_new} hold.
Assume that \ $\bmu\ne\bzero$, \ i.e., \ $\mu_k \ne 0$ \ for some \ $k\in\{1,\ldots,d\}$.
\ Then there exists a finite constant \ $r_0\in \big(0,\max_{k=1,\ldots,d}\,(\mu_k)_{\mathrm{sup}}\big]$ such that
 the multi-type CBI process \ $(\bX^{(\overline{K}_r^c)}_t)_{t\in\RR_+}$ \ is irreducible for all \ $r\geq r_0$.
\end{Lem}

\noindent{\bf Proof.}
{\sl Step 1.}
First, note that \ $\max_{k=1,\ldots,d}\,(\mu_k)_{\mathrm{sup}}>0$, \ which follows from \ $\bmu\ne\bzero$ \ and the discussion after \eqref{Def_k_sup}.
\ Further, for all $r>0$, the multi-type CBI process \ $(\bX^{(\overline{K}_r^c)}_t)_{t\in\RR_+}$ \ is well-defined,
 since $\nu(\overline{K}_r^c)+\sum_{k=1}^d\mu_k(\overline{K}_r^c) <\infty$ for all \ $r>0$,
 see the discussion at the beginning of the proof of Proposition \ref{Pro_max_jump}.
We need to prove that there exists a finite constant $r_0\in \big(0,\max_{k=1,\ldots,d}\,(\mu_k)_{\mathrm{sup}}\big]$ such that
 \ $\tbB^{(\overline{K}_r^c)}\in\RR_{(+)}^{d\times d}$ \ is irreducible for all $r\geq r_0$, where, due to \eqref{help21}
 (at this point, we hiddenly use that $\EE(\|\bX_0\|) < \infty$ and that the moment condition \eqref{moment_condition_m_new} hold),
 we have
 \begin{align}\label{help23}
   \tb_{i,j}^{(\overline{K}_r^c)} = b_{i,j} + \int_{\cU_d} z_i (1 - \bbone_{\overline K_r^c}(\bz)) \, \mu_j(\dd\bz)
                     = b_{i,j} + \int_{\cU_d} z_i \bbone_{\overline{K}_r}(\bz)  \, \mu_j(\dd\bz)
 \end{align}
 for \ $i\ne j$ \ with \ $i,j\in\{1,\ldots,d\}$.
\ Since \ $(\bX_t)_{t\in\RR_+}$ \ is irreducible, we have $\tbB\in\RR_{(+)}^{d\times d}$ is irreducible, where, due to \eqref{help5},
 \[
   \tb_{i,j} = b_{i,j} + \int_{\cU_d} (z_i -\delta_{i,j})^+\,\mu_j(\dd \bz),\qquad i,j\in\{1,\ldots,d\}.
 \]

{\sl  Step 2.}
Recall that for each $j\in\{1,\ldots d\}$, we have $(\mu_j)_{\mathrm{sup}}=0$ if $\mu_j=0$, and $(\mu_j)_{\mathrm{sup}}\in(0,\infty]$ if $\mu_j\ne 0$.
\ Further, if $(\mu_j)_{\mathrm{sup}}\in(0,\infty)$, then we check that \ $\mu_j\left(\overline{K}_{(\mu_j)_{\mathrm{sup}}}^c\right)=0$.
Indeed, by the definition of supremum, we have
 \ $\mu_j\left(\overline{K}_{(\mu_j)_{\mathrm{sup}}+\vare }^c\right) = 0$ \ for all $\vare>0$.
\ Using that $\overline{K}_{(\mu_j)_{\mathrm{sup}}+\vare}\downarrow \overline{K}_{(\mu_j)_{\mathrm{sup}}}$ as $\vare\downarrow 0$, we get
 $\overline{K}_{(\mu_j)_{\mathrm{sup}}+\vare}^c\uparrow \overline{K}_{(\mu_j)_{\mathrm{sup}}}^c$ as $\vare\downarrow 0$,
 and the continuity from below of the measure $\mu_j$ yields that
 \[
    0 = \mu_j\left(  \overline{K}_{(\mu_j)_{\mathrm{sup}}+\vare}^c \right)
        \uparrow \mu_j\left( \overline{K}_{(\mu_j)_{\mathrm{sup}}}^c \right)
        \qquad \text{as \ $\vare\downarrow 0$.}
 \]
Consequently, we get $ \mu_j\left( \overline{K}_{(\mu_j)_{\mathrm{sup}}}^c \right)=0$, as desired.

{\sl  Step 3.}
If $\mu_j\ne 0$ is such that $(\mu_j)_{\mathrm{sup}}=\infty$ for some $j\in\{1,\ldots,d\}$, then, since
 $\overline{K}_r\uparrow \cU_d$ as \ $r\uparrow  (\mu_j)_{\mathrm{sup}}=\infty$, \ by the continuity from below
 of the measure $\mu_j$, we get
 \begin{align}\label{help_irred_1}
    \mu_j(\{\bz \in \overline{K}_r : z_i\ne 0\})
       \uparrow \mu_j(\{\bz \in \cU_d : z_i\ne 0\})
       \qquad \text{as $r\uparrow (\mu_j)_{\mathrm{sup}}$.}
 \end{align}
If $\mu_j\ne 0$ is such that $(\mu_j)_{\mathrm{sup}}\in(0,\infty)$ for some $j\in\{1,\ldots,d\}$, then, since
 $\overline{K}_r\uparrow K_{(\mu_j)_{\mathrm{sup}}}$ as \ $r\uparrow (\mu_j)_{\mathrm{sup}}$, \
 by the continuity from below of the measure $\mu_j$, we get
 \begin{align}\label{help_irred_2}
    \mu_j(\{\bz\in \overline{K}_r : z_i\ne 0\}) \uparrow \mu_j(\{\bz\in K_{(\mu_j)_{\mathrm{sup}}} : z_i\ne 0\})
      \qquad \text{as \ $r\uparrow (\mu_j)_{\mathrm{sup}}$.}
 \end{align}

{\sl  Step 4.}
We check that there exists a finite constant \ $r_0\in\big(0,\max_{k=1,\ldots,d}\,(\mu_k)_{\mathrm{sup}}\big]$ \ such that
 if \ $\tb_{i,j}>0$ \ for some \ $i\ne j$, \ $i,j\in\{1,\ldots,d\}$, \ then
 \ $\tb_{i,j}^{(\overline{K}_r^c)}>0$ \ holds for all \ $r\geq r_0$.
\ In what follows, for each \ $i\ne j$, \ $i,j\in\{1,\ldots,d\}$ \ with \ $\tb_{i,j}>0$, \
 we define an appropriate constant $r_{i,j}$, we let $r_0$ be the maximum of these constants \ $r_{i,j}$,
 and then we check that it satisfies the property in question.

\noindent
In the remaining part of this step, let \ $i\ne j$, \ $i,j\in\{1,\ldots,d\}$ \ with
 \[
   \tb_{i,j} = b_{i,j} + \int_{\cU_d} z_i\,\mu_j(\dd \bz)>0,
 \]
 where the equality is due to \eqref{help5}.
Then \ $b_{i,j}>0$ \ or \ $\int_{\cU_d} z_i\,\mu_j(\dd \bz)>0$.

\noindent If $b_{i,j}>0$, then let $r_{i,j}:= \frac{1}{2}\big(1\wedge \max_{k=1,\ldots,d}\,(\mu_k)_{\mathrm{sup}} \big)$.

\noindent
If $b_{i,j}=0$, $\int_{\cU_d} z_i\,\mu_j(\dd \bz)>0$ and $(\mu_j)_{\mathrm{sup}}=\infty$,
 then $\mu_j\ne 0$ and $\mu_j(\{ \bz\in\cU_d : z_i\ne 0\})>0$, since otherwise
 \ $\int_{\cU_d} z_i\,\mu_j(\dd \bz) = \int_{\cU_d} z_i \bone_{\{\bz\in\cU_d : z_i\ne 0\}}\,\mu_j(\dd \bz)=0$ \ were true,
 and hence \eqref{help_irred_1} implies the existence of an $r_{i,j}\in(0,\infty)$ such that
 \ $\mu_j(\{\bz \in \overline{K}_r : z_i\ne 0\})>0$ \ for all \ $r\geq r_{i,j}$.

\noindent
If $b_{i,j}=0$, $\int_{\cU_d} z_i\,\mu_j(\dd \bz)>0$ and $(\mu_j)_{\mathrm{sup}}\in(0,\infty)$,
 then $\mu_j\ne 0$ and
 \[
  \mu_j(\{ \bz\in \overline{K}_{(\mu_j)_{\mathrm{sup}}}  : z_i\ne 0\}) = \mu_j(\{ \bz\in\cU_d : z_i\ne 0\})>0,
 \]
 where the equality follows from Step 2 and the inequality from $\int_{\cU_d}z_i\,\mu_j(\dd\bz)>0$ \ (detailed before).
Since
 \begin{align*}
   \mu_j(\{ \bz\in \overline{K}_{(\mu_j)_{\mathrm{sup}}}  : z_i\ne 0\})
      &= \mu_j(\{ \bz\in K_{(\mu_j)_{\mathrm{sup}}}  : z_i\ne 0\})\\
      &\quad  + \mu_j(\{ \bz\in\cU_d  : \Vert \bz \Vert=(\mu_j)_{\mathrm{sup}}, \, z_i\ne 0\}),
 \end{align*}
 we have $\mu_j(\{ \bz\in K_{(\mu_j)_{\mathrm{sup}}}  : z_i\ne 0\})>0$ or
 $\mu_j(\{ \bz\in\cU_d  : \Vert \bz \Vert=(\mu_j)_{\mathrm{sup}}, \, z_i\ne 0\})>0$ hold.
In case of $\mu_j(\{ \bz\in K_{(\mu_j)_{\mathrm{sup}}}  : z_i\ne 0\})>0$,
 using \eqref{help_irred_2}, we have the existence of an $r_{i,j}\in(0,(\mu_j)_{\mathrm{sup}})$ such that
 \ $\mu_j(\{\bz \in \overline{K}_r : z_i\ne 0\})>0$ \ for all \ $r\geq r_{i,j}$.
In case of \ $\mu_j(\{ \bz\in K_{(\mu_j)_{\mathrm{sup}}} : z_i\ne 0\})=0$ \ and
 $\mu_j(\{ \bz\in\cU_d  : \Vert \bz \Vert=(\mu_j)_{\mathrm{sup}}, \, z_i\ne 0\})>0$,
 let $r_{i,j}:=(\mu_j)_{\mathrm{sup}}$.
Therefore, in this case we also have that $\mu_j(\{\bz \in \overline{K}_r : z_i\ne 0\})>0$ for all $r\geq r_{i,j}$,
 since \ $\{ \bz\in\cU_d : \Vert \bz\Vert  = (\mu_j)_{\mathrm{sup}} \}\subseteq \overline{K}_r$ \ for all $r\geq r_{i,j}$.

Let us define
 \[
    r_0:=\max\big\{ r_{i,j} : i\ne j, \; i,j\in\{1,\ldots,d\}, \; \tb_{i,j}>0\big\}.
 \]
Note that \ $r_0$ \ is well-defined, since there exist $i\ne j$, $i,j\in\{1,\ldots,d\}$ \ for which $\tb_{i,j}>0$.
Indeed, otherwise $\tbB$ would be a diagonal matrix, and hence $\tbB$ would be reducible, leading us to a contradiction,
 since \ $(\bX_t)_{t\in\RR_+}$ \ is irreducible due to our assumption.

Next, we check that the finite constant $r_0$ defined above satisfies the following two properties:
 $r_0\in \big(0,\max_{k=1,\ldots,d}\,(\mu_k)_{\mathrm{sup}}\big]$
 and if \ $\tb_{i,j}>0$ \ for some \ $i\ne j$, \ $i,j\in\{1,\ldots,d\}$, \ then
 \ $\tb_{i,j}^{(\overline{K}_r^c)}>0$ \ holds for all \ $r\geq r_0$.

Since \ $(\mu_j)_{\mathrm{sup}}\leq \max_{k=1,\ldots,d}\,(\mu_k)_{\mathrm{sup}}$, \ $j\in\{1,\ldots,d\}$, \
 by the choices of $r_{i,j}$, $i\ne j$, $i,j\in\{1,\ldots,n\}$, \ we readily have that
 \ $r_0\in \big(0,\max_{k=1,\ldots,d}\,(\mu_k)_{\mathrm{sup}}\big]$.

Recall that if \ $\tb_{i,j}>0$ \ for some \ $i\ne j$, \ $i,j\in\{1,\ldots,d\}$, \ then, by \eqref{help5},
 we have \ $b_{i,j}>0$ \ or \ $\int_{\cU_d} z_i\,\mu_j(\dd \bz)>0$.

If \ $b_{i,j}>0$, \ then, by \eqref{help23}, we get \ $\tb_{i,j}^{(\overline{K}_r^c)}>0$ \ for all \ $r>0$ \ (in particular, for all $r\geq r_0$).

If $b_{i,j}=0$ and $\int_{\cU_d} z_i\,\mu_j(\dd \bz)>0$, then,
 as we already checked, \ $\mu_j(\{\bz \in \overline{K}_r : z_i\ne 0\})>0$ \ for all \ $r\geq r_0\geq r_{i,j}>0$.
This implies that $\int_{\cU_d} z_i \bone_{\overline K_r}(\bz)\,\mu_j(\dd \bz)>0$ for all $r\geq r_0$ \
 (since, otherwise, \ $\int_{\cU_d} z_i\bone_{\overline K_r}(\bz)\,\mu_j(\dd \bz) = \int_{\overline K_r} z_i \,\mu_j(\dd \bz) =0$ \ were true),
 which, as a consequence of \eqref{help23}, yields that \ $\tb_{i,j}^{(\overline{K}_r^c)}>0$ \ for all \ $r\geq r_0$.

{\sl Step 5.}
 We check that $\tbB^{(\overline{K}_r^c)}$ is irreducible for all \ $r\geq  r_0$, where $r_0$ is defined in Step 4.
On the contrary, let us assume that there exists an \ $r\geq r_0$ \ such that $\tbB^{(\overline{K}_r^c)}$
 is reducible.
Then there exist a permutation matrix \ $\bP \in \RR^{ d \times d}$ \ and an integer \ $p$ \ with \ $1 \leq p \leq d-1$ \ such that
 \[
  \bP^\top \tbB^{(\overline{K}_r^c)} \bP
   = \begin{pmatrix} \bA_1 & \bA_2 \\ \bzero & \bA_3 \end{pmatrix},
 \]
 where \ $\bA_1 \in \RR^{p\times p}$, \ $\bA_2 \in \RR^{p \times (d-p)}$, \ $\bA_3 \in \RR^{ (d-p) \times (d-p) }$, \ and \ $\bzero \in \RR^{(d-p)\times p}$ \ is
 a null matrix.
Let \ $\ell\ne m$, \ $\ell,m\in\{1,\ldots,d\}$ \ be arbitrary such that \ $(\bP^\top \tbB^{(\overline{K}_r^c)} \bP)_{\ell,m}=0$.
\ Then, since \ $\bP$ \ is a permutation matrix, there exist $i\ne j$, \ $i,j\in\{1,\ldots,d\}$ \ such that
 \ $(\bP^\top \tbB^{(\overline{K}_r^c)} \bP)_{\ell,m} = \tb_{i,j}^{(\overline{K}_r^c)}$, \ and hence
 \ $\tb_{i,j}^{(\overline{K}_r^c)}=0$.
\ Using Step 4 and that $\tb_{i,j}\geq 0$ (which holds due to $\tbB\in\RR_{(+)}^{d\times d}$), we have $\tb_{i,j} = 0$
 \ (indeed, otherwise, \ $\tb_{i,j}^{(\overline{K}_r^c)}$ \ would be positive).
\ Using again that $\bP$ is a permutation matrix, we get $(\bP^\top \tbB \bP)_{\ell,m} = \tb_{i,j}=0$.
In particular, we have $(\bP^\top \tbB \bP)_{\ell,m} =0$ for $\ell\in\{p+1,\ldots,d\}$ and $m\in\{1,\ldots,p\}$.
It implies that $\tbB$ is reducible, leading us to a contradiction.
\proofend

\begin{Pro}\label{Pro_sup_jump}
Let \ $(\bX_t)_{t\in\RR_+}$ \ be an irreducible multi-type CBI process with parameters \ $(d, \bc, \Bbeta, \bB, \nu, \bmu)$ \
 such that the moment condition \eqref{moment_condition_m_new} holds.
Let us suppose that
   \begin{itemize}
      \item[(i)] $\psi$ \ is not identically zero (i.e., \ $\bbeta\ne\bzero$ \ or \ $\nu\ne 0$),
      \item[(ii)] for each \ $i\ne j$, $i,j\in\{1,\ldots,d\}$ \ with \ $\tb_{i,j}>0$, \ $b_{i,j}=0$ and $(\mu_j)_{\mathrm{sup}}\in(0,\infty)$,
                   we have that $\mu_j(\{ \bz\in K_{(\mu_j)_{\mathrm{sup}}}  : z_i\ne 0\})>0$.
   \end{itemize}
Then for all \ $\bx\in\RR_+^d$, \ we have that
 \[
   \PP_\bx\left( \sup_{s\in(0,\infty)} \Vert \Delta\bX_s\Vert  = \Big(\nu+\sum_{i=1}^d \mu_i\Big)_{\mathrm{sup}} \right) = 1.
 \]
\end{Pro}

Before proving Proposition \ref{Pro_sup_jump}, in the next remark, we shed some light on the role of the assumption
 (ii) of Proposition \ref{Pro_sup_jump}.

\begin{Rem}
In the proof of Proposition \ref{Pro_sup_jump} in case of $d\geq 2$ and $\bmu\ne\bzero$, we will use Lemma \ref{Lem_irreducibility} in order to
 verify the existence of a finite constant $r_0\in \big(0, \max_{k=1,\ldots,d}\,(\mu_k)_{\mathrm{sup}} \big]$ such that
 the multi-type CBI process \ $\Big(\bX^{ (\overline{K}_r^c)}_t \Big)_{t\in\RR_+}$ \ is irreducible for all \ $r\geq r_0$.
In the case of $R:=\big(\nu+\sum_{i=1}^d \mu_i\big)_{\mathrm{sup}}\in(0,\infty)$, there is a subcase in the proof, where
 we will need that \ $\Big(\bX^{ (\overline{K}_{R-\vare}^c)}_t \Big)_{t\in\RR_+}$ \ is irreducible for sufficiently small $\vare>0$
 in order to be able to apply part (iii) of Corollary \ref{Cor_tau}.
For this, we need that \ $r_0<R$, \ which follows under the assumption (ii) of Proposition \ref{Pro_sup_jump}
 (see the forthcoming proof of Proposition \ref{Pro_sup_jump}).
Finally, we note that in case of $d=1$, the statement of Proposition \ref{Pro_sup_jump}
 gives back Corollary 4.1 in He and Li \cite{HeLi}, since in case of $d=1$, all the (single-type) CBI processes are
 irreducible and the assumption (ii) of Proposition \ref{Pro_sup_jump} holds automatically.
\proofend
\end{Rem}

\noindent{\bf Proof of Proposition \ref{Pro_sup_jump}.}
Recall the notation \ $R=\Big(\nu+\sum_{i=1}^d \mu_i\Big)_{\mathrm{sup}}$.
\ Let \ $\bx\in\RR_+^d$ \ be fixed arbitrarily.
First, note that
 \begin{align*}
   &\Big\{\sup_{s\in(0,t_2]} \Vert \Delta \bX_s\Vert \leq y \Big\}
           \subseteq \Big\{\sup_{s\in(0,t_1]} \Vert \Delta \bX_s\Vert \leq y \Big\},\qquad y>0, \;\; t_2\geq t_1 >0,\\
   &\Big\{\sup_{s\in(0,\infty)} \Vert \Delta \bX_s\Vert \leq y \Big\}
     = \bigcap_{t>0} \Big\{ \sup_{s\in(0,t]} \Vert \Delta \bX_s\Vert \leq y \Big\},
     \qquad y>0,
 \end{align*}
 and, by the continuity of probability, for all \ $y>0$, \ we get
 \begin{align}\label{help_sup_continuity}
  \PP_\bx\Big(\sup_{s\in(0,\infty)} \Vert \Delta \bX_s\Vert \leq y \Big)
   = \lim_{t\to\infty} \PP_\bx\Big(\sup_{s\in(0,t]} \Vert \Delta \bX_s\Vert \leq y\Big).
 \end{align}

{\sl Case I.}
Assume that \ $\nu+\sum_{i=1}^d \mu_i=0$.
\ Then \ $R=0$, \  and, by part (ii) of Corollary \ref{Cor_tau}, we have that
 \ $\PP_\bx(\tau_{\cU_d}=\infty)=1$.
Since
 \[
   \{\tau_{\cU_d}=\infty\} \subseteq \{ \Delta\bX_u = \bzero, \, \forall \; u\in\RR_{++}\}
      \subseteq \Big\{ \sup_{s\in(0,\infty)} \Vert \Delta \bX_s\Vert =0\Big\},
 \]
 this yields the statement in case of \ $\nu+\sum_{i=1}^d \mu_i=0$.

{\sl Case II.}
Assume that \ $\nu+\sum_{i=1}^d \mu_i\ne0$.
\ Then, as we discussed before Lemma \ref{Lem_irreducibility}, we have \ $R\in(0,\infty]$.

{\sl Case II/(a).}
Assume, in addition, that $R=\infty$.
Then $\left(\nu+\sum_{i=1}^d \mu_i\right)(\overline{K}_M^c)\in(0,\infty)$ \ for each \ $M\in\NN$.
Indeed, for each $M\in\NN$, the discussion at the beginning of the proof of
 Proposition \ref{Pro_max_jump} implies that $\left(\nu+\sum_{i=1}^d \mu_i\right)(\overline{K}_M^c)<\infty$,
 and, by the definition of supremum, there exists $r_M>M$ such that
 $\left(\nu+\sum_{i=1}^d \mu_i\right)(\overline{K}_{r_M}^c)>0$ and $\overline{K}_{r_M}^c\subseteq \overline{K}_M^c$.
Consequently, if $\bmu=\bzero$, then $\nu(\overline{K}_M^c)>0$ for each $M\in\NN$, and hence, by part (iii) of Corollary \ref{Cor_tau},
 we have \ $\PP_\bx(\tau_{\overline{K}_M^c} < \infty)=1$, $M\in\NN$.
If \ $\bmu\ne\bzero$, \ then we can also apply part (iii) of Corollary \ref{Cor_tau}, since in case of $d=1$,
 every (single-type) CBI process is irreducible; and in case of $d\geq 2$,
 using that \ $(\bX_t)_{t\in\RR_+}$ \ is irreducible (due to our assumption), by Lemma \ref{Lem_irreducibility},
 there exists a finite constant $r_0\in \big(0, \max_{k=1,\ldots,d}\,(\mu_k)_{\mathrm{sup}}\big]$ such that the multi-type
 the CBI process \ $\Big(\bX^{ (\overline{K}_M^c)}_t \Big)_{t\in\RR_+}$ \ is irreducible for each $M\geq r_0$, $M\in\NN$,
 and, in view of part (iii) of Corollary \ref{Cor_tau}, we get that
 \[
  \PP_\bx(\tau_{\overline{K}_M^c} < \infty)=1 \qquad \text{for each $M\geq r_0$, $M\in\NN$,}
 \]
All in all, in Case II/(a), we have that
  \[
    \PP_\bx(\tau_{\overline{K}_M^c} < \infty)=1
  \]
 for each $M\in\NN$ in case of $d=1$ or $\bmu=\bzero$, and  for each $M\geq r_0$, $M\in\NN$ in case of $d\geq 2$ and $\bmu\ne\bzero$.

In what follows, first, we handle the case $d\geq 2$ and \ $\bmu\ne\bzero$.
Then, using that \ $\{\tau_{\overline{K}_M^c} < \infty \}\subseteq \{\sup_{s\in(0,\infty)} \Vert \Delta \bX_s\Vert > M \}$, \ $M\in\NN$, \ we have
 \begin{align}\label{help57}
  \PP_\bx\Big(\sup_{s\in(0,\infty)} \Vert \Delta \bX_s\Vert > M \Big)=1
 \end{align}
 for each $M\geq r_0$, $M\in\NN$.
Since $\{\sup_{s\in(0,\infty)} \Vert \Delta \bX_s\Vert =\infty\} = \bigcap_{ M=\lfloor r_0\rfloor+1}^\infty \{ \sup_{s\in(0,\infty)} \Vert \Delta \bX_s\Vert > M\}$,
 the continuity of probability yields that
 \[
   \PP_\bx\Big( \sup_{s\in(0,\infty)} \Vert \Delta\bX_s\Vert  = \infty \Big) = 1.
 \]
The other two cases $d=1$ or $\bmu=\bzero$ can be handled in the same way (in these two cases \eqref{help57} holds for
 each $M\in\NN$).
This yields the statement in Case II/(a).

{\sl Case II/(b).}
Assume, in addition, that \ $R\in(0,\infty)$.
\ Using that $\nu\leq \nu + \sum_{i=1}^d \mu_i$ and $\mu_k \leq \nu+\sum_{i=1}^d \mu_i$, $k\in\{1,\ldots,d\}$, \ we have
  \ $\nu_{\mathrm{sup}}\leq R<\infty$ \ and \ $(\mu_k)_{\mathrm{sup}}\leq R<\infty$, $k\in\{1,\ldots,d\}$.
Further, by the definition of supremum, we have \ $\Big(\nu+\sum_{i=1}^d \mu_i\Big)(\overline{K}_{R+\vare}^c)=0$ \ for all \ $\vare>0$.
Hence, using part (i) of Proposition \ref{Pro_max_jump}, for all \ $t>0$ \ and \ $\vare>0$, we get
 \begin{align*}
   \PP_\bx\Big(\sup_{s\in(0,t]} \Vert \Delta \bX_s\Vert \leq R+\vare \Big) = \ee^0 =1,
 \end{align*}
 where we used that the solution \ $(\tbv^{(\overline{K}_{R+\vare}^c)}(t,\bzero))_{t\in\RR_+}$ \ of the system \eqref{help_DE_tv} of differential equations
 with \ $\bmu(\overline{K}_{R+\vare}^c)=\bzero$ \ is the identically \ $\bzero$ \ function and \ $\psi^{(\overline{K}_{R+\vare}^c)}(\bzero) = 0$.
Hence, by \eqref{help_sup_continuity}, for all \ $\vare>0$, \ we have
 \begin{align}\label{help19}
 \PP_\bx\Big(\sup_{s\in(0,\infty)} \Vert \Delta \bX_s\Vert \leq R + \vare \Big) = 1.
 \end{align}
Moreover, using again the definition of supremum and Remark \ref{Rem_A_condition}, one can see that
 \begin{align}\label{help24}
    \Big(\nu+\sum_{i=1}^d \mu_i\Big)( \overline{K}_{R-\vare}^c )\in(0,\infty), \qquad \vare\in(0,R).
 \end{align}
Indeed, by the definition of supremum, for all $\vare\in(0,R)$, there exists \ $R_\vare\in(R-\vare,R)$ \
 such that \ $\big(\nu+\sum_{i=1}^d \mu_i\big)(\overline{K}_{R_\vare}^c)>0$, \
 and, since \ $\overline{K}_{R_\vare}^c\subseteq \overline{K}_{R - \vare}^c$, \ it yields that \ $\big(\nu+\sum_{i=1}^d \mu_i\big)( \overline{K}_{R-\vare}^c) >0$.
Further, the discussion at the beginning of the proof of Proposition \ref{Pro_max_jump} implies
 that $\big(\nu+\sum_{i=1}^d \mu_i\big)( \overline{K}_{R-\vare}^c)< \infty$.

\noindent If \ $\bmu=\bzero$, \ then, as a consequence of \eqref{help24}, \ $\nu(\overline{K}_{R-\vare}^c)\in(0,\infty)$, $\vare\in(0,R)$, and hence,
 using part (iii) of Corollary \ref{Cor_tau}, for all $\vare\in(0,R)$, we have
 \begin{align}\label{help35}
  \PP_\bx\big( \tau_{ \overline{K}_{R-\vare}^c } <\infty \big) = 1.
 \end{align}
If $d=1$ and $\bmu\ne\bzero$, then using that every single-type CBI process is irreducible,
 part (iii) of Corollary \ref{Cor_tau} implies that \eqref{help35} holds  for all $\vare\in(0,R)$.

\noindent If \ $d\geq 2$ \ and \ $\bmu\ne\bzero$, \ then, using that \ $(\bX_t)_{t\in\RR_+}$ \ is irreducible (due to our assumption),
 by Lemma \ref{Lem_irreducibility},
 there exists a finite constant $r_0\in \big(0, \max_{k=1,\ldots,d}\,(\mu_k)_{\mathrm{sup}} \big]  \subseteq (0,R]$ such that
  the multi-type CBI process \ $\Big(\bX^{ (\overline{K}_r^c)}_t \Big)_{t\in\RR_+}$ \  is irreducible for all $r\geq r_0$.
Next, we check that in this case we have \ $r_0<R$.
By Step 4 in the proof of Lemma \ref{Lem_irreducibility} and the assumption (ii), for each $i\ne j$, $i,j\in\{1,\ldots,d\}$
 with $\tb_{i,j}>0$ we have that $r_{i,j}= \frac{1}{2}\big(1\wedge \max_{k=1,\ldots,d}\,(\mu_k)_{\mathrm{sup}}\big)$ or $r_{i,j}\in(0,(\mu_j)_{\mathrm{sup}})$,
 and hence $r_{i,j}\in(0,\max_{k=1,\ldots,d}\,(\mu_k)_{\mathrm{sup}})$.
This yields that \ $r_0\in(0,\max_{k=1,\ldots,d}\,(\mu_k)_{\mathrm{sup}})$ (following from the definition of $r_0$ in Step 4 in the proof of Lemma \ref{Lem_irreducibility}),
 and consequently,
 \[
   0 < r_0<\max_{k=1,\ldots,d}\,(\mu_k)_{\mathrm{sup}} \leq \Big(\nu+\sum_{i=1}^d \mu_i\Big)_{\mathrm{sup}}=R.
 \]
Therefore, we have $r_0<R$, as desired.
Since $r_0<R$, using \eqref{help24}, for all $\vare\in(0,R-r_0)$, we have that \ $\big(\nu+\sum_{i=1}^d \mu_i\big)( \overline{K}_{R-\vare}^c)\in(0,\infty)$
 \ and \ $\Big(\bX^{ (\overline{K}_{R-\vare}^c)}_t \Big)_{t\in\RR_+}$ \ is irreducible.
Consequently, using again part (iii) of Corollary \ref{Cor_tau}, we get that \eqref{help35} holds
 for all $\vare\in(0,R-r_0)$.

\noindent All in all, in Case II/(b), we get that
 \begin{align*}
  1 = \PP_\bx\big( \tau_{ \overline{K}_{R-\vare}^c } <\infty \big)
    \leq \PP_\bx\left( \sup_{s\in(0,\infty)} \Vert \Delta\bX_s\Vert > R-\vare \right)
 \end{align*}
 for all $\vare\in(0,R)$ in case of $d=1$ or $\bmu=\bzero$, and for all $\vare\in(0,R-r_0)$ in case of $d\geq 2$ and $\bmu\ne\bzero$.

In what follows, we handle the case $d\geq 2$ and \ $\bmu\ne\bzero$, the other two cases
 $d=1$ or $\bmu=\bzero$ can be handled in the same way (replace $\vare\in(0,R-r_0)$ by $\vare\in(0,R)$ in the following argument).
Then for all $\vare\in(0,R-r_0)$, we have
 \begin{align}\label{help20}
   \PP_\bx\Big( \sup_{s\in(0,\infty)} \Vert \Delta\bX_s\Vert > R - \vare \Big)=1.
 \end{align}
Using \eqref{help19} and \eqref{help20}, for all $\vare\in(0,R-r_0)$, we get
 \[
    \PP_\bx\Big( \sup_{s\in(0,\infty)} \Vert \Delta\bX_s\Vert \in (R - \vare, R+\vare] \Big)=1.
 \]
By taking the limit as \ $\vare\downarrow 0$, \ the continuity of probability implies that
  \[
    \PP_\bx\Big( \sup_{s\in(0,\infty)} \Vert \Delta\bX_s\Vert = R \Big)=1.
 \]
This yields the statement in Case II/(b) as well.
\proofend

In what follows, we will study the relationship
 between the distributional properties of jumps for a multi-type CBI process \ $(\bX_t)_{t\in\RR_+}$ \ and
 its total L\'evy measure \ $\nu+\sum_{i=1}^d\mu_i$ \ on \ $(\cU_d, \cB(\cU_d))$.

For a c\`{a}dl\`{a}g,  $\RR^d$-valued stochastic process \ $(\bxi_t)_{t\in\RR_+}$, \ let us introduce the notation
 \[
    \sup_{s\in(0,t]} \bxi_s := \Bigg( \sup_{s\in(0,t]} \xi_{s,1},\,\ldots, \sup_{s\in(0,t]} \xi_{s,d} \Bigg)^\top \in\RR^d, \qquad t>0,
 \]
 where \ $\bxi_t=(\xi_{t,1},\ldots,\xi_{t,d})^\top$, $t\in\RR_+$.
Recall the notation
 \begin{align*}
 \cR_d=\left\{ \left(\prod_{i=1}^d [0, w_i] \right) \setminus \{\bzero\}:  w_1,\ldots,w_d\in \RR_{++} \right\} .
 \end{align*}

Given a multi-type CBI process \ $(\bX_t)_{t\in\RR_+}$, \ recall that
 for all $t>0$ and \ $\bx\in\RR_+^d$, \ the probability measure \ $\pi_{t,\bx}$ \ on \ $(\RR_+^d,\cB(\RR_+^d))$ \ is given by
 \[
  \pi_{t,\bx}(A) = \PP_\bx\Bigg(\sup_{s\in(0,t]} \Delta \bX_s\in A \Bigg), \qquad A\in\cB(\RR_+^d).
 \]

The forthcoming Theorem \ref{Thm_sup_jump_equiv} can be considered as a multi-type counterpart
 of Theorem 4.2 in He and Li \cite{HeLi},
 which is about single-type CBI processes.
Theorem 4.2 in He and Li \cite{HeLi} is also contained as Theorem 10.23 in Li \cite{Li}.
In the proof of part (ii) of Theorem \ref{Thm_sup_jump_equiv},
 we use quite different arguments from those used in the proof of Theorem 4.2 in He and Li \cite{HeLi}.

\begin{Thm}\label{Thm_sup_jump_equiv}
Let \ $(\bX_t)_{t\in\RR_+}$ \ be a multi-type CBI process with parameters \ $(d, \bc, \Bbeta, \bB, \nu, \bmu)$
 \ such that the moment condition \eqref{moment_condition_m_new} holds.
 \begin{itemize}
   \item[(i)] If \ $\left( \nu+\sum_{i=1}^d\mu_i \right)(A)=0$ \ with some \ $A\in\cR_d$, \
              then \ $\pi_{t,\bx}(A)=0$ \ for all \ $t>0$ \ and \ $\bx\in\RR_+^d$.
   \item[(ii)] If \ $\pi_{t,\bx}(A)=0$ \ for some \ $t>0$, \ $\bx\in\RR_{++}^d$ and \ $A\in\cR_d$, \
               then \ $\left( \nu+\sum_{i=1}^d\mu_i \right)(A)=0$.
               Furthermore, in case of \ $\Bbeta\in\RR_{++}^d$, \
               we can extend it to $\bx\in\RR_+^d$, i.e.,
               if \ $\pi_{t,\bx}(A)=0$ \ for some \ $t>0$, \ $\bx\in\RR_+^d$ and \ $A\in\cR_d$, \
               then \ $\left( \nu+\sum_{i=1}^d\mu_i \right)(A)=0$.
 \end{itemize}
\end{Thm}

\noindent{\bf Proof.}
(i). Let $t>0$, $\bx\in\RR_+^d$, and $A\in\cR_d$ be such that $\nu(A)+\sum_{i=1}^d \mu_i(A)=0$.
We check that
 \begin{align}\label{help44}
   \Big\{ \sup_{s\in(0,t]} \Delta \bX_s\in A \Big\}\subseteq \{\tau_A\leq t\}.
 \end{align}
If $\sup_{s\in(0,t]} \Delta \bX_s\in A$ holds,
 then there exists $\ba=(a_1,\ldots,a_d)^\top\in A$ such that $\sup_{s\in(0,t]}\Delta\bX_s=\ba$.
Then, by the definition of $\sup_{s\in(0,t]}\Delta\bX_s$, we get that $\sup_{s\in(0,t]}\Delta X_{s,i} = a_i$, $i=1,\ldots,d$.
Since $\ba\in A$ and $\bzero\notin A$, we have that $a_i\geq 0$, $i\in \{1,\ldots,d\}$, and there exists $i_0\in\{1,\ldots,d\}$
 such that $a_{i_0}>0$.
Further, by the definition of supremum of a set of real numbers, there exists a sequence $(s_n^{(i_0)})_{n\in\NN}$ in $(0,t]$
 such that $\Delta X_{s_n^{(i_0)},i_0}\uparrow a_{i_0}$ as $n\to\infty$, and hence, since $a_{i_0}>0$, we also have
  $\Delta X_{s_n^{(i_0)},i_0}>0$ (in particular, $\Delta \bX_{s_n^{(i_0)}}\ne \bzero$) for sufficiently large $n\in\NN$.
Using that $\sup_{s\in(0,t]}\Delta\bX_s=\ba$, we get that $\Delta \bX_{s_n^{(i_0)}} \leq \ba$, $n\in\NN$.
Consequently, since $A\in\cR_d$, it yields that $\Delta \bX_{s_n^{(i_0)}}\in A$ for sufficiently large $n\in\NN$.
Hence, we obtain that $J_t(A)\geq 1$ and $\tau_A\leq t$, yielding \eqref{help44}.

Consequently, part (iii) of Theorem \ref{Thm_jump_times_dCBI} implies that for all $t>0$ and $\bx\in\RR_+^d$, we have that
 \begin{align*}
   \pi_{t,\bx}(A)
      \leq \PP_\bx(\tau_A\leq t) = 1 - \PP_\bx(\tau_A>t) = 1- \ee^0 =0,
 \end{align*}
 since \ $\bmu(A)=0$ \ implies that \ $\tbv^{(A)}(t,\bmu(A))=\bzero$, \ $t\in\RR_+$, \ and \ $\psi^{(A)}(\bzero)=0$.
This yields the assertion of part (i).

(ii).
We will prove it by contradiction.
Let \ $t>0$ \ and \ $A\in\cR_d$ \ be such that \ $\pi_{t,\bx}(A)=0$, \ where \ $\bx\in\RR_{++}^d$ \ or \ $\Bbeta\in\RR_{++}^d$ \
 (and $\bx\in\RR_+^d$), \ and, on the contrary, let us assume that \ $\big(\nu+\sum_{\ell=1}^d\mu_\ell\big)(A)>0$.
\ We check that there exists an $r_0>0$ such that $K_r\setminus \{\bzero\} \subseteq A$ and \ $\big(\nu+\sum_{\ell=1}^d\mu_\ell\big)(A\cap K_r^c)\in(0,\infty)$
 for all $r\in(0,r_0)$.
First, note that if $A$ has  the form $\left(\prod_{i=1}^d [0, w_i]\right)\setminus\{\bzero\}$
 \ with some $w_1,\ldots,w_d\in \RR_{++}$, then $K_{\frac{1}{2}(w_1\wedge \cdots\wedge w_d)}\setminus \{\bzero\} \subseteq A$.
\ Further, by Remark \ref{Rem_A_condition}, we have
 \[
     \Big(\nu+\sum_{\ell=1}^d\mu_\ell\Big)(A\cap K_r^c)
        \leq \Big(\nu+\sum_{\ell=1}^d\mu_\ell\Big)(K_r^c)<\infty \qquad \text{for all $r>0$.}
 \]
Moreover, since $A$ is a nondegenerate rectangle in $\RR_+^d$ anchored at $\bzero$, we have
 that $A = \bigcup_{r>0} (A\cap K_r^c)$, and hence the continuity below of the measure \ $\nu+\sum_{\ell=1}^d\mu_\ell$ \ implies that
 \[
     \lim_{r\downarrow 0} \Big(\nu+\sum_{\ell=1}^d\mu_\ell\Big)(A\cap K_r^c) = \Big(\nu+\sum_{\ell=1}^d\mu_\ell\Big)(A)>0.
 \]
This yields the existence of an $r_0\in(0,\frac{1}{2}(w_1\wedge \cdots\wedge w_d))$ such that $\big(\nu+\sum_{\ell=1}^d\mu_\ell\big)(A\cap K_r^c) \in(0,\infty)$
 for all $r\in(0,r_0)$, as desired.

Using part (i) of Corollary \ref{Cor_tau}, \ $K_r^c\setminus A\subseteq K_r^c$,
 \ and that \ $\big(\sum_{\ell=1}^d\mu_\ell\big)(K_r^c)<\infty$ (see Remark \ref{Rem_A_condition}), we have
 \begin{align}\label{help36}
    \tbv^{(K_r^c\setminus A)}(t,\bmu(K_r^c\setminus A)) \leq \tbv^{(K_r^c)}(t,\bmu(K_r^c)).
 \end{align}
 By \eqref{formula3}, we obtain that
 \begin{align*}
   \PP_\bx(\tau_{K_r^c}>t)
           = \exp\left\{ -\nu(K_r^c)t - \sum_{\ell=1}^d x_\ell\,\tv_\ell^{(K_r^c)}(t,\bmu(K_r^c))
                           - \int_0^t \psi^{(K_r^c)}\big(\tbv^{(K_r^c)}(s,\bmu(K_r^c))\big)\,\dd s \right\}.
 \end{align*}
Here for all \ $\blambda\in\RR_+^d$, \ by the definitions of \ $\psi^{(K_r^c)}$ \ and \ $\psi^{(K_r^c\setminus A)}$ \ (see \eqref{formula_psi_A}),
 we have that
 \begin{align*}
   \psi^{(K_r^c)}(\blambda)
   & = \psi(\blambda) - \int_{K_r^c} \bigl( 1 - \ee^{- \langle\blambda, \br\rangle} \bigr) \,\nu(\dd \br) \\
   & = \psi(\blambda) - \int_{K_r^c\setminus A} \bigl( 1 - \ee^{- \langle\blambda, \br\rangle} \bigr) \,\nu(\dd \br)
      - \int_{K_r^c\cap A} \bigl( 1 - \ee^{- \langle\blambda, \br\rangle} \bigr) \,\nu(\dd \br) \\
   & = \psi^{(K_r^c\setminus A)}(\blambda) - \int_{K_r^c\cap A} \bigl( 1 - \ee^{- \langle\blambda, \br\rangle} \bigr) \,\nu(\dd \br),
 \end{align*}
 and hence
  \begin{align}\label{help37}
   \begin{split}
         &\PP_\bx(\tau_{K_r^c}>t)\\
         &\qquad  = \exp\Bigg\{ - \nu(K_r^c\setminus A)t - \sum_{\ell=1}^d x_\ell\,\tv_\ell^{(K_r^c)}(t,\bmu(K_r^c))
                           - \int_0^t \psi^{(K_r^c\setminus A)}\big(\tbv^{(K_r^c)}(s,\bmu(K_r^c))\big)\,\dd s \\
        &\phantom{\qquad  = \exp\Bigg\{ }
                           - \nu(K_r^c\cap A)t
                           + \int_0^t \left( \int_{K_r^c\cap A}  \Big( 1 -  \ee^{- \big\langle \tbv^{(K_r^c)}(s,\bmu(K_r^c)) , \br\big\rangle} \Big)\nu(\dd\br) \right)\dd s
                           \Bigg\} =
  \end{split}
  \end{align}
 \begin{align*}
        &\qquad = \exp\Bigg\{ - \nu(K_r^c\setminus A)t - \sum_{\ell=1}^d x_\ell\,\tv_\ell^{(K_r^c)}(t,\bmu(K_r^c))
                           - \int_0^t \psi^{(K_r^c\setminus A)}\big(\tbv^{(K_r^c)}(s,\bmu(K_r^c))\big)\,\dd s\\
        &\phantom{\qquad = \exp\Bigg\{}
                           - \int_0^t \left( \int_{K_r^c\cap A}  \ee^{- \big\langle \tbv^{(K_r^c)}(s,\bmu(K_r^c)) , \br\big\rangle}\, \nu(\dd\br) \right)\dd s
                           \Bigg\},
 \end{align*}
 where the second equality follows from \ $\int_0^t\left(\int_{K_r^c\cap A} 1\,\nu(\dd r)\right)\dd s = t\nu(K_r^c\cap A)<\infty$.

{\sl Case I:} First, we consider the case $\nu(A\cap K_r^c) >0$.
Using \eqref{formula_psi_A} and that $\bbeta\in\RR_+^d$, we have $\psi^{(K_r^c\setminus A)}$ is monotone increasing.
This together with \eqref{help36} and \eqref{help37} yield that
 \begin{align*}
  &\PP_\bx(\tau_{K_r^c}>t)\\
  &< \exp\left\{ - \nu(K_r^c\setminus A)t - \sum_{\ell=1}^d x_\ell\,\tv_\ell^{(K_r^c)}(t,\bmu(K_r^c))
                           - \int_0^t \psi^{(K_r^c\setminus A)}\big(\tbv^{(K_r^c)}(s,\bmu(K_r^c))\big)\,\dd s \right\}\\
  & \leq \exp\left\{ - \nu(K_r^c\setminus A)t - \sum_{\ell=1}^d x_\ell\,\tv_\ell^{(K_r^c\setminus A)}(t,\bmu(K_r^c\setminus A))
                           - \int_0^t \psi^{(K_r^c\setminus A)}\big(\tbv^{(K_r^c\setminus A)}(s,\bmu(K_r^c\setminus A))\big)\,\dd s \right\},
 \end{align*}
 where the strict inequality follows from \ $\nu(A\cap K_r^c) >0$ \
 and \ $\ee^{- \big\langle \tbv^{(K_r^c)}(s,\bmu(K_r^c)) , \br\big\rangle}\in(0,1]$ \ for \ $s>0$
 \ (following from \ $\sum_{\ell=1}^d \mu_\ell(K_r^c)<\infty$ and $\tbv^{(K_r^c)}(s,\bmu(K_r^c))\in\RR_+^d$, $s\in\RR_+$).
Using again \eqref{formula3}, the right hand side of the above inequality coincides with
 \ $\PP_\bx(\tau_{K_r^c\setminus A}>t)$, \ and hence we have
 \begin{align}\label{help38}
    \PP_\bx(\tau_{K_r^c}>t) < \PP_\bx(\tau_{K_r^c\setminus A}>t).
 \end{align}
Furthermore, we get
 \begin{align}\label{help42}
  \begin{split}
   \Big\{ \sup_{s\in(0,t]} \Delta \bX_s\in A \cap K_r^c  \Big\}
       & \supseteq \{ \tau_{K_r^c\cap A}\leq t\} \cap \{ \tau_{K_r^c\setminus A} > t\}
                    \cap \{  J_{t+1}(K_r^c\cap A)<\infty \} \\
       &= \Big(\{ \tau_{K_r^c}\leq t\} \setminus \{ \tau_{K_r^c\setminus A} \leq t\}\Big)
                    \cap \{  J_{t+1}(K_r^c\cap A)<\infty \},
  \end{split}
 \end{align}
 where the equality follows from the (not necessarily disjoint) decomposition
  \[
   \{ \tau_{K_r^c} \leq t \}
     = \{ \tau_{K_r^c\cap A} \leq t \} \cup  \{ \tau_{K_r^c\setminus A} \leq t \},
  \]
 and the inclusion can be checked as follows.
If \ $J_{t+1}(K_r^c\cap A)<\infty$, \ then $(\bX_s)_{s\in\RR_+}$ has at most a finite number of jumps of which the size vectors
 belong to the set $K_r^c\cap A$ during the time interval \ $(0,t+1]$.
If, in addition, $\tau_{K_r^c\cap A}\leq t$ and $\tau_{K_r^c\setminus A} > t$ hold,
 then, on the one hand, $(\bX_s)_{s\in\RR_+}$ has at least one jump, at time point $\tau_{K_r^c\cap A}$,
 of which the size vector belongs to $K_r^c \cap A$ during the time interval \ $(0,t]$; \
 and, on the other hand, there is no jump of which the size vector belongs to $K_r^c\setminus A = K_r^c\cap A^c = (K_r\cup A)^c = (A\cup\{\bzero\})^c$
 (due to $K_r\setminus \{\bzero\} \subseteq A$) during the time interval \ $(0,t]$.
Consequently, on the one hand, we get \ $\sup_{s\in(0,t]} \Delta \bX_s\in K_r^c$,
 \ since, with the notation \ $\by:= \Delta\bX_{\tau_{K_r^c\cap A}}$ \ and using that \ $\tau_{K_r^c\cap A}\in(0,t]$,
 \ we have \ $\sup_{s\in(0,t]} \Delta X_{s,i} \geq y_i\geq 0$, \ $i=1,\ldots,d$, \ and hence
 \[
  \left\Vert \sup_{s\in(0,t]}  \Delta\bX_s\right\Vert
    = \left( \sum_{i=1}^d \Big( \sup_{s\in(0,t]} \Delta X_{s,i} \Big)^2 \right)^{1/2}
    \geq \left( \sum_{i=1}^d y_i^2\right)^{1/2} = \Vert \by\Vert \geq r.
 \]
On the other hand, taking into account that $A$ has the form $\left(\prod_{i=1}^d [0, w_i]\right)\setminus\{\bzero\}$
 \ with some $w_1,\ldots,w_d\in \RR_{++}$, we get \ $\sup_{s\in(0,t]} \Delta \bX_s\in A$, yielding the inclusion in \eqref{help42}.

Since $\big(\nu+\sum_{\ell=1}^d\mu_\ell\big)(A\cap K_r^c)<\infty$, by part (i) of Theorem \ref{Thm_jump_times_dCBI}, we have
 \[
   \PP_\bx(J_s(K_r^c\cap A)<\infty)=1,\qquad s>0,
 \]
 and consequently, by \eqref{help42}, we obtain that
 \begin{align}\label{help41}
  \begin{split}
   \pi_{t,\bx}(K_r^c\cap A)
     &\geq \PP_\bx\big(\{ \tau_{K_r^c}\leq t\} \setminus \{ \tau_{K_r^c\setminus A} \leq t\} \big) \\
     &= \PP_\bx( \tau_{K_r^c}\leq t)  -  \PP_\bx( \tau_{K_r^c\setminus A}\leq t) \\
     & = 1- \PP_\bx( \tau_{K_r^c}> t)  -  \big(1 - \PP_\bx( \tau_{K_r^c\setminus A} >  t) \big) \\
     & = \PP_\bx( \tau_{K_r^c\setminus A}> t) - \PP_\bx( \tau_{K_r^c}> t),
  \end{split}
 \end{align}
 where the second equality follows from \ $\tau_{K_r^c}\leq \tau_{K_r^c\setminus A}$ \ yielding that
 \ $\{ \tau_{K_r^c\setminus A} \leq t\} \subseteq \{ \tau_{K_r^c} \leq t\}$.
\ Hence, by \eqref{help38}, we get \ $\pi_{t,\bx}(K_r^c\cap A)>0$.
\ This is a contradiction, since the assumption  \ $\pi_{t,\bx}(A)=0$ \ yields that $\pi_{t,\bx}(K_r^c\cap A)=0$.

 \smallskip

{\sl Case II:}
Next, we consider the case \ $\nu(A\cap K_r^c) =0$.
\ Then, since \ $\big(\nu+\sum_{\ell=1}^d\mu_\ell\big)(A\cap K_r^c)\in(0,\infty)$, \
 we must have \ $\sum_{\ell=1}^d\mu_\ell(A\cap K_r^c)\in(0,\infty)$.
\ Hence we obtain that
 \begin{align*}
  \sum_{\ell=1}^d\mu_\ell(K_r^c)
    = \sum_{\ell=1}^d\mu_\ell(K_r^c\cap A) + \sum_{\ell=1}^d\mu_\ell(K_r^c\setminus A)
    >  \sum_{\ell=1}^d\mu_\ell(K_r^c\setminus A).
 \end{align*}
We are going to prove that for all \ $\bz=(z_1,\ldots,z_d)\in\RR_{++}^d$, \ there exists a sufficiently small $t_0\in(0,t)$
 (may depend on \ $\bz$, $t$, $r$ \ and $A$) \ such that
 \begin{align}\label{help39}
   \sum_{\ell=1}^d z_\ell \tv_\ell^{(K_r^c)}(s,\bmu(K_r^c)) > \sum_{\ell=1}^d z_\ell \tv_\ell^{(K_r^c\setminus A)}(s,\bmu(K_r^c\setminus A)),
      \qquad s\in(0,t_0].
 \end{align}
By \eqref{help_DE_tv}, we have that
 \begin{align}\label{help40}
  \begin{split}
  &\sum_{\ell=1}^d
           z_\ell \Big(\partial_1 \tv_\ell^{(K_r^c)}(s, \bmu(K_r^c)) -  \partial_1 \tv_\ell^{(K_r^c\setminus A)}(s, \bmu(K_r^c\setminus A)) \Big)\\
    &=\sum_{\ell=1}^d z_\ell\Big(\mu_\ell(K_r^c)  - \mu_\ell(K_r^c\setminus A)
                                   - \varphi_\ell^{(K_r^c)}( \tbv^{(K_r^c)}(s, \bmu(K_r^c)))\\
    &\phantom{=\sum_{\ell=1}^d z_\ell\Big(\;} + \varphi_\ell^{(K_r^c\setminus A)}( \tbv^{(K_r^c\setminus A)}(s, \bmu(K_r^c\setminus A)) ) \Big) =
  \end{split}
 \end{align}
 \begin{align*}
  \begin{split}
           = \sum_{\ell=1}^d z_\ell \mu_\ell(K_r^c\cap A)
             - \sum_{\ell=1}^d z_\ell \Big( \varphi_\ell^{(K_r^c)}( \tbv^{(K_r^c)}(s, \bmu(K_r^c)))
                                      - \varphi_\ell^{(K_r^c\setminus A)}( \tbv^{(K_r^c\setminus A)}(s, \bmu(K_r^c\setminus A)) ) \Big)
  \end{split}
 \end{align*}
 for $s\in\RR_{++}$ and
 \begin{align}\label{help_s_to_zero}
    \tbv^{(K_r^c)}(0, \bmu(K_r^c)) = \bzero = \tbv^{(K_r^c\setminus A)}(0, \bmu(K_r^c\setminus A)).
 \end{align}
Using that $\varphi_\ell^{(K_r^c)}(\bzero) = \varphi_\ell^{(K_r^c\setminus A)}(\bzero) = 0$, $\ell\in\{1,\ldots,d\}$,
 \eqref{help_s_to_zero} yields that
 \begin{align}\label{help45}
  \sum_{\ell=1}^d z_\ell \Big( \varphi_\ell^{(K_r^c)}( \tbv^{(K_r^c)}(0, \bmu(K_r^c)))
                                      - \varphi_\ell^{(K_r^c\setminus A)}( \tbv^{(K_r^c\setminus A)}(0, \bmu(K_r^c\setminus A)) ) \Big)=0.
 \end{align}
By part (iii) of Theorem \ref{Thm_jump_times_dCBI} and the continuity of $\bvarphi$ (see the proof of part (iv) in Proposition \ref{Pro_tv}),
 for all \ $S\in\cB(\cU_d)$ \ with $\sum_{\ell=1}^d \mu_\ell(S)<\infty$, we get that the functions
 \ $\RR_+^d \ni\blambda\mapsto \bvarphi^{(S)}(\blambda)$ \ and
 \ $\RR_+\ni t\mapsto\tbv^{(S)}(t,\bmu(S))$ \ are continuous.
Hence, using also \eqref{help45},
 we have that
 \[
  \lim_{s\downarrow 0} \Bigg[\sum_{\ell=1}^d z_\ell \Big( \varphi_\ell^{(K_r^c)}( \tbv^{(K_r^c)}(s, \bmu(K_r^c)))
                                      - \varphi_\ell^{(K_r^c\setminus A)}( \tbv^{(K_r^c\setminus A)}(s, \bmu(K_r^c\setminus A)) ) \Big) \Bigg] =0.
 \]
Since \ $\sum_{\ell=1}^d \mu_\ell(K_r^c\cap A) \in (0,\infty)$ \ and \ $\bz\in\RR_{++}^d$, \
 we have $\sum_{\ell=1}^d z_\ell \mu_\ell(K_r^c\cap A) \in(0,\infty)$.
Consequently, the limit of the right hand side of \eqref{help40} as $s\downarrow 0$ is
 \ $\sum_{\ell=1}^d z_\ell \mu_\ell(K_r^c\cap A)\in(0,\infty)$.
\ Hence, using that $t>0$ and the definition of the limit of a function at a point,
 there exists a sufficiently small \ $t_0\in(0,t)$ \ such that the right hand side of
 \eqref{help40} is positive for all $s\in[0,t_0]$.
Consequently, using that the right hand side of \eqref{help40} as a function of $s$ is continuous, its integral
 on $[0,s]$ is positive for all $s\in(0,t_0]$.
\ Then, integrating the left and right hand sides of \eqref{help40} and using \eqref{help_s_to_zero}, we get
 \begin{align*}
 &\sum_{\ell=1}^d
           z_\ell \Big(\tv_\ell^{(K_r^c)}(s, \bmu(K_r^c)) -  \tv_\ell^{(K_r^c\setminus A)}(s, \bmu(K_r^c\setminus A)) \Big) \\
 &= \int_0^s  z_\ell \sum_{\ell=1}^d  \Bigg[ \mu_\ell(K_r^c\cap A)
                   -  \Big( \varphi_\ell^{(K_r^c)}( \tbv^{(K_r^c)}(u, \bmu(K_r^c)))
                                         - \varphi_\ell^{(K_r^c\setminus A)}( \tbv^{(K_r^c\setminus A)}(u, \bmu(K_r^c\setminus A)) ) \Big) \Bigg]\dd u
 \end{align*}
 is positive for all $s\in(0,t_0]$, implying \eqref{help39}, as desired.

Using \eqref{formula3} (similarly, as we derived \eqref{help37} for which we did not use that $\pi_{t,\bx}(A)=0$)
 together with $\nu(A\cap K_r^c)=0$, we can get that
  \begin{align}\label{help43}
   \begin{split}
        &\PP_\bx(\tau_{K_r^c}>t_0)\\
        & = \exp\Big\{ - \nu(K_r^c\setminus A)t_0 - \sum_{\ell=1}^d x_\ell\,\tv_\ell^{(K_r^c)}(t_0,\bmu(K_r^c))
                           - \int_0^{t_0} \psi^{(K_r^c\setminus A)}\big(\tbv^{(K_r^c)}(s,\bmu(K_r^c))\big)\,\dd s\Big\}<
   \end{split}
  \end{align}
  \begin{align*}
      & < \exp\Big\{ \!\!- \nu(K_r^c\setminus A)t_0 - \! \sum_{\ell=1}^d x_\ell\,\tv_\ell^{(K_r^c\setminus A)}(t_0,\bmu(K_r^c\setminus A))\\
        &\phantom{<\exp\Big\{}
           - \int_0^{t_0} \psi^{(K_r^c\setminus A)}\big(\tbv^{(K_r^c\setminus A)}(s,\bmu(K_r^c\setminus A))\big)\,\dd s\Big\}\\
        &= \PP_\bx(\tau_{K_r^c\setminus A}>t_0),
 \end{align*}
 where, at the second step, the strict inequality can be checked as follows.
In case of $\bx\in\RR_{++}^d$, it directly follows from \eqref{help39} together with the facts
 that $\tbv^{(K_r^c\setminus A)}(s,\bmu(K_r^c\setminus A))\leq \tbv^{(K_r^c)}(s,\bmu(K_r^c))$, $s\in[0,t_0]$
 \ (due to part (i) of Corollary \ref{Cor_tau}, \ $K_r^c\setminus A\subseteq K_r^c$,
 \ and \ $\big(\sum_{\ell=1}^d\mu_\ell\big)(K_r^c)<\infty$) and that \ $\psi^{(K_r^c\setminus A)}$ \ is monotone increasing.
In case of $\Bbeta\in\RR_{++}^d$, using \eqref{formula_psi_A},  we have
\begin{align*}
  &\int_0^{t_0} \psi^{(K_r^c\setminus A)}\big(\tbv^{(K_r^c)}(s,\bmu(K_r^c))\big)\,\dd s \\
  & = \int_0^{t_0} \big\langle \Bbeta, \tbv^{(K_r^c)}(s,\bmu(K_r^c)) \big\rangle \,\dd s
     + \int_0^{t_0} \int_{\cU_d\setminus (K_r^c\setminus A)} (1-\ee^{-\langle \tbv^{(K_r^c)}(s,\bmu(K_r^c)) , \br\rangle })\,\nu(\dd\br)\,\dd s.
 \end{align*}
Using $\Bbeta\in\RR_{++}^d$, \eqref{help39} and the continuity of the functions $\RR_+\ni s\mapsto \tbv^{(K_r^c)}(s,\bmu(K_r^c))$ and
 $\RR_+\ni s\mapsto \tbv^{(K_r^c\setminus A)}(s,\bmu(K_r^c\setminus A))$, we have
 \begin{align*}
  &\int_0^{t_0} \big\langle \Bbeta, \tbv^{(K_r^c)}(s,\bmu(K_r^c)) \big\rangle \,\dd s
    - \int_0^{t_0} \big\langle \Bbeta, \tbv^{(K_r^c\setminus A)}(s,\bmu(K_r^c\setminus A)) \big\rangle \,\dd s \\
  &\qquad = \int_0^{t_0} \left(  \sum_{\ell=1}^d \beta_\ell\Big( \tv_\ell^{(K_r^c)}(s,\bmu(K_r^c)) - \tv_\ell^{(K_r^c\setminus A)}(s,\bmu(K_r^c\setminus A)) \Big)  \right)\dd s>0.
 \end{align*}
Consequently, using also \ $\tbv^{(K_r^c\setminus A)}(s,\bmu(K_r^c\setminus A))\leq \tbv^{(K_r^c)}(s,\bmu(K_r^c))$, $s\in[0,t_0]$, we obtain that
  \begin{align*}
  &\int_0^{t_0} \psi^{(K_r^c\setminus A)}\big(\tbv^{(K_r^c)}(s,\bmu(K_r^c))\big)\,\dd s \\
  & > \int_0^{t_0} \big\langle \Bbeta, \tbv^{(K_r^c\setminus A)}(s,\bmu(K_r^c\setminus A)) \big\rangle \,\dd s
      + \int_0^{t_0} \int_{\cU_d\setminus (K_r^c\setminus A)} (1-\ee^{-\langle \tbv^{(K_r^c)}(s,\bmu(K_r^c)) , \br\rangle })\,\nu(\dd\br)\,\dd s\\
  &\geq \int_0^{t_0} \big\langle \Bbeta, \tbv^{(K_r^c\setminus A)}(s,\bmu(K_r^c\setminus A)) \big\rangle \,\dd s
     + \int_0^{t_0} \int_{\cU_d\setminus (K_r^c\setminus A)} (1-\ee^{-\langle \tbv^{(K_r^c\setminus A)}(s,\bmu(K_r^c\setminus A)) , \br\rangle })\,\nu(\dd\br)\,\dd s\\
  & = \int_0^{t_0} \psi^{(K_r^c\setminus A)}\big(\tbv^{(K_r^c\setminus A)}(s,\bmu(K_r^c\setminus A))\big)\,\dd s.
 \end{align*}

Just as in Case I (see \eqref{help38} and \eqref{help41}), we can see that \eqref{help43} implies that \ $\pi_{t_0,\bx}(K_r^c \cap A)>0$.
The aim of the following discussion is to show that \ $\pi_{t,\bx}(K_r^c\cap A)>0$ \ holds as well.
Let $A_r\in\cR_d$ be such that \ $A_r\subseteq K_r\setminus\{\bzero\}\subseteq A$.
By the Markov property of \ $(\bX_s)_{s\in\RR_+}$, \ we have
 \begin{align}\label{help46}
  \begin{split}
   &\pi_{t,\bX_0}(K_r^c\cap A) = \PP\Big( \sup_{s\in(0,t]}  \Delta\bX_s \in K_r^c\cap A  \mid \bX_0\Big)\\
   &\qquad \geq\PP\Big( \Big\{ \sup_{s\in(0,t_0]}  \Delta\bX_s \in K_r^c\cap A  \Big\}
              \cap \Big\{ \sup_{s\in(t_0,t]}  \Delta\bX_s \in A_r\cup\{\bzero\}  \Big\} \;\Big\vert\; \bX_0\Big) \\
   &\qquad=\EE\Big( \bbone_{\{\sup_{s\in(0,t_0]}  \Delta\bX_s \in K_r^c\cap A  \}}
            \PP\Big(\sup_{s\in(t_0,t]}  \Delta\bX_s \in A_r\cup\{\bzero\} \;\big\vert\;  \sigma(\bX_u, u\in[0,t_0])\Big)  \;\Big\vert\;  \bX_0\Big)\\
   &\qquad=\EE\Big( \bbone_{\{\sup_{s\in(0,t_0]}  \Delta\bX_s \in K_r^c\cap A  \}}
            \PP\Big(\sup_{s\in(t_0,t]}  \Delta\bX_s \in A_r\cup\{\bzero\}  \;\big\vert\;  \bX_{t_0}\Big)   \;\Big\vert\; \bX_0\Big).
  \end{split}
 \end{align}
Recall that the conditional distribution of $(\bX_s)_{s\in[t_0,t]}$ given $\bX_{t_0}=\by$ (where $\by\in\RR_+^d$) coincides with that of
 $(\bX_s)_{s\in[0,t-t_0]}$ given $\bX_0=\by$ \ as a consequence of the time-homogeneous Markov property of $(\bX_s)_{s\in\RR_+}$.
Consequently, we obtain that
 \begin{align}\label{help59}
   \PP\Big(\sup_{s\in(t_0,t]}  \Delta\bX_s \in A_r\cup\{\bzero\}  \;\big\vert\;  \bX_{t_0}=\by\Big)
      = \PP_\by\Big(\sup_{s\in(0,t-t_0]}  \Delta\bX_s \in A_r\cup\{\bzero\} \Big),
      \quad \by\in\RR_+^d.
 \end{align}
Using that \ $A_r$ \ is of the form $\left(\prod_{i=1}^d [0, w_i]\right)\setminus\{\bzero\}$
 \ with some $w_1,\ldots,w_d\in\RR_{++}$, $i=1,\ldots,d$, \ we have
 \begin{align*}
   &\Big\{ \sup_{s\in(0,t-t_0]}  \Delta\bX_s \in A_r\cup\{\bzero\}  \Big\}
      \cap\{J_{t-t_0+1}(A_r^c\cap \cU_d)<\infty \}\\
   &\qquad  = \big\{ \tau_{A_r^c\cap\, \cU_d} > t-t_0 \big\} \cap\{J_{t-t_0+1}(A_r^c\cap \cU_d)<\infty \}.
 \end{align*}
By Remark \ref{Rem_A_condition}, we have \ $\Big(\nu+\sum_{i=1}^d \mu_i\Big)(A_r^c\cap \cU_d)<\infty$, \ and hence
 part (i) of Theorem \ref{Thm_jump_times_dCBI} yields that \ $\PP(J_{t-t_0+1}(A_r^c\cap \cU_d)<\infty )=1$.
This together with \eqref{help59} implies that
 \begin{align}\label{help47}
  \begin{split}
    \PP\Big(\sup_{s\in(t_0,t]}  \Delta\bX_s \in A_r\cup\{\bzero\}  \;\big\vert\;  \bX_{t_0}=\by\Big)
      = \PP_\by\Big( \tau_{A_r^c\cap\,\cU_d} > t-t_0\Big),
      \qquad \by\in\RR_+^d.
   \end{split}
 \end{align}
Using again $\left(\nu+\sum_{i=1}^d \mu_i\right)(A_r^c\cap \cU_d)<\infty$,
 by part (iii) of Theorem \ref{Thm_jump_times_dCBI}, we get \ $\PP_\by(\tau_{A_r^c\cap \cU_d} > t-t_0)>0$, \ $\by\in\RR_+^d$.
Hence, by \eqref{help47}, we obtain that
 \begin{align}\label{help48}
  \PP\Big(\sup_{s\in(t_0,t]}  \Delta\bX_s \in A_r\cup\{\bzero\}  \;\big\vert\;  \bX_{t_0}\Big) > 0 \qquad \text{$\PP$-almost surely.}
 \end{align}
Consequently, taking into account that \ $\pi_{t_0,\bx}(K_r^c\cap A)= \PP_\bx(\sup_{s\in(0,t_0]}  \Delta\bX_s \in K_r^c\cap A)>0$, \
 we check that $\pi_{t,\bx}(K_r^c\cap A)>0$.
On the contrary, let us suppose that $\pi_{t,\bx}(K_r^c\cap A)=0$.
 Then, by \eqref{help46}, we would have that
 \[
   \bbone_{\{\sup_{s\in(0,t_0]}  \Delta\bX_s \in K_r^c\cap A  \}}
            \PP_\bx\Big(\sup_{s\in(t_0,t]}  \Delta\bX_s \in A_r\cup\{\bzero\}  \;\big\vert\;  \bX_{t_0}\Big) =0
      \qquad \text{$\PP_\bx$-almost surely.}
 \]
Then, taking into account \eqref{help48}, we would get \ $\PP_\bx(\sup_{s\in(0,t_0]}  \Delta\bX_s \in K_r^c\cap A)=0$,
\ leading us to a contradiction with the fact that $\PP_\bx(\sup_{s\in(0,t_0]}  \Delta\bX_s \in K_r^c\cap A)>0$,
 as desired.

All in all, we derived that $\pi_{t,\bx}(K_r^c\cap A)>0$.
This leads us to the desired contradiction in Case II, since the assumption  \ $\pi_{t,\bx}(A)=0$ \ yields that $\pi_{t,\bx}(K_r^c\cap A)=0$.
\proofend

The next remark is devoted to a consequence of Theorem \ref{Thm_sup_jump_equiv}.
We check that if the total L\'evy measure \ $\nu+\sum_{i=1}^d\mu_i$ \ is not zero, then, for all $t>0$ and $\bx\in\RR_{++}^d$,
 we have \ $\sup_{s\in(0,t]}\Delta\bX_s$ \ can not be $\bzero$ with probability $1$ given that $\bX_0=\bx$.

\begin{Rem}
If the total L\'evy measure \ $\nu+\sum_{i=1}^d\mu_i$ \ is not zero (i.e., \ $\nu\ne0$ \ or \ $\mu_i\ne 0$ \ for some \ $i\in\{1,\ldots,d\}$),
 \ then there exists an \ $A\in\cR_d$ \ such that $\Big(\nu+\sum_{i=1}^d\mu_i\Big)(A)>0$, and, in this case,
 part (ii) of Theorem \ref{Thm_sup_jump_equiv} yields that \ $\pi_{t,\bx}(A)>0$ \ for all $t>0$ and $\bx\in\RR_{++}^d$.
\ Using that $A\subseteq \cU_d$, \ for all \ $t>0$ \ and \ $\bx\in\RR_{++}^d$, \ we get
 \begin{align*}
  0 < \pi_{t,\bx}(A)\leq \PP_\bx\Big(\sup_{s\in(0,t]}\Delta\bX_s \in\cU_d\Big)
    = 1 - \PP_\bx\Big(\sup_{s\in(0,t]}\Delta\bX_s = \bzero \Big),
 \end{align*}
 which yields that
 \ $\PP_\bx\big(\sup_{s\in(0,t]}\Delta\bX_s = \bzero \big) < 1$, \ $t>0$, \ $\bx\in\RR_{++}^d$.
\proofend
\end{Rem}

In the next remark, we point out the fact that in case of $d\geq 2$,
 part (i) of Theorem \ref{Thm_sup_jump_equiv} does not hold for a general Borel set $A\in\cB(\cU_d)$.
This also shows that, in general, the total L\'evy measure of a multi-type CBI process $(\bX_t)_{t\in\RR_+}$
 is not equivalent to the probability measure $\pi_{t,\bx}$, where $\bx\in\RR_+^d$ and $t>0$.

\begin{Rem}\label{Rem_HeLI_comparison}
As a consequence of Theorem 4.2 in He and Li \cite{HeLi}, in case of $d=1$
 both statements (i) and (ii) of our Theorem \ref{Thm_sup_jump_equiv}
 remain true for a general \ $A\in \cB(\cU_1) = \cB((0,\infty))$.
In case of $d\geq 2$, the following counterexample shows that
 part (i) of Theorem \ref{Thm_sup_jump_equiv} does not hold for a general \ $A\in \cB(\cU_d)$.
Let \ $d:=2$ \ and \ $(\bX_t)_{t\in\RR_+}$ \ be a two-type CBI process with parameters
 \ $(2, \bc, \Bbeta, \bB, \nu, \bmu)$ \ such that \ $\EE(\|\bX_0\|) < \infty$, \ $\Bbeta\in\RR_{++}^2$,
  \ $\nu$ \ is the point mass at $(\frac{3}{4}, \frac{1}{4})$,
 \ $\mu_1$ is the point mass at $(\frac{1}{4}, \frac{3}{4})$ \ and \ $\mu_2:=0$.
Then the moment condition \eqref{moment_condition_m_new} holds.
We also have \ $(\nu+\mu_1+\mu_2)(B) =0$ \ for any $B\in\cB(\cU_2)$ such that $(\frac{3}{4}, \frac{1}{4}), (\frac{1}{4}, \frac{3}{4})\notin B$,
 and in this case, by part (ii) of Corollary \ref{Cor_tau}, we have
 $\PP_{\bx}(J_t(B)=0)=1$, $\bx\in\RR_+^2$, $t>0$.
Consequently, we get
 \begin{align}\label{help49}
    \PP_\bx\left(\Delta\bX_u \in\left\{ (0,0), \Big(\frac{1}{4}, \frac{3}{4}\Big), \Big(\frac{3}{4}, \frac{1}{4}\Big) \right\}, \;\; u\in(0,t]\right)=1,
     \qquad \bx\in\RR_+^2, \quad t>0.
 \end{align}
Let $A:=[\frac{1}{2},1]\times[\frac{1}{2},1]$.
Then $A\in \cB(\cU_2)$, but $A$ does not belong to $\cR_2$.
Further, $(\nu+\mu_1+\mu_2)(A)=0$, and we show that it cannot hold that $\pi_{t,\bx}(A)=0$ for all $t>0$ and $\bx\in\RR_+^2$,
 yielding that part (i) of Theorem \ref{Thm_sup_jump_equiv} does not hold for the given 2-type CBI process and $A\in\cB(\cU_2)$.
Let $B:=([0,\frac{1}{2}]\times[0,1])\setminus\{\bzero\}$ and $C:=([0,1]\times[0,\frac{1}{2}])\setminus\{\bzero\}$.
Then $B,C\in\cR_2$, \ $(\nu+\mu_1+\mu_2)(B)=1$ \ and \ $(\nu+\mu_1+\mu_2)(C)=1$, \ and hence, using that $\Bbeta\in\RR_{++}^2$,
 part (ii) of Theorem \ref{Thm_sup_jump_equiv} implies that \ $\pi_{t,\bx}(B)>0$ \ and \ $\pi_{t,\bx}(C)>0$ \ for all $t>0$ and $\bx\in\RR_+^2$.
Consequently, using \eqref{help49}, we have
 \begin{align}\label{help_counterexample}
   0<\pi_{t,\bzero}(B) = \PP_\bzero\left(\sup_{s\in(0,t]} \Delta\bX_s \in B  \right)
       \leq \PP_\bzero\big(\tau_{\{(\frac{1}{4}, \frac{3}{4})\}} \leq t \big), \qquad t>0,
 \end{align}
 and
 \[
   0<\pi_{t,\bzero}(C) = \PP_\bzero\left(\sup_{s\in(0,t]} \Delta\bX_s \in C \right)
      \leq \PP_\bzero\big(\tau_{\{(\frac{3}{4}, \frac{1}{4})\}} \leq t \big), \qquad t>0.
 \]
Using again \eqref{help49}, we have
  \begin{align}\label{help_counterexample2}
   \begin{split}
  \PP_\bzero\Big( \tau_{\{(\frac{1}{4}, \frac{3}{4})\}} \leq t,\; \tau_{\{(\frac{3}{4}, \frac{1}{4})\}} \leq s \Big)
       &\leq \PP_\bzero\left( \sup_{u\in(0,s]} \Delta \bX_u  = \left(\frac{3}{4}, \frac{3}{4}\right) \right)\\
       & \leq \PP_\bzero\left( \sup_{u\in(0,s]} \Delta \bX_u  \in A \right)
         =  \pi_{s,\bzero}(A)
   \end{split}
 \end{align}
 for all  $0<t<s$.
If part (i) of Theorem \ref{Thm_sup_jump_equiv} were true for the given two-type CBI process \ $(\bX_t)_{t\in\RR_+}$ \ and
 \ $A=[\frac{1}{2},1]\times[\frac{1}{2},1]\in\cB(\cU_2)$, then,
 since $(\nu+\mu_1+\mu_2)(A)=0$, we would get that $\pi_{s,\bx}(A)=0$ for all $s>0$ and $\bx\in\RR_+^2$.
In particular, $\pi_{s,0}(A)=0$ were true for all $s>0$, and hence \eqref{help_counterexample2} would imply that
  \begin{align}\label{help_counterexample3}
   \PP_\bzero\left(\tau_{\{(\frac{1}{4}, \frac{3}{4})\}} \leq t, \, \tau_{\{(\frac{3}{4}, \frac{1}{4})\}} \leq s \right) =0,
         \qquad 0<t<s.
 \end{align}
Using the continuity of probability
 and that \ $\bigcup_{s>0}\{\tau_{(\frac{3}{4}, \frac{1}{4})} \leq s\} = \{ \tau_{(\frac{3}{4}, \frac{1}{4})} <\infty \}$, \
  by taking the limit in \eqref{help_counterexample3} as $s\uparrow \infty$ for any fixed $t>0$, we would get
 \begin{align}\label{help_counterexample4}
   \PP_\bzero\left(\tau_{\{(\frac{1}{4}, \frac{3}{4})\}} \leq t, \; \tau_{\{(\frac{3}{4}, \frac{1}{4})\}}
        <\infty \right)  =0,
        \qquad t>0.
 \end{align}
Since \ $\nu(\{(\frac{3}{4}, \frac{1}{4})\})=1\in(0,\infty]$, \ part (iii) of Corollary \ref{Cor_tau} implies
 that \ $\PP_\bzero( \tau_{\{(\frac{3}{4}, \frac{1}{4})\}} < \infty)=1$.
Hence, taking into account \eqref{help_counterexample4}, we have
 \begin{align*}
   \PP_\bzero\left(\tau_{\{(\frac{1}{4}, \frac{3}{4})\}} \leq t \right) =0, \qquad t>0.
 \end{align*}
This leads us to a contradiction (see \eqref{help_counterexample}).
\proofend
\end{Rem}

In the next remark, we highlight why the case \ $\bx=\Bbeta=\bzero$ \ is excluded in  part (ii) of Theorem \ref{Thm_sup_jump_equiv}.

\begin{Rem}\label{Rem_nulla}
Let \ $(\bX_t)_{t\in\RR_+}$ \ be a multi-type CBI process with parameters \ $(d, \bc, \Bbeta, \bB, \nu, \bmu)$
 \ such that \ $\bX_0:=\bzero$, \ $\Bbeta:=\bzero$, \ $\nu$ is the point mass at $\br_0$ with $\br_0:=(d+1,0,\ldots,0)\in\cU_d$,
 \ $\mu_1$ is the point mass at $\bz_0$ \ with $\bz_0:=(1,0,\ldots,0)\in\cU_d$ \ and \ $\mu_i:=0$, $i=2,\ldots,d$.
We readily have that \ $\EE(\|\bX_0\|) < \infty$ \ and the moment condition \eqref{moment_condition_m_new} hold.
Then $\big(\nu+\sum_{i=1}^d\mu_i\big)(B)=0$ for any $B\in\cB(\cU_d)$ such that $\br_0,\bz_0\notin B$, and in this case,
 by part (ii) of Corollary \ref{Cor_tau}, we have $\PP_{\bzero}(J_t(B)=0)=1$, $t>0$.
Consequently, we get
 \begin{align}\label{help50}
    \PP_\bzero\big(\Delta\bX_u \in \{ \bzero, \br_0,\bz_0 \}, \;\; u\in (0,t]\big)=1,\qquad t>0.
 \end{align}
Further, for all $\blambda=(\lambda_1,\ldots,\lambda_d)\in\RR_+^d$, we get that
 \begin{align*}
   &\varphi_1(\blambda) = c_1\lambda_1^2 - \langle \bB \be_1^{(d)}, \blambda \rangle
                          + \ee^{-\lambda_1} - 1
                          + \lambda_1,\\
   &\varphi_i(\blambda) = c_i\lambda_i^2 - \langle \bB \be_i^{(d)}, \blambda \rangle, \qquad i=2,\ldots,d,\\
   &\psi(\blambda) = 1 - \ee^{- (d+1)\lambda_1 },
 \end{align*}
 and, using formulae \eqref{formula_psi_A}  and \eqref{help54},
\begin{align*}
   &\varphi_1^{(\{\bz_0\})}(\blambda)
                         = c_1\lambda_1^2 - \langle \bB \be_1^{(d)}, \blambda \rangle
                              + \lambda_1, \\
   &\varphi_i^{(\{\bz_0\})}(\blambda) = c_i\lambda_i^2 - \langle \bB \be_i^{(d)}, \blambda \rangle, \qquad i=2,\ldots,d,\\
   &\psi^{(\{\bz_0\})}(\blambda) = 1 - \ee^{- (d+1)\lambda_1 }  =  \psi(\blambda).
 \end{align*}
In what follows, let \ $t\in\RR_{++}$ \ be fixed arbitrarily.
Next, we show that if \ $\tau_{\{\bz_0\}}\leq t$, \ then \ $\tau_{\{\br_0\}}\leq t$.
For this it is enough to check that the first jump time of the given CBI process $(\bX_t)_{t\in\RR_+}$ is $\tau_{\{\br_0\}}$,
 \ i.e., $\tau_{\{\br_0\}} \leq \tau_{\{\bz_0\}}$.
First, note that a multi-type CB process starting from \ $\bzero$ \ is the identically \ $\bzero$ \ process.
Further, using that \ $\bX_0=\bzero$ \ and \ $\Bbeta=\bzero$, \ we have \ $(\bX_s)_{s\in[0,\tau_{\{\br_0\}})}$ \ coincides with
 a multi-type CB process starting from $\bzero$ on the time interval \ $[0,\tau_{\{\br_0\}})$.
Indeed, pathwise uniqueness holds for the SDE \eqref{SDE_atirasa_dimd} and both branching processes in question
 are strong solutions of this SDE on the time interval \ $[0,\tau_{\{\br_0\}})$.
This is due to \ $\int_0^s \int_{\cU_d} \br \, M(\dd u,\dd \br) = \bzero$, \
 $0<s<\tau_{\{\br_0\}}$ \ almost surely, where \ $M$ \ is the Poisson random measure on
 \ $\cU_1\times\cU_d$ \ with intensity measure \ $\dd u\,\nu(\dd\br)$ \ appearing in the SDE \eqref{SDE_atirasa_dimd} for $(\bX_s)_{s\in\RR_+}$.
This can be seen as a consequence of the definition of the stochastic integral with respect to $M$,
 more precisely, of the fact that $D_{p_M}\cap (0,s]=\emptyset$, $0<s<\tau_{\{\br_0\}}$ \ almost surely,
 for more details (in particular, for the notation $D_{p_M}$), see, e.g., Step 2 in the proof of Theorem \ref{Thm_jump_times_dCBI}.
To summarize, we get
 \[
   \PP_\bzero\big( \tau_{\{\br_0\}} \leq t \mid  \tau_{\{\bz_0\}}\leq t \big)=1.
 \]
Using \eqref{help50} and that \ $\bz_0\leq \br_0$, \ it implies that
 \ $\PP_\bzero\big( \sup_{s\in(0,t]}\Delta \bX_s \in\{ \bzero,\br_0\} \big)=1$.

Let $A:=([0,1]\times\cdots\times[0,1])\setminus\{\bzero\}\in\cR_d$.
Since $\bzero$ and $\br_0$ \ are not contained in $A$, we get
 \begin{align*}
     \pi_{t,\bzero}(A)
       = \PP_\bzero\Big( \sup_{s\in(0,t]}\Delta \bX_s\in A \Big) =0.
 \end{align*}
However, \ $\big(\nu+\sum_{i=1}^d\mu_i\big)(A) = \mu_1(A) = 1$, \ and hence we obtain that
 part (ii) of Theorem \ref{Thm_sup_jump_equiv} does not hold for the given multi-type CBI process $(\bX_s)_{s\in\RR_+}$ and $A$.
It is not a contradiction, since in part (ii) of Theorem \ref{Thm_sup_jump_equiv} we excluded the case that $\bX_0=\bzero$ and $\Bbeta=\bzero$.
\proofend
\end{Rem}

\section*{Acknowledgements}
We would like to thank Professor Zenghu Li for his constant help in discussing our questions about distributional properties
 of single-type CBI processes, and for sending us a preliminary draft version of the new second edition of his book \cite{Li}.
We would like to thank the referee for his/her comments that helped us to improve the paper.


\addcontentsline{toc}{section}{References}

\begin{thebibliography}{99}

\bibitem{App}
\textsc{Applebaum, D.} (2009).
\textit{L\'evy Processes and Stochastic Calculus, 2nd ed.}
Cambridge University Press, Cambridge.

\bibitem{BarLiPap2}
\textsc{Barczy, M., Li, Z.} and \textsc{Pap, G.} (2015).
Stochastic differential equation with jumps for multi-type continuous state
 and continuous time branching processes with immigration.
\textit{ALEA. Latin American Journal of Probability and Mathematical Statistics}
\textbf{12(1)} 129--169.

\bibitem{BarPalPap}
\textsc{Barczy, M., Palau, S.} and \textsc{Pap, G.} (2020).
Almost sure, $L_1$- and $L_2$-growth behavior of supercritical multi-type
 continuous state and continuous time branching processes with immigration.
\textit{Science China Mathematics}
\textbf{63(10)} 2089--2116.

\bibitem{BarPalPap2}
\textsc{Barczy, M., Palau, S.} and \textsc{Pap, G.}  (2021).
Asymptotic behavior of projections of supercritical multi-type continuous-state and continuous-time branching processes with immigration.
\textit{Advances in Applied Probability}
\textbf{53(4)} 1023--1060.

\bibitem{BarPap}
\textsc{Barczy, M.} and \textsc{Pap, G.} (2016).
Asymptotic behavior of critical, irreducible multi-type continuous state and
 continuous time branching processes with immigration.
\textit{Stochastics and Dynamics}
\textbf{16(4)} Article ID 1650008, 30 pages.

\bibitem{CabPerUri}
\textsc{Caballero, M. E.}, \textsc{P\'erez Garmendia, J. L.}, and \textsc{Uribe Bravo, G.} (2017).
Affine processes on $\mathbb{R}^m_+\times \mathbb{R}^n$ and multiparameter time changes.
\textit{Annales de l'Institut Henri Poincaré - Probabilités et Statistiques}
\textbf{3(53)} 1280--1304.

\bibitem{ChaMar2}
\textsc{Chaumont, L.} and \textsc{Marolleau, M.} (2023).
Extinction times of multitype continuous-state branching processes.
\textit{Annales de l’Institut Henri Poincaré - Probabilités et Statistiques}
\textbf{59(2)} 563--577.

\bibitem{CheLi}
\textsc{Chen, S.} and \textsc{Li, Z.} (2021).
Continuous time mixed state branching processes and stochastic equations.
\textit{Acta Mathematica Scientia}
\textbf{41B(5)} 1445--1473.

\bibitem{DufFilSch}
\textsc{Duffie, D., Filipovi\'{c}, D.} and \textsc{Schachermayer, W.} (2003).
Affine processes and applications in finance.
\textit{Annals of Applied Probability}
\textbf{13} 984--1053.

\bibitem{FilMaySch}
\textsc{Filipovi\'c, D.,  Mayerhofer, E.} and \textsc{Schneider, P. } (2013).
 Density approximations for multivariate affine jump-diffusion processes.
\textit{Journal of Econometrics}
\textbf{176(2)} 93--111.

\bibitem{FriJinRud}
\textsc{Friesen, M., Jin, P.} and \textsc{R\"udiger, B.} (2020).
On the boundary behavior of multi-type continuous-state branching processes with immigration.
\textit{Electronic Communications in Probability}
\textbf{25} 1--14.

\bibitem{FriJinRud2}
\textsc{Friesen, M., Jin, P.} and \textsc{R\"udiger, B.} (2020).
Existence of densities for multi-type continuous-state branching processes with immigration.
\textit{Stochastic Processes and their Applications}
\textbf{130(9)} 5426--5452.


\bibitem{HeLi}
\textsc{He, X.} and \textsc{Li, Z.} (2016).
Distributions of jumps in a continuous-state branching process with immigration.
\textit{Journal of Applied Probability}
\textbf{53} 1166--1177.

\bibitem{HorJoh}
\textsc{Horn, R. A.} and \textsc{Johnson, Ch.\ R.} (2013).
\textit{Matrix Analysis, 2nd ed.}
Cambridge University Press, Cambridge.


\bibitem{HorXu}
\textsc{Horst, U.} and \textsc{Xu, W.} (2022).
The microstructure of stochastic volatility models with self-exciting jump dynamics.
\textit{The Annals of Applied Probability}
\textbf{32(6)} 4568--4610.

\bibitem{IkeWat}
\textsc{Ikeda, N.} and \textsc{Watanabe, S.} (1989).
\textit{Stochastic Differential Equations and Diffusion Processes, 2nd ed.}
 North-Holland Publishing Co., Amsterdam, Kodansha, Ltd., Tokyo.

\bibitem{JiLi}
\textsc{Ji, L.} and \textsc{Li, Z.} (2020).
Moments of continuous-state branching processes with or without immigration.
 \textit{Acta Mathematicae Applicatae Sinica, English Series}
 \textbf{36(2)} 361--373.

\bibitem{JiaMaSco}
\textsc{Jiao, Y., Ma, C.} and \textsc{Scotti, S.} (2017).
Alpha-CIR model with branching processes in sovereign interest rate modeling.
\textit{Finance and Stochastics}
\textbf{21} 789--813.

\bibitem{K2}
\textsc{Kallenberg, O.} (2017).
\textit{Random Measures, Theory and Applications}.
Springer, Cham.

\bibitem{Kel}
\textsc{Keller-Ressel, M.} (2008).
\textit{Affine Processes -- Theory and Applications in Finance (Ph.D. thesis).}
Vienna University of Technology, 110 pages.

\bibitem{Kin}
\textsc{Kingman, J. F. C.} (1993).
\textit{Poisson Processes}.
Oxford University Press, New York.

\bibitem{Kyp}
\textsc{Kyprianou, A. E.} (2014).
\textit{Fluctuations of L\'evy Processes with Applications, 2nd ed.}
 Springer-Verlag, Heidelberg.

\bibitem{KypPal}
\textsc{Kyprianou, A.} and \textsc{Palau, S.} (2018).
Extinction properties of multi-type continuous-state branching processes.
\textit{Stochastic Processes and their Applications} \textbf{128(10)} 3466--3489.

\bibitem{KypPalRen}
\textsc{Kyprianou, A., Palau, S.} and \textsc{Ren, Y. X.} (2018).
Almost sure growth of supercritical multi-type continuous-state branching process.
\textit{ALEA} \textbf{15} 409--428.

\bibitem{Li}
\textsc{Li, Z.} (2022).
\textit{Measure-Valued Branching Markov Processes, 2nd ed. }
  Springer-Verlag GmbH, Berlin.

\bibitem{Li3}
\textsc{Li, Z.} (2020).
Continuous-state branching processes with immigration.
 From probability to finance--lecture notes of BICMR Summer School on Financial Mathematics.
\textit{Mathematical Lectures from Peking University,}
\textbf{1--69.}
Springer, Singapore.

\bibitem{Ma}
\textsc{Ma, R. G.} (2013).
Stochastic equations for two-type continuous-state branching processes with immigration.
\textit{Acta Mathematica Sinica, English Series}
\textbf{29} 287--294.

\bibitem{Wa}
\textsc{Watanabe, S.} (1969).
On two dimensional Markov processes with branching property.
\textit{Transactions of the American Mathematical Society}
\textbf{136} 447--466.



\end{thebibliography}
\end{document}